\documentclass[11pt]{article}

\usepackage{amsfonts}
\usepackage{amsthm}
\usepackage{amsmath}
\usepackage{amssymb}

\usepackage{overpic}
\usepackage{graphics}
\usepackage{subfig}

\usepackage{placeins} 

\usepackage{booktabs} 

\usepackage{verbatim}
\usepackage{xcolor} 
\definecolor{blun}{cmyk}{0.8, 0.5, 0, 0.7}

\usepackage{multicol}
\usepackage{listings}

\usepackage[all]{xy}
\usepackage{makeidx}
\usepackage{enumerate}
\usepackage{enumitem}
\usepackage[letterpaper,margin=0.9in]{geometry}

\usepackage[nottoc]{tocbibind}
\usepackage[titletoc,title]{appendix}
\numberwithin{equation}{section}

\usepackage{hyperref}
\hypersetup{bookmarks,colorlinks,pdfauthor={Cinzia Soresina},linkcolor={blun},citecolor={blun},urlcolor={blun},hyperfootnotes={false}}

\let\oldFootnote\footnote
\newcommand\nextToken\relax

\renewcommand\footnote[1]{%
    \oldFootnote{#1}\futurelet\nextToken\isFootnote}

\newcommand\isFootnote{%
    \ifx\footnote\nextToken\textsuperscript{,}\fi}

\let\oldbibliography\thebibliography
\renewcommand{\thebibliography}[1]{%
  \oldbibliography{#1}%
  \setlength{\itemsep}{-1.2mm}%
}

\theoremstyle{plain}

\theoremstyle{definition}

\theoremstyle{remark}

\newtheoremstyle{myremark}
  {3pt}
  {3pt}
  {\small \rmfamily}
  {5pt}
  {\rmfamily}
  {:}
  {.5em}
  {}

\theoremstyle{myremark}
\newtheorem*{remark}{\textit{Remark}}

\def\R{\mathbb{R}}

\def\txtd{{\textnormal{d}}}

\def\txtD{{\textnormal{D}}}

\def\ra{\rightarrow}

\newcommand{\Lap}{\Delta}
\newcommand{\tr}{\textnormal{tr}}

\newcommand{\be}{\begin{equation}}
\newcommand{\ee}{\end{equation}}
\newcommand{\benn}{\begin{equation*}}
\newcommand{\eenn}{\end{equation*}}
\newcommand{\bea}{\begin{eqnarray}}
\newcommand{\eea}{\end{eqnarray}}
\newcommand{\beann}{\begin{eqnarray*}}
\newcommand{\eeann}{\end{eqnarray*}}

\newcommand{\myendex}{$\blacklozenge$\end{ex}}
\newcommand{\myendexerc}{$\lozenge$\end{exerc}}
\newcommand{\myendpexerc}{$\lozenge$\end{pexerc}}

\begin{document}

\author{Christian Kuehn\footnote{
Zentrum Mathematik, Technische Universit\"at M\"unchen, Boltzmannstr. 3, 85748 Garching bei M\"unchen, Germany} \footnote{email: \href{mailto: ckuehn@ma.tum.de}{ckuehn@ma.tum.de}}, Cinzia Soresina\footnotemark[1] \footnote{email: \href{mailto: soresina@ma.tum.de}{soresina@ma.tum.de}}}

\title{Numerical continuation for a fast-reaction system \\and its cross-diffusion limit}

\maketitle

\begin{abstract}
\noindent
In this paper we investigate the bifurcation structure of the triangular SKT model in the weak competition regime and of the corresponding fast-reaction system in 1D and 2D domains via numerical continuation methods. We show that the software \texttt{pde2path} can be adapted to treat cross-diffusion systems, reproducing the already computed bifurcation diagrams on 1D domains. We show the convergence of the bifurcation structure obtained selecting the growth rate as bifurcation parameter. Then, we compute the bifurcation diagram on a 2D rectangular domain providing the shape of the solutions along the branches and linking the results with the Turing instability analysis. In 1D and 2D, we pay particular attention to the fast-reaction limit by always computing sequences of bifurcation diagrams as the time-scale separation parameter tends to zero. We show that the bifurcation diagram undergoes major deformations once the fast-reaction systems limits onto the cross-diffusion singular limit. Furthermore, we find evidence for time-periodic solutions by detecting Hopf bifurcations, we characterize several regions of multi-stability, and improve our understanding of the shape of patterns in 2D for the SKT model.  
\end{abstract}







\section{Introduction}\label{sec_int}

Systems with multiple time-scales appear in a wide variety of mathematical areas but also in many applications~\cite{KuehnBook}. The basic class of multiple time-scale systems has two time-scales, which are so-called fast--slow systems given by
\begin{equation}
\begin{array}{lcl}
\dfrac{\txtd u}{\txtd t} &=& \dfrac{1}{\varepsilon} f(u,v,\varepsilon),\\[0.2cm]
\dfrac{\txtd v}{\txtd t} &=& g(u,v,\varepsilon),\\
\end{array}
\end{equation}
where $u,v$ are the unknowns and $\varepsilon>0$ is a small parameter so that $u$ is fast while $v$ is slow. There is quite substantial theory for the case of fast--slow ordinary differential equations (ODEs). Although there are many links of the ODE theory to multiscale partial differential equations (PDEs)~\cite{KuehnBook,KuehnBook1}, a lot less is known about fast--slow PDEs. One important sub-class are fast-reaction PDEs given by
\begin{equation}
\label{eq:fsPDE}
\begin{array}{lcl}
\partial_t u  &=& \kappa_u \Delta u + \dfrac1\varepsilon f(u,v,\varepsilon),\\[0.2cm]
\partial_t v &=& \kappa_v \Delta v + g(u,v,\varepsilon),\\
\end{array}
\end{equation}
posed on a domain $\R_+\times \Omega$, where $\Delta$ is the usual spatial Laplacian with respect to $x\in\Omega$, $t\in[0,T)$ for some $T>0$, $u=u(t,x)$, and $v=v(t,x)$. In this work, one of the two key equations we study is a particular version of~\eqref{eq:fsPDE} arising in mathematical ecology \cite{iida2006diffusion}
\begin{equation}
\label{fast}
\begin{array}{lcl}
\partial_t u_1&=&d_1\Lap u_1+(r_1-a_1 u-b_1 v)u_1+
\dfrac{1}{\varepsilon}\left(u_2\left( 1-\dfrac v M \right)-u_1\dfrac v M \right),\\[0.2cm]
\partial_t u_2&=&(d_1+d_{12} M)\Lap u_2+(r_1-a_1 u-b_1 v)u_2-\dfrac{1}{\varepsilon}\left(u_2\left( 1-\dfrac v M \right)-u_1\dfrac v M \right),\\[0.2cm]
\partial_t v&=&d_2\Lap v+(r_2-b_2 (u_1+u_2)-a_2 v)v.\\
\end{array}
\end{equation}
where the quantities $u_{1,2}(t,x),\, v(t,x)\geq 0$ represent the population densities of two species at time~$t$ and position~$x$, confined and competing for resources on a bounded domain $\Omega \subset\mathbb{R}^N$. The positive coefficients $d_i,\,r_i,\,a_i,\,b_i \, (i=1,2)$ describe the diffusion, the intrinsic growth, the intra-specific competition and the inter-specific competition rates respectively. The population $u$ is split into two states, distinguishing between less active and active individuals, respectively denoted by $u_1$ and $u_2$ (thus $u=u_1+u_2$). We assume throughout that $0\leq u_2(t,x)\leq M$ on $[0,T]\times\Omega$ for a fixed constant $M>0$. The coefficient $d_{12}> 0$ stands for competition pressure between the two sub-classes. Each individual in the sub-classes converts its state into the other one depending on the spatial distribution of the competitor $v$. 

Since $\varepsilon>0$ is assumed to be small, it is natural to consider the fast-reaction limit $\varepsilon\rightarrow 0$. This limit aims to model the dynamics of the fast variables as instantaneous and reduce the analysis to the slow time-scale dynamics. Yet, one easily sees that this limit is very singular. To motivate this observation, suppose we discard the diffusion terms in~\eqref{fast} and linearize the fast-reaction terms, then we get a matrix
\begin{equation}
\txtD_u f(u_1,u_2,v,0)\in\R^{2\times 2}
\end{equation}
which always has a zero eigenvalue. This means that the classical normal hyperbolicity condition~\cite{Fenichel4} for fast--slow ODEs fails. Yet, the intuition is that the diffusion terms should help to still obtain a well-defined fast-reaction limit as $\varepsilon\ra 0$. 

In the literature, several results for such fast-reaction limits for various PDEs exist. One of the first works is the paper~\cite{hilhorst1996fast} presenting the fast-reaction limit in a system of one parabolic and one ordinary differential equation. A reaction--diffusion system which models a fast reversible reaction between two mobile reactants was then treated in~\cite{bothe2003reaction} and the limiting problem yields a nonlinear diffusion term. Also fast irreversible reactions (in which two chemical components form a product) were considered, where the limiting system is a Stefan-type problem with a moving interface at which the chemical reaction front is localized \cite{bothe2012instantaneous}. Furthermore, in~\cite{henneke2016fast} a dynamical boundary condition has been interpreted as a fast-reaction limit of a volume--surface reaction--diffusion system. It turns out that when the fast-reaction system has three components, the limiting system are often of two types: free boundary problems \cite{murakawa2011fast} and cross-diffusion systems \cite{bothe2012cross}, which also arise in population dynamics \cite{conforto2018about,desvillettes2018non,iida2006diffusion}. In this framework, individuals of one or more species exist in different states and the small parameter represents the average switching time. In our context, it is known that the limiting system of~\eqref{fast} is a cross-diffusion PDE~\cite{hilhorst2009fast}. The limiting triangular cross-diffusion system~\cite{iida2006diffusion} is given by
\begin{equation}\label{cross}
\begin{array}{ll}
\partial_t u=\Lap((d_1+d_{12} v)u)+(r_1-a_1 u-b_1 v)u,& \textnormal{on } \mathbb{R}_+\times \Omega,\\[0.1cm]
\partial_t v=d_2\Lap v+(r_2-b_2 u-a_2 v)v,& \textnormal{on }  \mathbb{R}_+\times \Omega,\\[0.1cm]
\dfrac{\partial u}{\partial n}=\dfrac{\partial v}{\partial n}=0,&\textnormal{on } \mathbb{R}_+\times \partial\Omega,\\[0.2cm]
u(0,x)=u_{in}(x),\; v(0,x)=v_{in}(x),& \textnormal{on } \Omega,\\
\end{array}
\end{equation}
where the quantities $u(t,x),\, v(t,x)\geq 0$ again represent the population densities of the same two species at time $t$ and position $x$, confined and competing for resources on a bounded domain $\Omega \subset\mathbb{R}^N$. The coefficients as above and are all supposed to be positive constants. The model~\eqref{cross} is known as Shigesada--Kawasaki--Teramoto (SKT) model as it was initially proposed in \cite{shigesada1979spatial} in 1979 to account for stable non-homogeneous steady states in certain ecological systems. These states describe, for suitable parameters sets, \emph{spatial segregation} that is a situation of coexistence of two competitive species on a bounded domain. Historically, the SKT model was first postulated without any reference to a larger system such as~\eqref{fast}. In 2006, the fast-reaction system~\eqref{fast} was introduced in \cite{iida2006diffusion} to approximate bounded solutions to the cross-diffusion system \eqref{cross} in $[0,T]\times\Omega$. 

In this paper, we aim to contribute the understanding of both variants of the SKT model: the three-component fast-reaction PDE as well as its cross-diffusion singular limit. A primary focus from the ecological and mathematical viewpoints is evidently played by the steady states (or equilibria) of both of the models. So we briefly recall the well-known structure of \emph{homogeneous} steady states starting from the cross-diffusion SKT model~\eqref{cross}. 

The homogeneous solutions are the total extinction $(0,0)$ (always unstable), two non-coexistence states $(\bar u,0)=(r_1/a_1,0)$ and $(0,\bar v)=(0,r_2/a_2)$, and one coexistence state (when it exists)
$$(u_*,v_*)=\left( \dfrac{r_1a_2-r_2b_1}{a_1a_2-b_1b_2},\dfrac{r_2a_1-r_1b_2}{a_1a_2-b_1b_2}\right).$$
While the non-coexistence equilibria exist for all the parameter values, the coexistence is admissible (positive coordinates) only in two cases.
\begin{itemize}[leftmargin=0.3cm]
\item[-] \textit{weak competition} or \textit{strong intra-specific competition}, namely $a_1a_2-b_1b_2>0$.\\
Under the additional condition on the growth rates $b_1/a_2<r_1/r_2<a_1/b_2$, without diffusion the coexistence homogeneous steady state is stable, while the non-coexistence ones are unstable. With standard diffusion in a convex domain and with zero-flux boundary conditions, any non-negative solution generically converges to homogeneous one, and this implies that the two species coexist but their densities are homogeneous in the whole domain \cite{kishimoto1985spatial}. With cross-diffusion, the model exhibits stable non-homogeneous steady states if $d_1,\; d_2$ are small enough or $d_{12}$ is large enough \cite{iida2006diffusion}. \\
If $r_1/r_2<b_1/a_2$ or $r_1/r_2>a_1/b_2$, the coexistence homogeneous steady state is not positive and the non-coexistence ones are stable for the reaction part (in absence of diffusion).
\item[-] \textit{strong competition} or \textit{strong inter-specific competition}, namely $a_1a_2-b_1b_2<0$.\\
When $a_1/b_2<r_1/r_2<b_1/a_2$, the coexistence homogeneous steady state is unstable, while the non-coexistence ones are stable. With only standard diffusion, in a convex domain and with zero-flux boundary conditions, it has been shown that if positive and non-constant steady states exist, they must be unstable \cite{kishimoto1985spatial}, and numerical simulations suggest that any non-negative solution generically converges to either $(\bar u,0)$ or $(0,\bar v)$, that is, the competitive exclusion principle occurs between the two species. However, adding the cross-diffusion term, it is numerically shown that even if the domain is convex, there exists stable non-homogeneous solutions exhibiting spatially segregating coexistence when $d_{12}$ is suitable large \cite{izuhara2008reaction}.\\
Analogously to the previous regime, if $r_1/r_2<a_1/b_2$ or $r_1/r_2>b_1/a_2$, the coexistence homogeneous steady state is not positive and the non-coexistence ones are stable for the reaction part (in absence of diffusion).
\end{itemize}
We remark that the stationary problem of the classical Lotka--Volterra competition model with linear diffusion endowed with Dirichlet boundary conditions has been extensively studied (see \cite{blat1984bifurcation, conti2005asymptotic, cosner1984stable} and references therein), and it shows different features to the SKT cross-diffusion model. For the fast-reaction PDE~\eqref{fast}, the homogeneous steady state is 
$$(u_{1*},u_{2*},v_*)=\left(u_*\left(1-\dfrac{v_*}{M}\right),u_*\dfrac{v_*}{M},v_*\right)$$
where $u_*,\;v_*$ have the same expression as in the cross-diffusion case. In particular, the homogeneous steady state turns out to be independent of $\varepsilon$, which is a nice starting point for a comparative study as it remains to investigate also the non-homogeneous steady states.

Before proceeding to summarize the main results of our work, we briefly put it into context with the existing literature. Both systems, the cross-diffusion model and the fast-reaction model, are interesting from different mathematical viewpoints. From the modeling point of view the justification of cross-diffusion terms by means of semilinear reaction--diffusion systems including simple reactions and linear diffusion is fundamental to the understanding of the hidden processes that they can capture. Vice versa, their approximation with simpler models (not in the number of equations but due to the possibility to remove the nonlinearity in the diffusion) is useful both for the analysis and the numerics.

On the one hand, theoretical results require sophisticated techniques. Regarding the cross-diffusion system, existence, smoothness and uniqueness of solutions have been widely investigated (see~\cite{Ama88, Ama90, desvillettes2015entropic, galiano2003semi, jungel2016book} and the references therein). The convergence of the solution of the stationary problem have been treated in~\cite{izuhara2008reaction}, while the convergence of the solutions of the fast-reaction system~\eqref{fast} towards the solutions of the cross-diffusion system~\eqref{cross} was shown in~\cite{conforto2014rigorous} in dimension one, and then generalized to a wider set of admissible reaction terms and in any dimension~\cite{desvillettes2015new}. Similar results were obtained for a class of non-triangular cross-diffusion systems~\cite{murakawa2012relation}, assuming the same time-scale for all the species. Finally, it has been proven that the limiting system inherits the entropy structure with an entropy that is not the sum of the entropies of the components~\cite{daus2017cross}.

On the other hand, the capability of the cross-diffusion system~\eqref{cross} to model the spatial segregation of competing species is related to the appearance of non-homogeneous solutions. Although the system~\eqref{cross} does not have an activator--inhibitor structure (the most important mechanism in the Turing instability theory for pattern formation), the cross-diffusion term turns out to be the key ingredient to destabilize the homogeneous equilibrium~\cite{gambino2012turing, iida2006diffusion, zou2017bifurcation}. In this framework bifurcations diagrams are useful to present and explore the bifurcating branches and the non-homogeneous steady state solutions. Numerical continuation techniques are going to allow us to obtain a more global picture far from the homogeneous branch. The bifurcation diagram of the triangular cross-diffusion system on a 1D domain was presented in~\cite{iida2006diffusion}; the existence of these non-homogeneous steady states significantly far from being perturbations of the homogeneous solutions, was proved in~\cite{breden2018existence} applying a computer assisted method. A mathematically rigorous construction of the bifurcation structure of the three-component system was obtained in~\cite{breden2013global}, while the convergence of the bifurcation structure on a 1D domain was also investigated in~\cite{izuhara2008reaction} with respect to two different bifurcation parameters appearing in the diffusion part, both in the weak and strong competition case. Only recently, the influence of the additional cross-diffusion term in the system has been investigated, combining a detailed linearized analysis and numerical continuation. To the best of the authors' knowledge, not much has been done on 2D domains: the possible pattern admitted by the non-triangular SKT model was explored in~\cite{gambino2013pattern}, but the bifurcation structure on 2D domain is not known.

To study the bifurcation structure, we use and extend the continuation software for PDEs \texttt{pde2path}~\cite{dohnal2014pde2path, uecker2018hopf, uecker2014pde2path}, based on a FEM discretization of the stationary problem. Since the general class of problems which can be numerically analyzed by the software does not include the cross-diffusion term appearing in~\eqref{cross}, it requires an additional setup to be able to compare the fast-reaction PDE results with the singular limit cross-diffusion systems. Yet, it has been shown in the past for other classes of PDEs related to classical elliptic problems that \texttt{pde2path} can be extended beyond its standard setting~\cite{KuehnEllipticCont}. This is the key reason, why we have chosen to carry out numerical continuation within this framework.

\medskip
Here is a summary of the main contributions of this paper.
\begin{itemize}[leftmargin=0.3cm]
\item[-] We explain how to set up the continuation software \texttt{pde2path} in order to treat cross-diffusion terms. Then we cross-validate and extend previous computations for various regimes of the time-scale separation parameter for 1D domains.
\item[-] We compute a new bifurcation diagram with respect to a parameter appearing in the reaction part. We show how the bifurcation structure of the fast-reaction system modifies when $\varepsilon \to 0$. In particular, a novel ``broken-heart'' structure of non-homogeneous steady state bifurcation branches is observed in the singular limit.
\item[-] Then we compute and interpret the various bifurcation diagrams and the non-homogeneous solutions of the triangular SKT model on a 2D rectangular domain. Also in this case the convergence towards the singular limit is analyzed carefully.
\item[-] We show a link between the computed bifurcation diagrams in 1D and 2D domains, and the Turing instability analysis as a tool to fully understand and validate the numerical continuation calculations. 
\end{itemize}

\medskip
The paper is organized as follows. In Section~\ref{sec:nc1D} we present the the numerical continuation results obtained with \texttt{pde2path} in the 1D case: we provide a detailed picture of the different types of stable unstable non-homogeneous solutions to the cross-diffusion system arising at each bifurcation point, and we quantify the convergence of the fast--slow system to the cross-diffusion one. Section~\ref{sec:nc2D} is devoted to the 2D case: we show the bifurcation diagram and different patterns, as well as the bifurcation diagrams of the fast-reaction system. In Section \ref{sec:concl} some concluding remarks can be found. Appendix \ref{app:setup} contains the \texttt{pde2path} setup for cross-diffusion systems, while in Appendix \ref{app:TI} we report the Turing instability analysis of the cross-diffusion system for reference and validation.

\FloatBarrier
\section{Numerical continuation on a 1D domain}
\label{sec:nc1D}

The numerical analysis of systems \eqref{cross} and \eqref{fast} is performed using the continuation software \texttt{pde2path}, originally developed to treat standard reaction--diffusion systems and here adapted to investigate cross-diffusion systems. The software setup required for system \eqref{cross} can be found in Appendix \ref{app:setup}. For the numerical results we use the set of parameter values widely used in literature for the weak competition regime \cite{breden2018existence, breden2013global,iida2006diffusion, izuhara2008reaction}, here reported in Table \ref{tab:param}. We set $d_1=d_2=:d$ and use $d$ as one main bifurcation parameter.
It follows that
$$u_{1*}=\dfrac{91}{64},\;u_{2*}=\dfrac{13}{64},\quad u_*=\dfrac{13}{8},\quad v_*=\dfrac{1}{8},\qquad  \dfrac{2}{3}<\dfrac{r_1}{r_2}<6,$$
which is one possible starting point for the continuation.

\begin{table}
\centering
\begin{tabular}{cccccc|cc}
\toprule
$r_1$&$r_2$&$a_1$&$a_2$&$b_1$&$b_2$&$d_{12}$&$M$\\
\midrule
5&2&3&3&1&1&3&1\\
\bottomrule
\end{tabular}
\caption{Parameter values used in numerical continuation. The set $r_i,\,a_i,\,b_i,\, (i=1,2)$ corresponds to the weak competition case (or strong intra-specific competition).}
\label{tab:param}
\end{table}

In this section we consider a 1D domain (interval) $\Omega=(0,1)$, as in \cite{breden2018existence, MBCKCS, breden2013global, iida2006diffusion, izuhara2008reaction}. We provide in Section \ref{subsec:1Dsol} a detailed characterization of different steady state types of the cross-diffusion system \eqref{cross}. In Section \ref{subsec:1D_conv_d} we study the convergence of the bifurcation structure of the fast-reaction system \eqref{fast} to the cross-diffusion system \eqref{cross}, taking the standard diffusion coefficient $d$ as bifurcation parameter. In this framework, reproducing such diagrams is crucial to test the numerical continuation software; we add to the existing literature a quantification of the convergence of the bifurcation structure when $\varepsilon\to 0$, a clear picture of the behavior of the non-homogeneous solutions along the bifurcation branches and their stability properties, and the corresponding bifurcation diagram in the $L_1$-norm. Finally, in Section \ref{subsec:conv_r1} we consider the growth rate $r_1$ as bifurcation parameter: also in this case we present the bifurcation diagram of the cross-diffusion system and the stable steady states appearing beyond the usual range of parameters, as well as the convergence of the bifurcation structure of the fast-reaction system to the cross-diffusion one.

\subsection{Bifurcation diagram of the cross-diffusion system}
\label{subsec:1Dsol}

We numerically compute the bifurcation diagram of the cross-diffusion system \eqref{cross} using the continuation software \texttt{pde2path} and we provide a clear picture of the different solution types.

\begin{figure}
\subfloat[\label{bifdiagv0}]{
       \begin{overpic}[width=0.5\textwidth,tics=10,trim=90 240 100 250,clip]{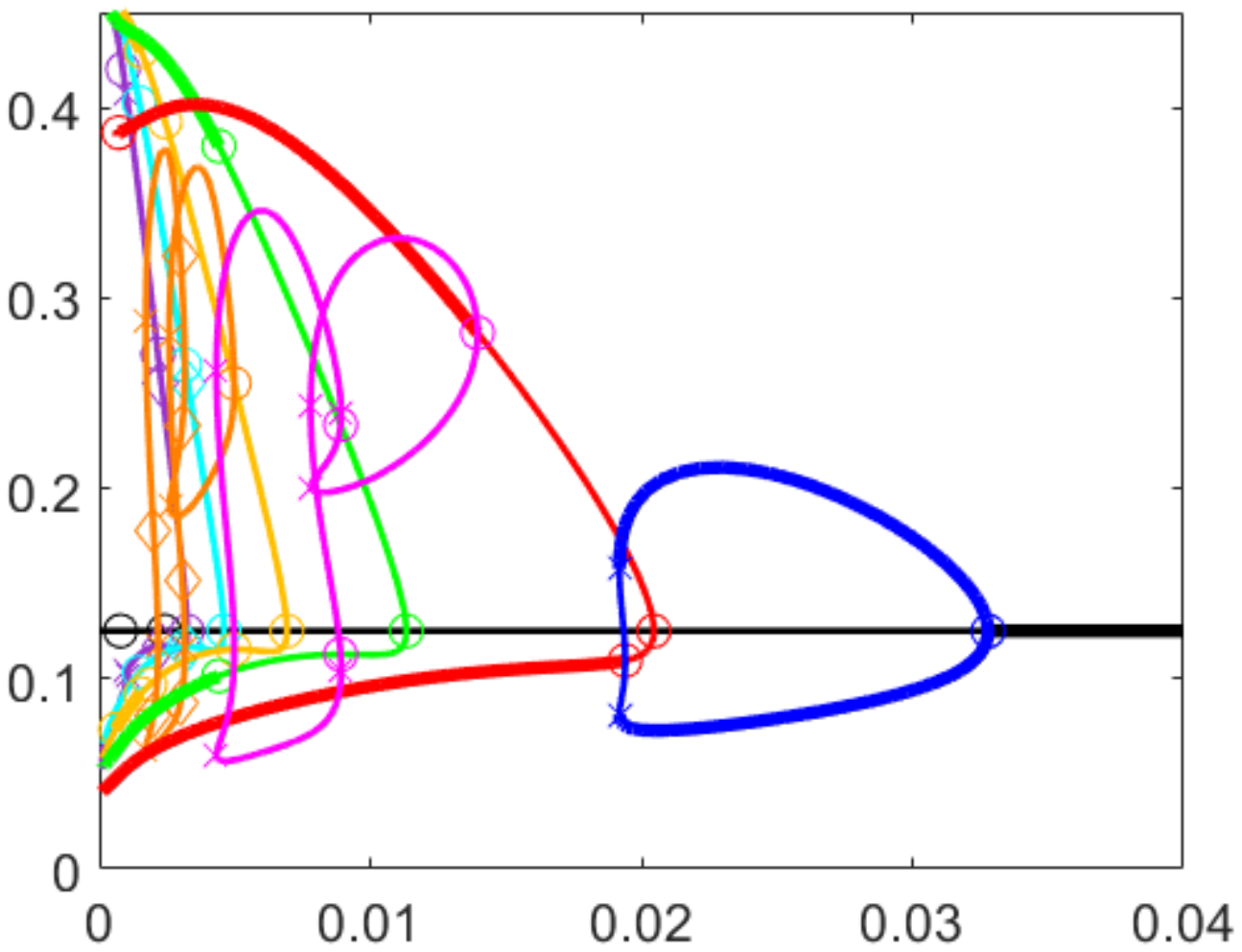}
         \put(0,30){\rotatebox{90}{$v(0)$}}
         \put(80,-2){$d$}
         \put(79,27){$\mathtt{B_1}$}
         \put(56,27){$\mathtt{B_2}$}
         \put(37,27){$\mathtt{B_3}$}
       \end{overpic}}
\subfloat[\label{bifdiag_L2}]{
       \begin{overpic}[width=0.5\textwidth,tics=10,trim=90 240 100 250,clip]{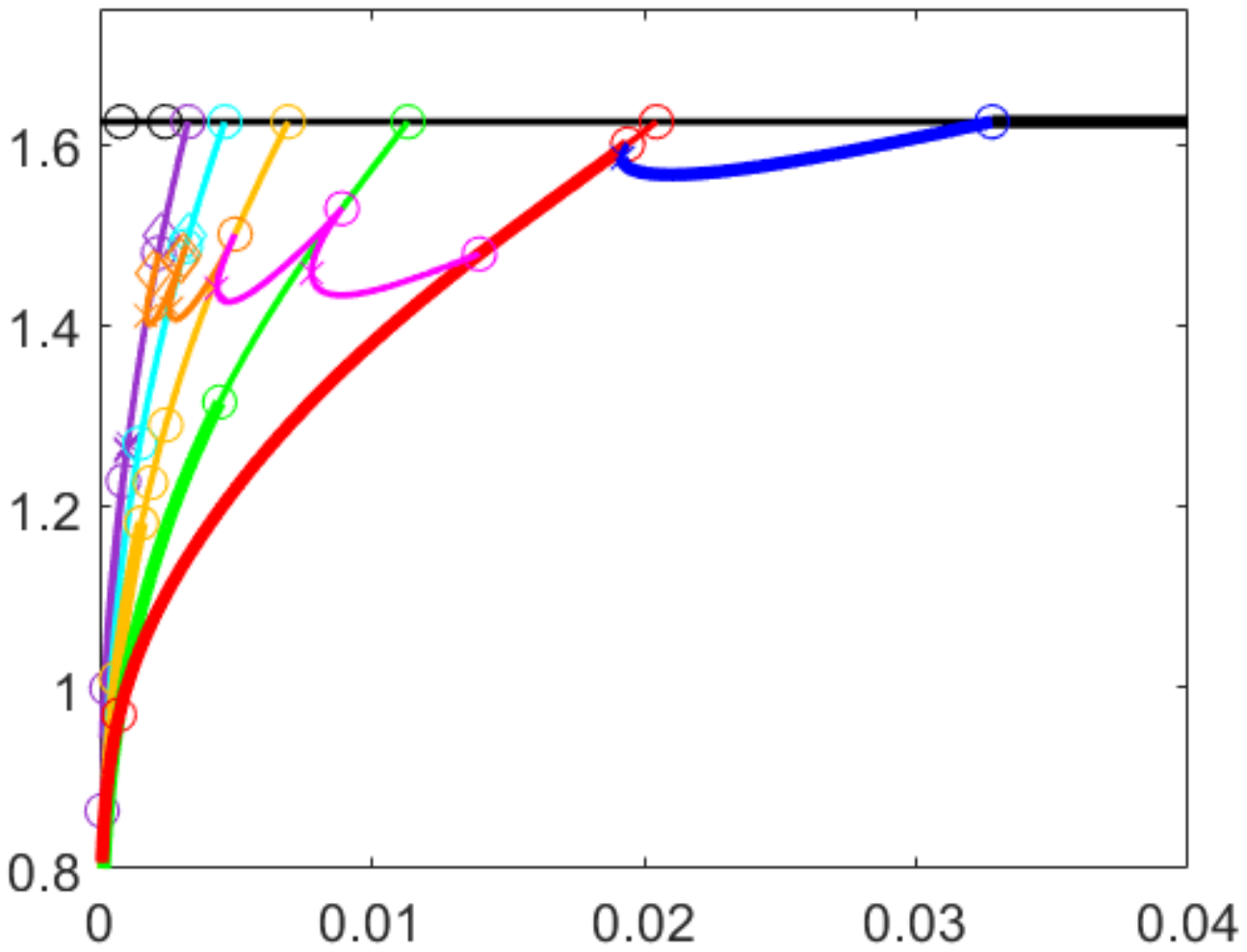}
         \put(0,30){\rotatebox{90}{$||u||_{L^1}$}}
         \put(80,-2){$d$}
         \put(77,63){$\mathtt{B_1}$}
         \put(54,63){$\mathtt{B_2}$}
         \put(37,63){$\mathtt{B_3}$}
       \end{overpic}}
\caption{Bifurcation diagram of the cross-diffusion system with respect to different quantities on the vertical axis: \protect\subref{bifdiagv0} $v(0)$, \protect\subref{bifdiag_L2} $||u||_{L_1}$. }
\label{DiagBifCross_y}
\end{figure}

\begin{figure}
\subfloat[blue branch\label{solb}]{
       \begin{overpic}[width=0.3\textwidth,tics=10,trim=90 240 100 250,clip]{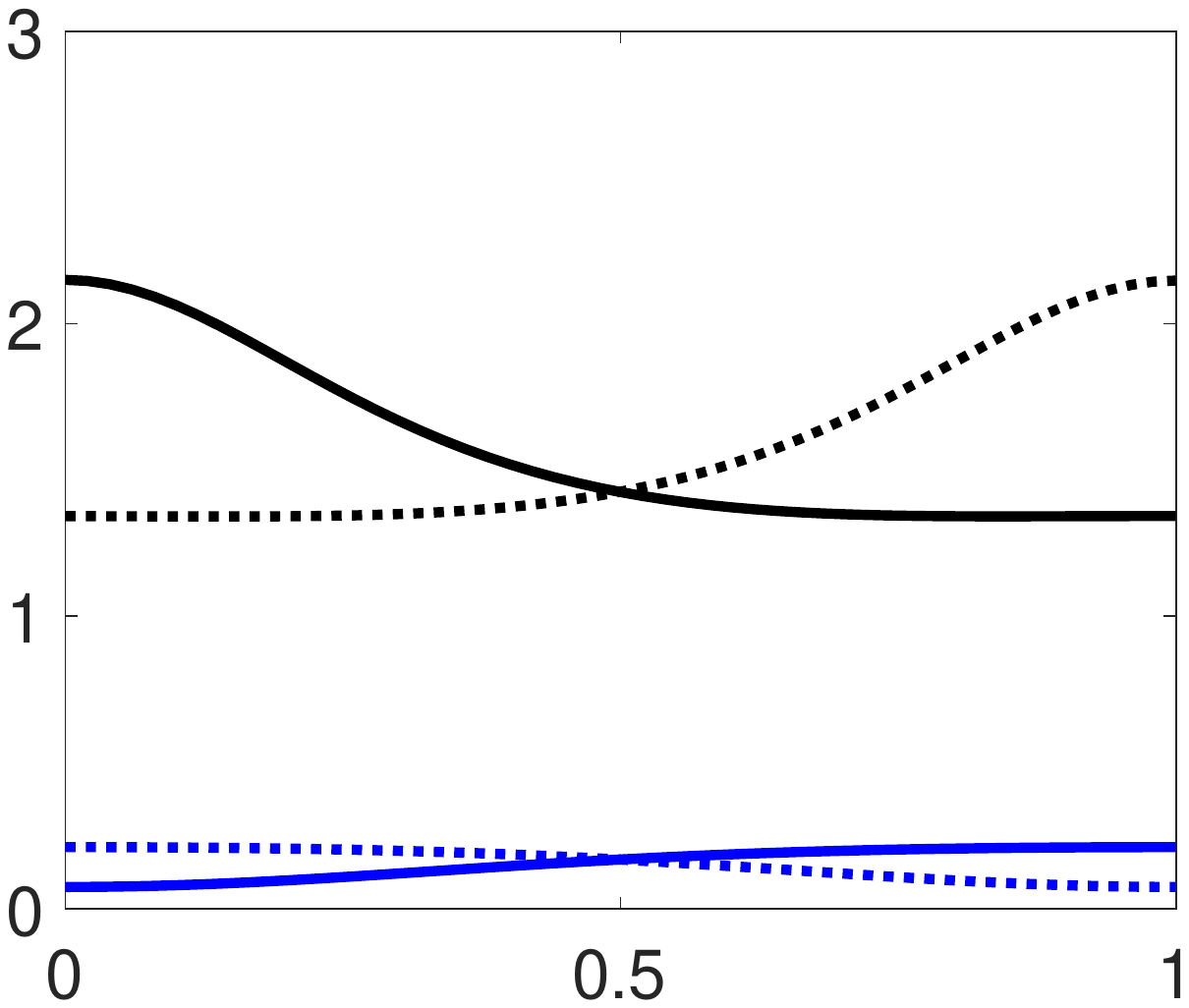}
         \put(-2,30){\rotatebox{90}{$u,\,v$}}
         \put(80,-2){$x$}
         \end{overpic}}
\subfloat[red branch\label{sor}]{
       \begin{overpic}[width=0.3\textwidth,tics=10,trim=90 240 100 250,clip]{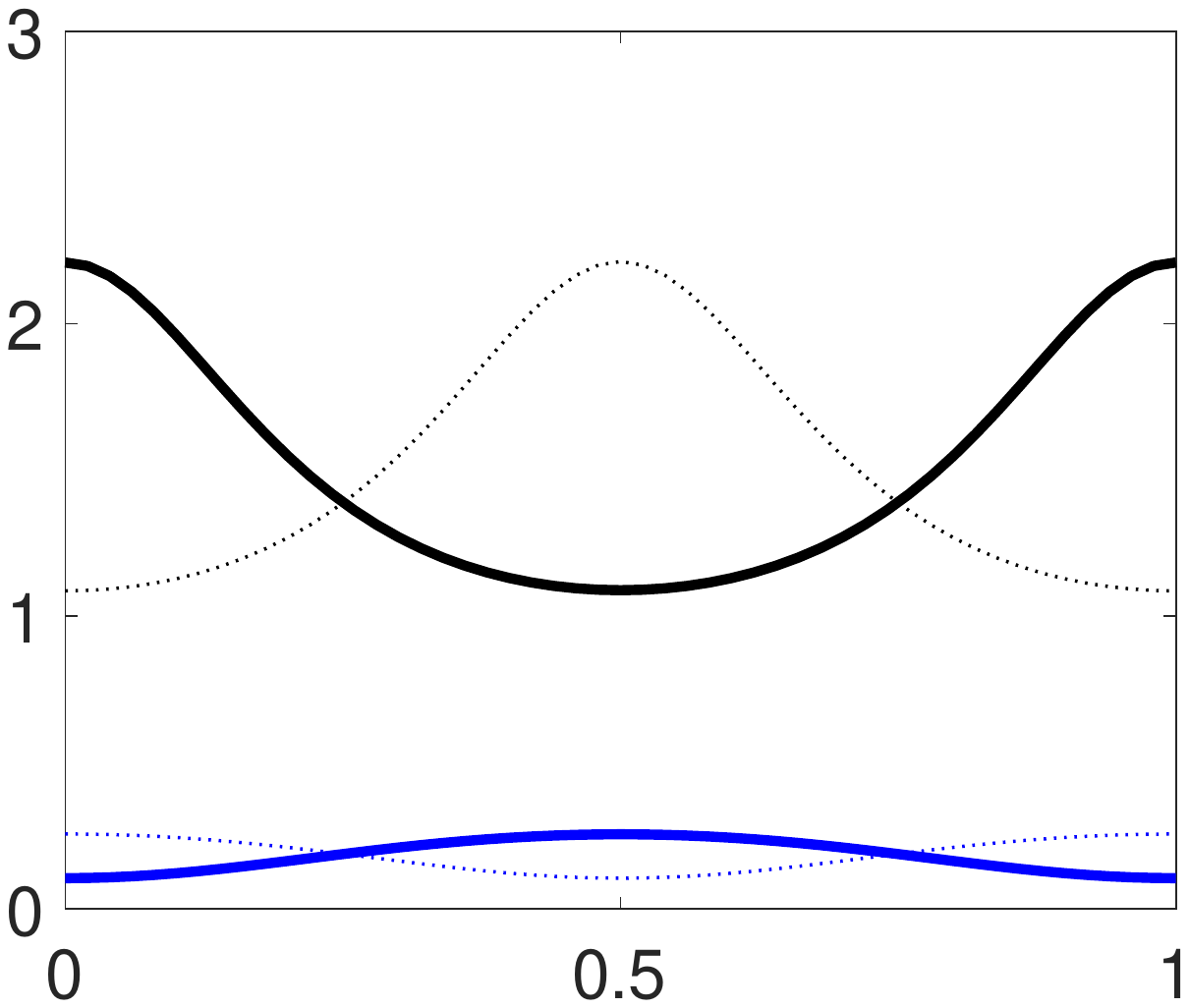}
         \put(-2,30){\rotatebox{90}{$u,\,v$}}
         \put(80,-2){$x$}
         \end{overpic}}
\subfloat[green branch\label{solg}]{
       \begin{overpic}[width=0.3\textwidth,tics=10,trim=90 240 100 250,clip]{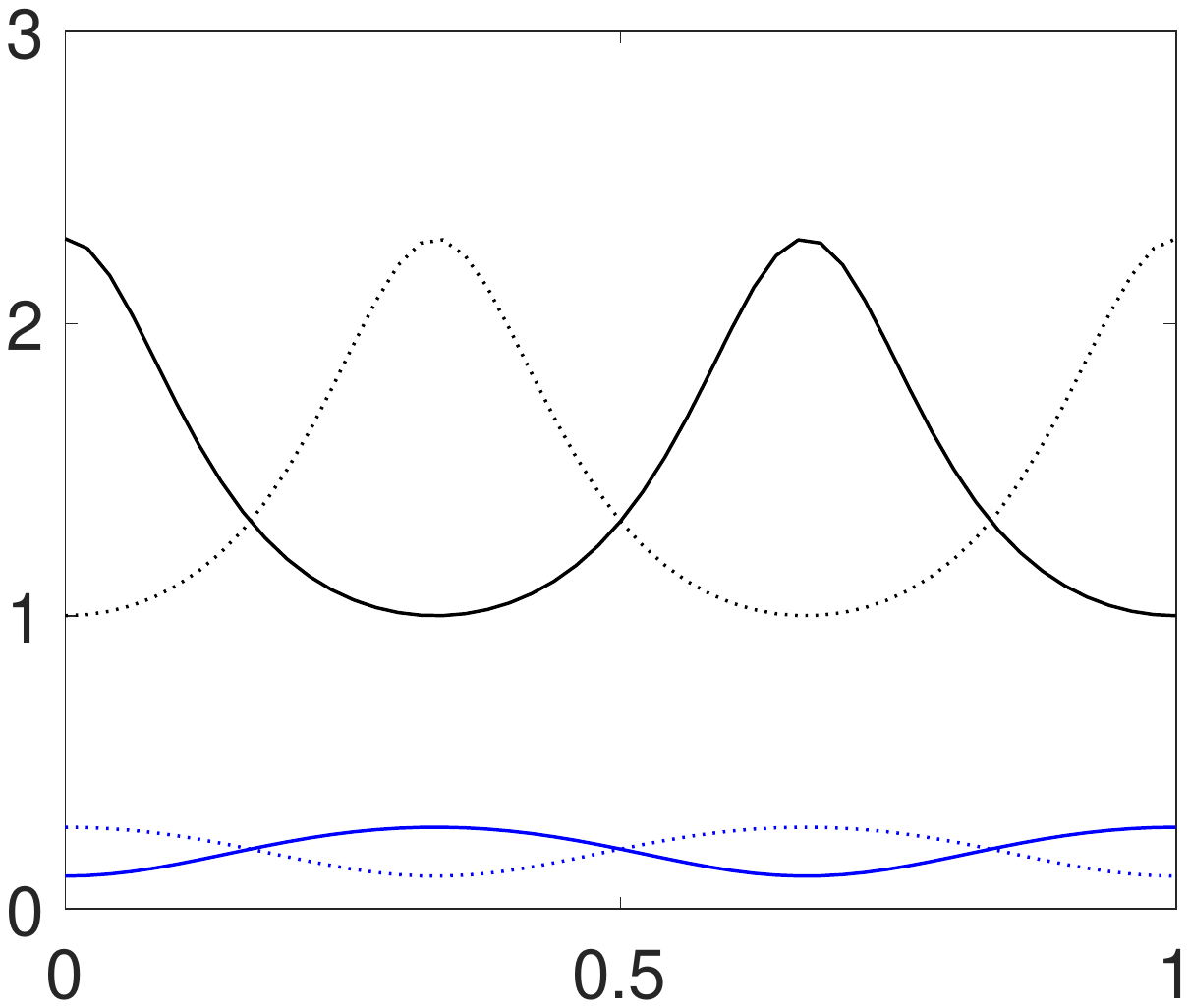}
         \put(-2,30){\rotatebox{90}{$u,\,v$}}
         \put(80,-2){$x$}
         \end{overpic}}\\
\subfloat[yellow branch\label{soly}]{
       \begin{overpic}[width=0.3\textwidth,tics=10,trim=90 240 100 250,clip]{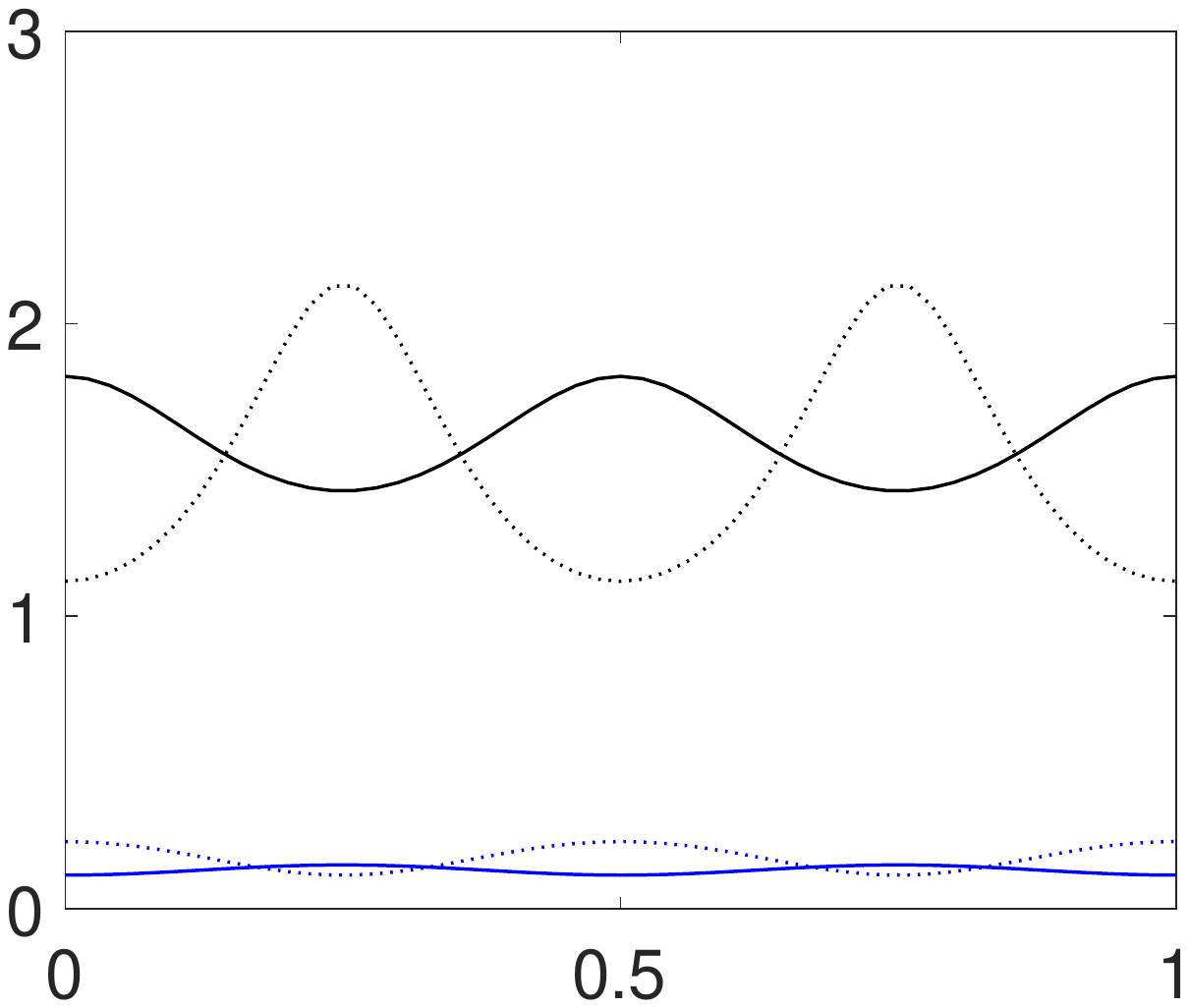}
         \put(-2,30){\rotatebox{90}{$u,\,v$}}
         \put(80,-2){$x$}
         \end{overpic}}
\subfloat[cyan branch\label{solc}]{
       \begin{overpic}[width=0.3\textwidth,tics=10,trim=90 240 100 250,clip]{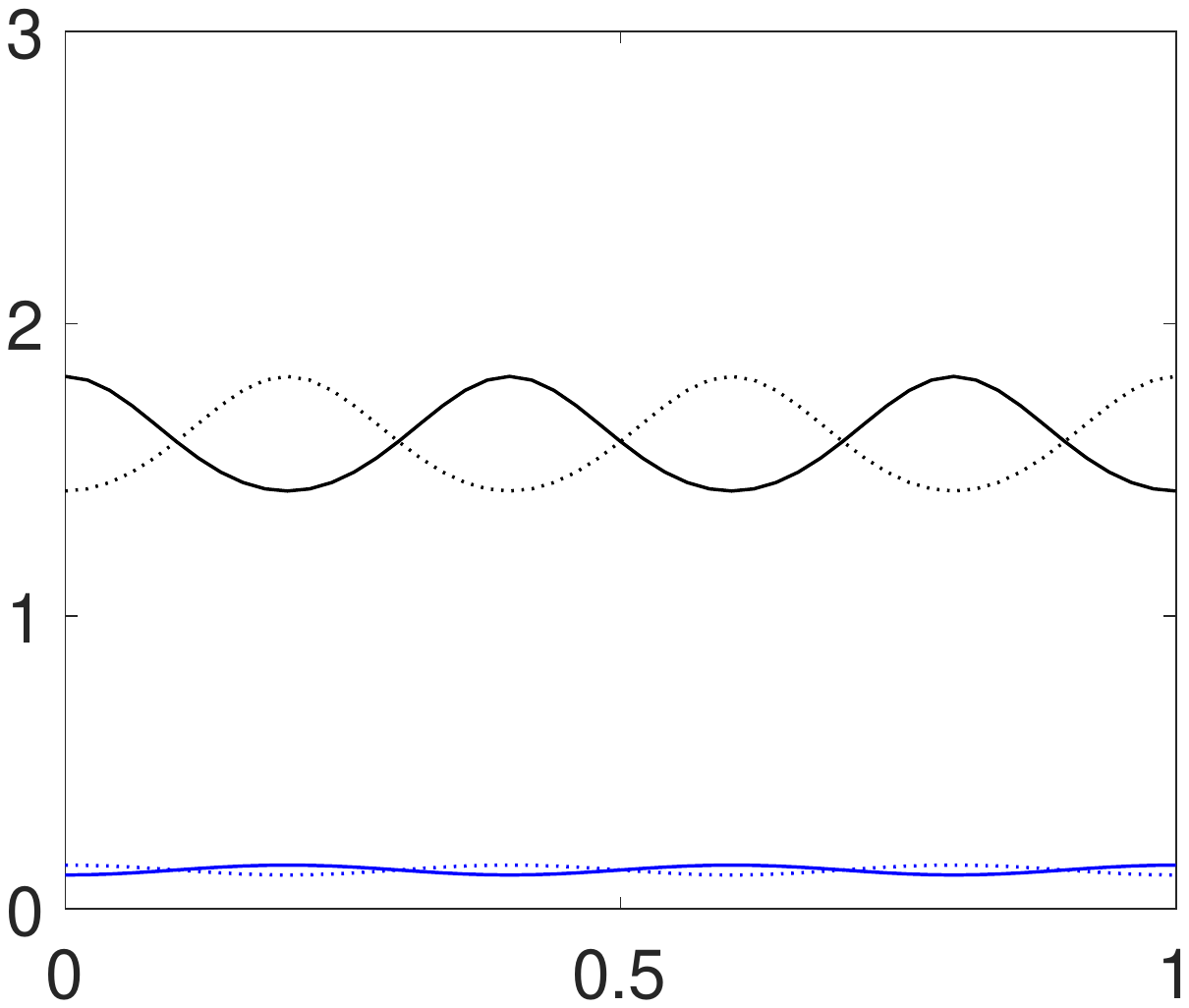}
         \put(-2,30){\rotatebox{90}{$u,\,v$}}
         \put(80,-2){$x$}
         \end{overpic}}
\subfloat[violet branch\label{solv}]{
       \begin{overpic}[width=0.3\textwidth,tics=10,trim=90 240 100 250,clip]{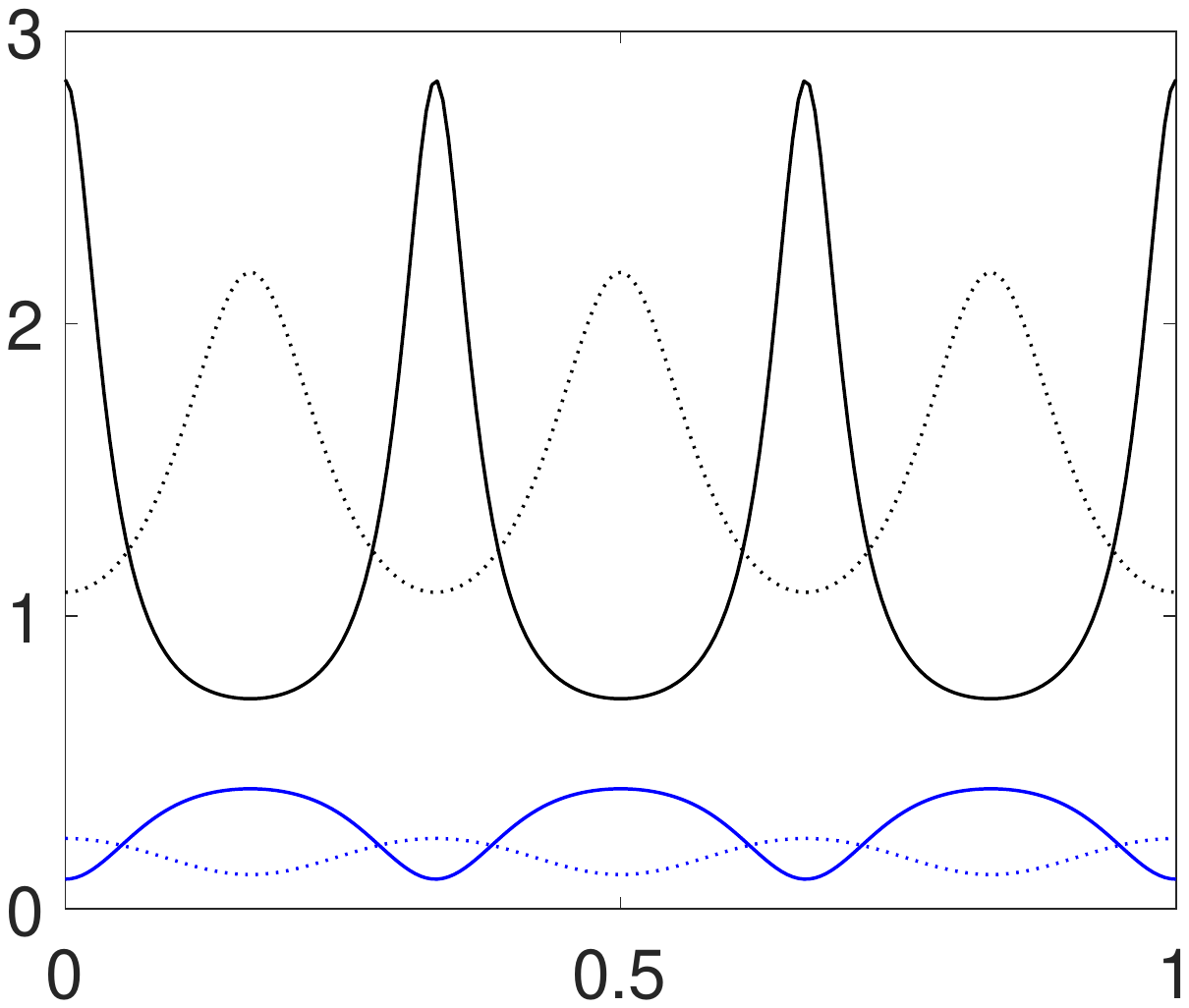}
         \put(-2,30){\rotatebox{90}{$u,\,v$}}
         \put(80,-2){$x$}
         \end{overpic}}\\[-0.5cm]
\begin{multicols}{2}    
\subfloat[magenta branch\label{solm}]{
       \begin{overpic}[width=0.3\textwidth,tics=10,trim=90 240 100 250,clip]{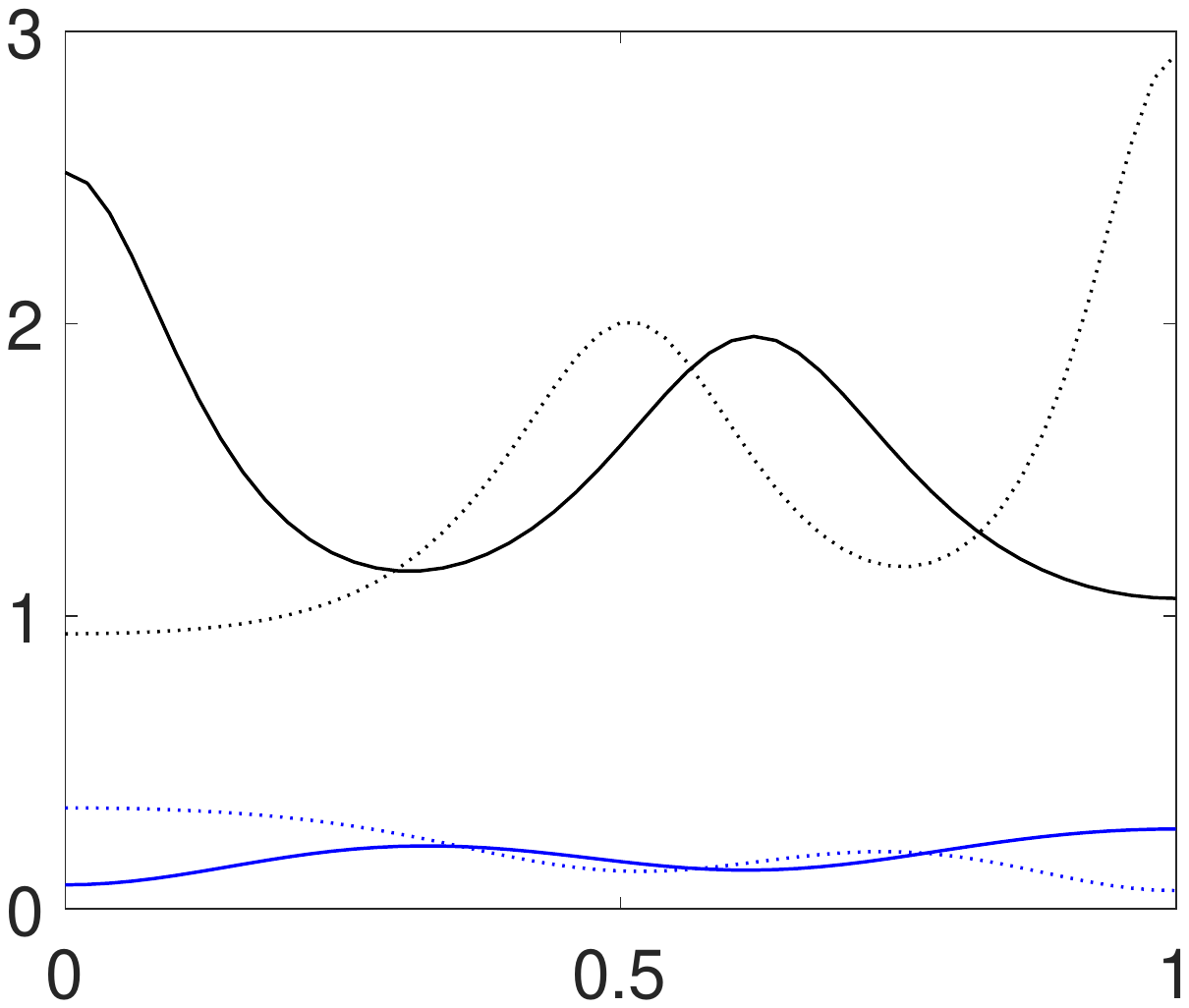}
         \put(-2,30){\rotatebox{90}{$u,\,v$}}
         \put(80,-2){$x$}
         \end{overpic}}\\
\subfloat[orange branch\label{solo}]{
       \begin{overpic}[width=0.3\textwidth,tics=10,trim=90 240 100 250,clip]{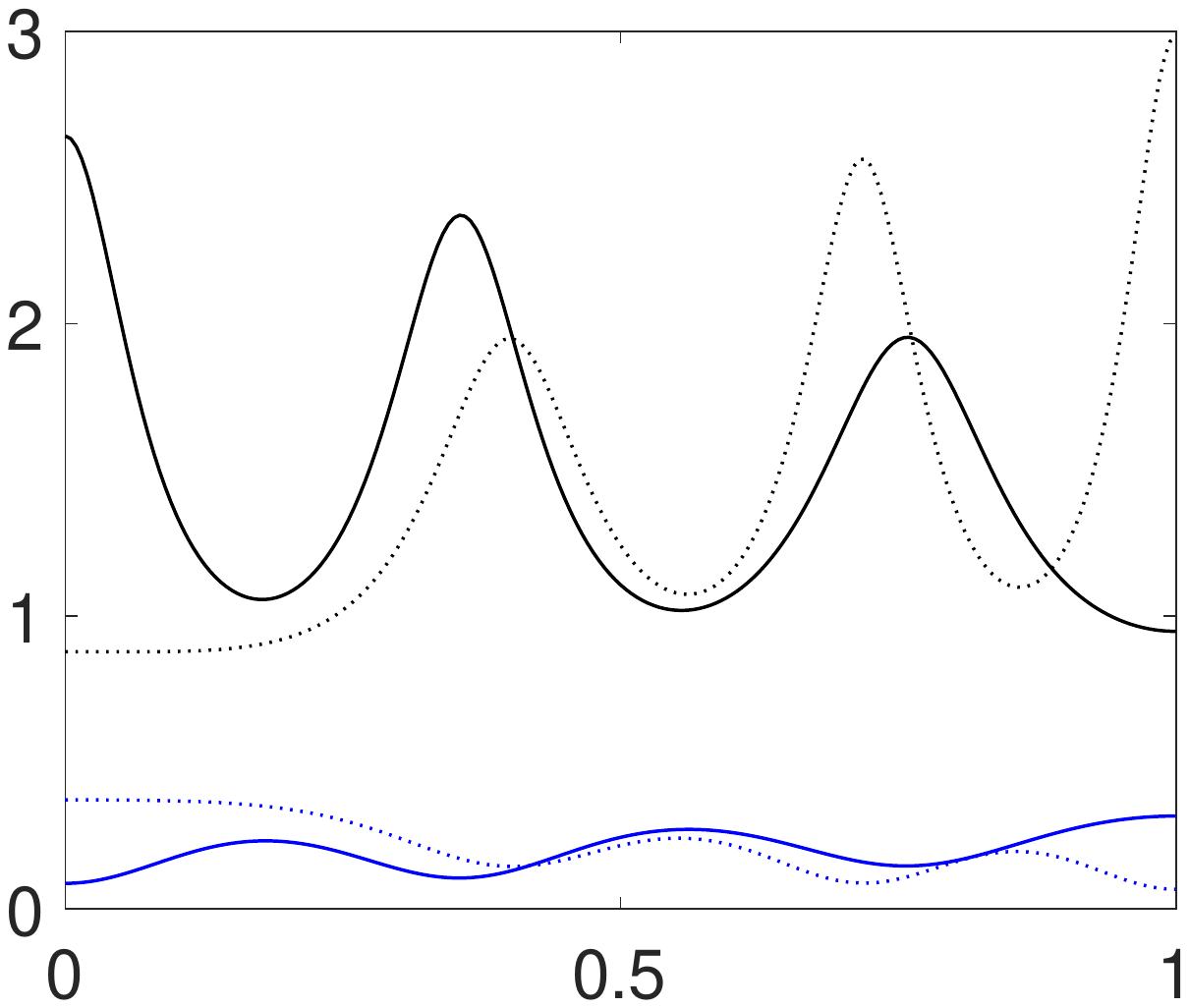}
         \put(-2,30){\rotatebox{90}{$u,\,v$}}
         \put(80,-2){$x$}
         \end{overpic}}\newpage\vspace*{\fill}        
\hspace{-3cm}
\subfloat[\label{bifdiagv0_p}]{
       \begin{overpic}[width=0.6\textwidth,tics=10,trim=90 240 100 250,clip]{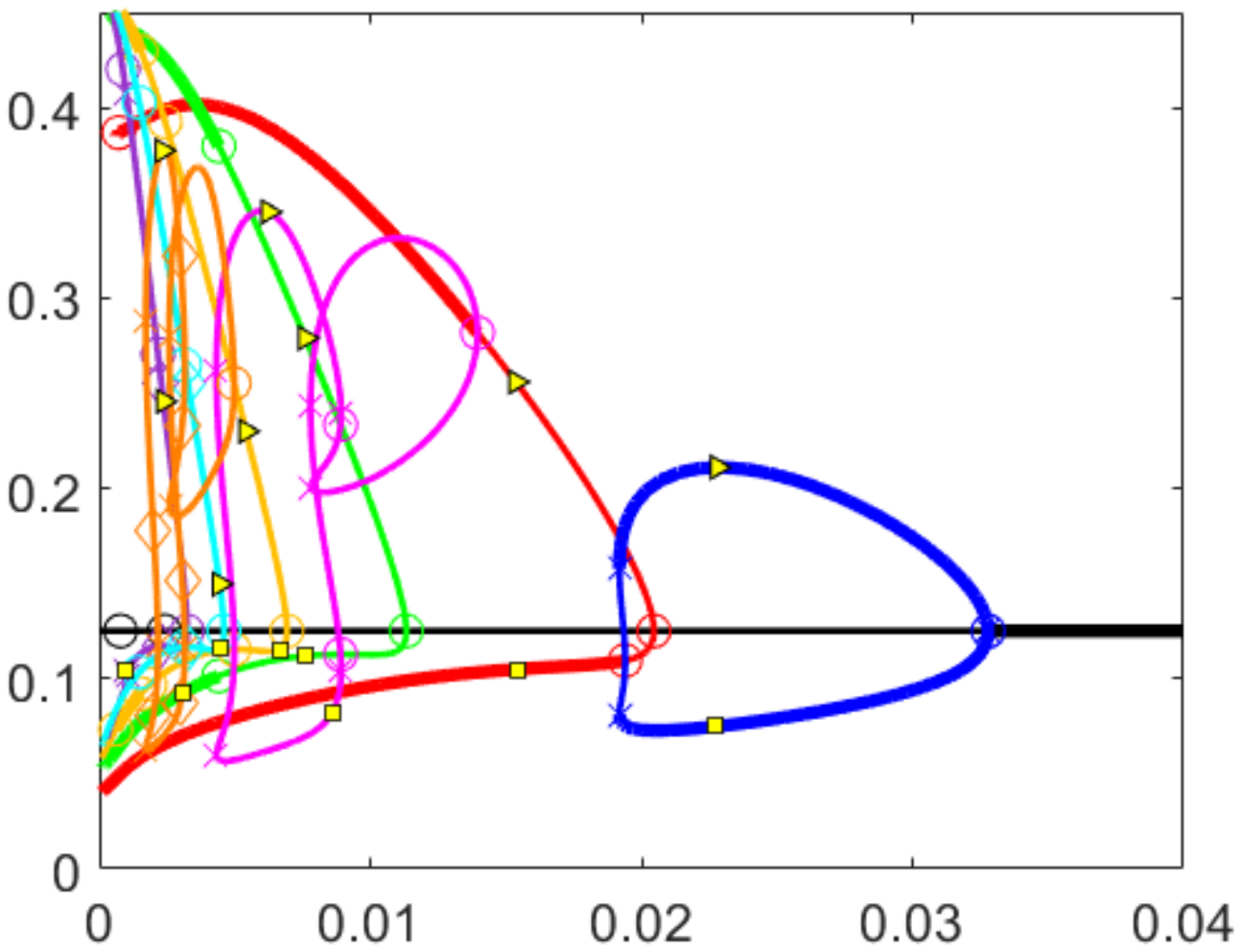}
         \put(0,30){\rotatebox{90}{$v(0)$}}
         \put(80,-2){$d$}
         \put(79,26){$\mathtt{B_1}$}
         \put(56,26){$\mathtt{B_2}$}
         \put(37,26){$\mathtt{B_3}$}
       \end{overpic} }\vspace*{\fill}        
 \end{multicols}
\caption{\protect\subref{solb}--\protect\subref{solo} Different solution types along the branches in the bifurcation diagram \protect\subref{bifdiagv0_p}; the concentrations of $u,v$ are shown in black and blue respectively. Solid lines correspond to points (marked with yellow triangles) located above the homogeneous branch in the bifurcation diagram, while dotted lines to points located below (marked with yellow squares). In the bifurcation diagram, thick lines correspond to stable solutions and thin lines to unstable ones.}
\label{DiagBifCross_sol}
\end{figure}

In Figure~\ref{DiagBifCross_y} the bifurcation diagram of the cross-diffusion system is plotted with respect to different quantities on the $y$-axis, namely $v(0)$, corresponding to the density of the second species at the left boundary value of the domain, and the $L^1$-norm of species $u$. From now on, the homogeneous solution is denoted by the black line, while the other branches correspond to non-homogeneous solutions originating by successive bifurcations. The bifurcations corresponding to branch points are marked with circles, while fold/limit points are marked with crosses. Thick and thin lines denote stable and unstable solutions, respectively. Note that in the $(d,v(0))$-plane at each bifurcation point, two separate branches appear, corresponding to two different solutions. In the $(d,||u||_{L^1})$-plane, the two branches are overlayed, since the solutions on the branches are \emph{symmetric} on the domain. The shape of the non-homogeneous steady states originating along the branches are reported in Figure~\ref{DiagBifCross_sol}. Markers on the branches indicate the positions on the bifurcation diagram of different solutions reported in Figures~\ref{solb}--\ref{solo} (squares corresponds to solid lines, triangles to dashed ones). 

Starting from $d=0.04$ and decreasing its value, we can see that the homogeneous solution is stable and no other solutions are present. At $\mathtt{B_1}$ the homogeneous solution undergoes to a primary bifurcation losing its stability, and two stable non-homogeneous solutions appear (blue lines). Along those branches the density of species $v$ is greater in a part of the domain, while species $u$ occupies the other one. Note that the stable solutions on those branches are symmetric (solid and dashed lines in Figure \ref{solb}). At~$\mathtt{B_2}$, further non-homogeneous solutions appear (red lines), initially unstable; the solutions are again symmetric but now the species density is concentrated either in the central part of the domain or close to the boundary. Those solutions become stable for smaller values of $d$ at a further bifurcation point. We observe the same behavior at the successive bifurcation points from the homogeneous branch: new branches appear (green, yellow and cyan), one each branch the new solutions add half a bump to their shape (Figures \ref{solg}-\ref{solc}). Along the bifurcation branches, the differences between peaks and valleys increases, as the bifurcation parameter $d$ becomes smaller. Finally, for small values of the bifurcation parameter $d$ there can be many different locally stable non-homogeneous solutions. Moreover, there are bifurcation curves connecting three different branches of the homogeneous solutions (magenta and orange lines): along the branches the solution changes shape in order to match the solution profile on the primary branches (Figures \ref{solm}, \ref{solo}). Black circles on the homogeneous branch for small values of $d$ indicate the presence of further bifurcation points that we have not continued.

\subsection{Convergence of the bifurcation diagram}
\label{subsec:1D_conv_d}

The idea to study bifurcation diagrams and the associated singular limit bifurcation diagram as $\varepsilon\ra 0$ has been successfully carried out for fast--slow ODEs in several examples~\cite{GuckenheimerKuehn1,IuorioKuehnSzmolyan}. Yet, more systematic studies for fast-reaction limits for PDEs are missing. Here we numerically compute the convergence of the bifurcation structure of the fast-reaction system \eqref{fast} to the one of the cross-diffusion system~\eqref{cross}, as already reported in~\cite{izuhara2008reaction}, quantifying the convergence of the bifurcation points.

\begin{figure}
\subfloat[$\varepsilon=0.05$\label{fast1D_e0p05}]{
       \begin{overpic}[width=0.5\textwidth,tics=10,trim=90 240 100 250,clip]{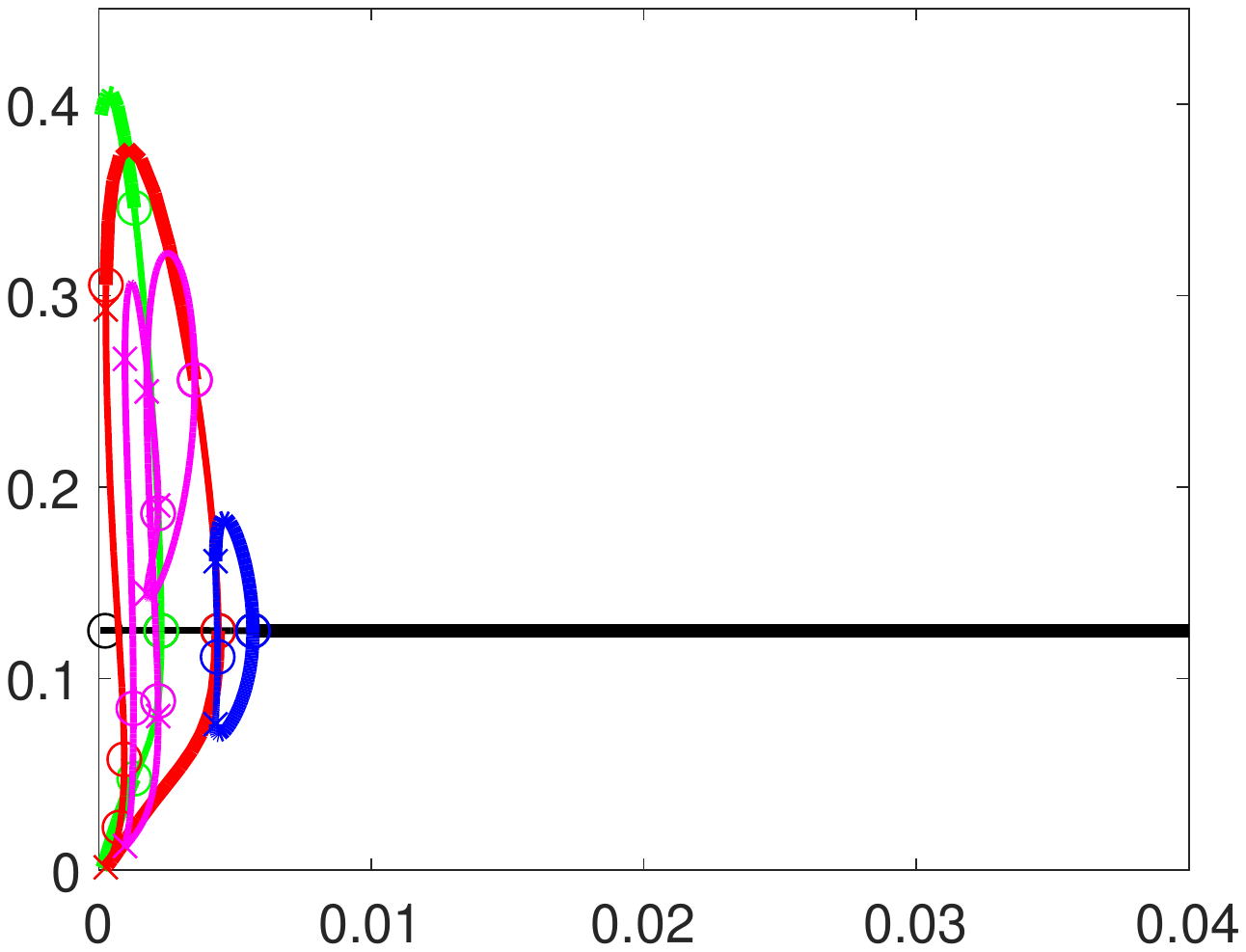}
         \put(0,30){\rotatebox{90}{$v(0)$}}
         \put(80,-2){$d$}
         \put(30,27){$\mathtt{B_1}$}
         \put(27,38){$\mathtt{B_2}$}
         \put(25,48){$\mathtt{B_3}$}
         \put(26,37){\color{gray!30!white}\vector(-1,-3){3.5}}
         \put(26,47){\color{gray!30!white}\vector(-1,-3){6.5}}
       \end{overpic}
}
\subfloat[$\varepsilon=0.01$\label{fast1D_e0p01}]{
       \begin{overpic}[width=0.5\textwidth,tics=10,trim=90 240 100 250,clip]{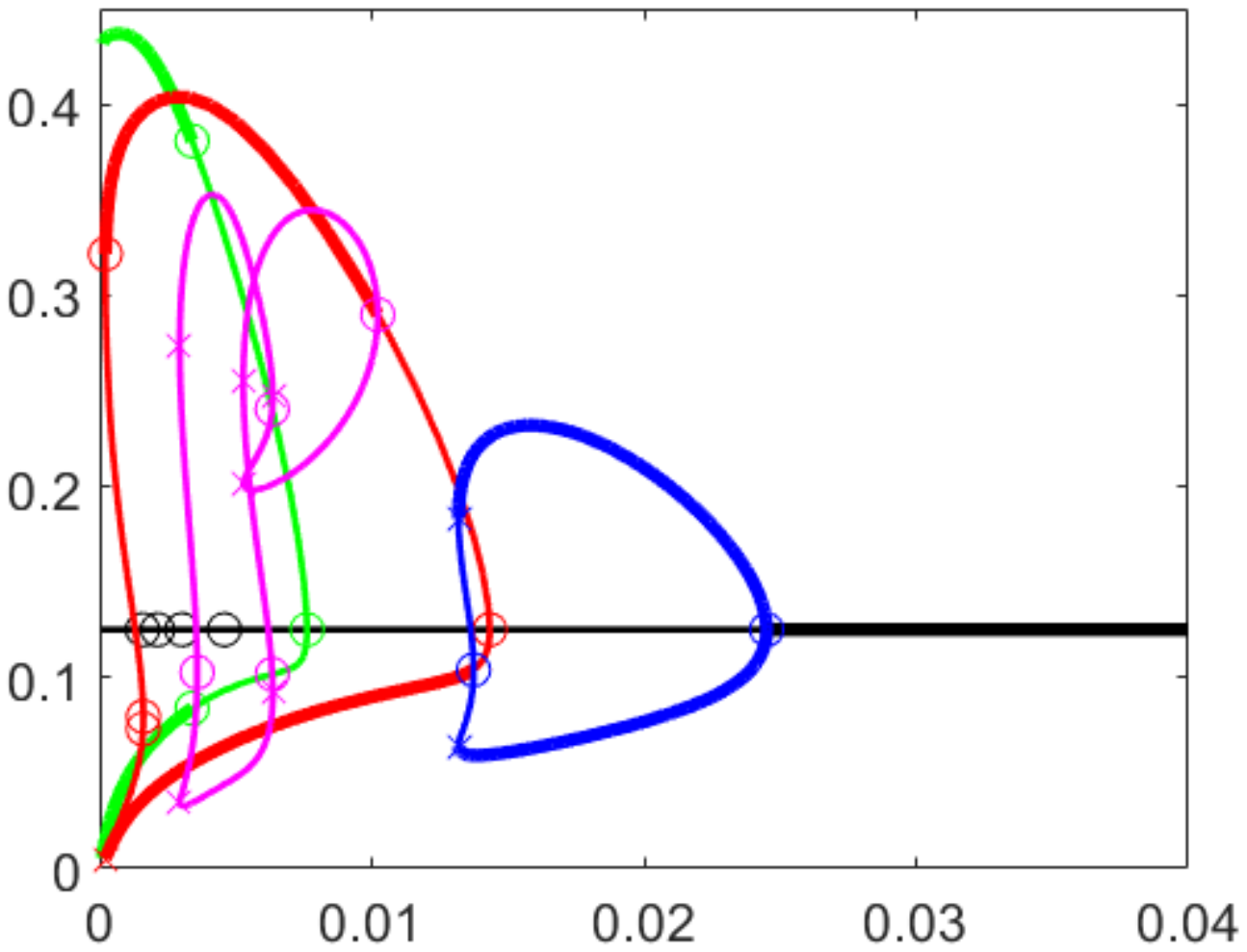}
         \put(0,30){\rotatebox{90}{$v(0)$}}
         \put(80,-2){$d$}
         \put(62,27){$\mathtt{B_1}$}
         \put(43,27){$\mathtt{B_2}$}
         \put(30,27){$\mathtt{B_3}$}
       \end{overpic}
}\\
\subfloat[$\varepsilon=0.001$\label{fast1D_e0p001}]{
       \begin{overpic}[width=0.5\textwidth,tics=10,trim=90 240 100 250,clip]{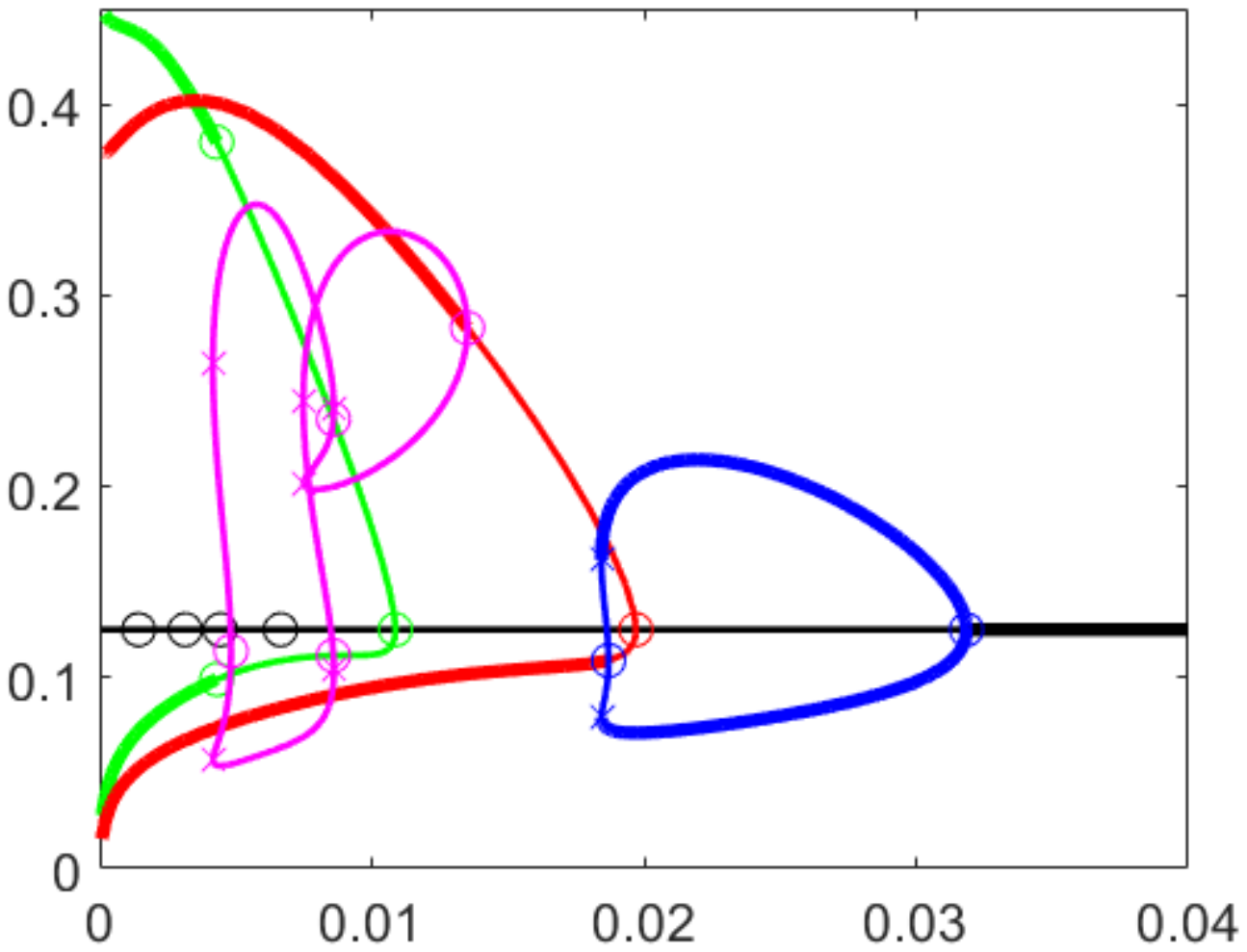}
         \put(0,30){\rotatebox{90}{$v(0)$}}
         \put(80,-2){$d$}
         \put(77,27){$\mathtt{B_1}$}
         \put(54,27){$\mathtt{B_2}$}
         \put(37,27){$\mathtt{B_3}$}
       \end{overpic}
}
\subfloat[$\varepsilon=0$ (CD)\label{cross1D}]{
       \begin{overpic}[width=0.5\textwidth,tics=10,trim=90 240 100 250,clip]{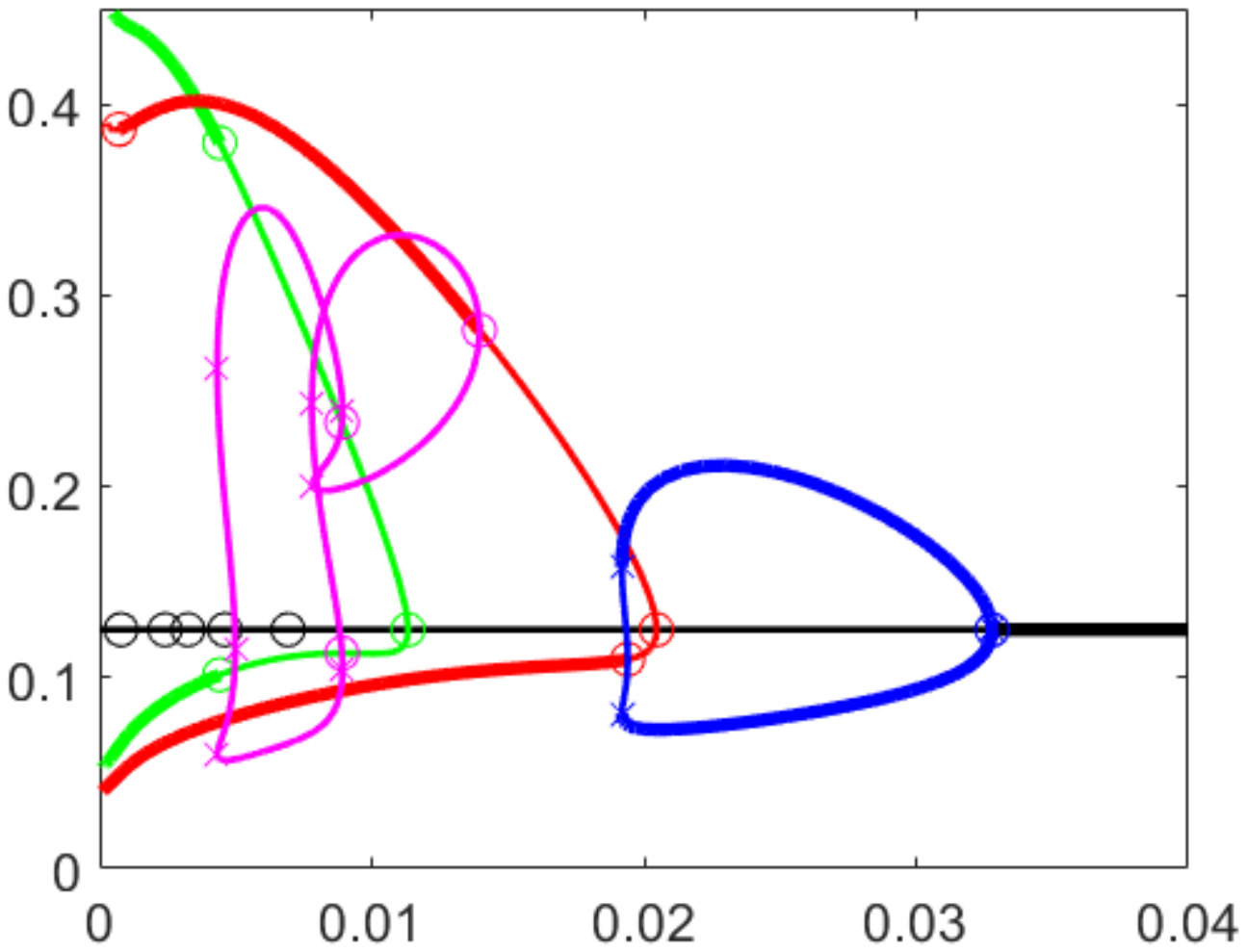}
         \put(0,30){\rotatebox{90}{$v(0)$}}
         \put(80,-2){$d$}
         \put(79,27){$\mathtt{B_1}$}
         \put(56,27){$\mathtt{B_2}$}
         \put(37,27){$\mathtt{B_3}$}
       \end{overpic}
}
\caption{Bifurcation diagrams: \protect\subref{fast1D_e0p05}--\protect\subref{fast1D_e0p001} correspond to the fast--slow system \eqref{fast} for different and decreasing values of $\varepsilon$, while \protect\subref{cross1D} corresponds to the cross-diffusion system \eqref{cross}. The first three bifurcation points on the homogeneous (black) branch are indicated by $\mathtt{B_1}$, $\mathtt{B_2}$, $\mathtt{B_3}$.}
\label{DiagBifConv}
\end{figure}

In Figure \ref{DiagBifConv} the bifurcation diagrams corresponding to the fast-reaction system~\eqref{fast} for different values of $\varepsilon$ are reported, showing the convergence of the bifurcation structure of system~\eqref{fast} to the one of the cross-diffusion system \eqref{cross}. For the sake of clarity of the visualization, we have reported here just the first few branches. In detail, for $\varepsilon=0.1$, system \eqref{fast} does not exhibit non-constant solutions (or they exist for small values of the parameter $d$ and it is difficult to numerically detect them), and the homogeneous steady state is stable for all the values of the bifurcation parameter~$d$. For $\varepsilon=0.05$, the homogeneous steady state becomes unstable for small values of $d$ but we can see that the bifurcation structure corresponding to the fast-reaction system is already qualitatively similar to the one of the cross-diffusion system, but it is squeezed into a small region near $d=0$. The stability properties also match with the cross-diffusion bifurcation structure. For smaller values of $\varepsilon$ the bifurcation structure is expanding to the right, towards the bifurcation structure of the cross-diffusion system. Note that Figure~\ref{fast1D_e0p001}, obtained with $\varepsilon=0.001$, is almost indistinguishable from Figure~\ref{cross1D} corresponding to the cross-diffusion system ($\varepsilon=0$). Hence, from the viewpoint of the bifurcation structure, the three-component fast-reaction system is indeed a good approximation for the cross-diffusion one (at least concerning stationary steady states).

\begin{figure}
\centering
\begin{overpic}[width=0.5\textwidth,tics=10,trim=45 180 55 200,clip]{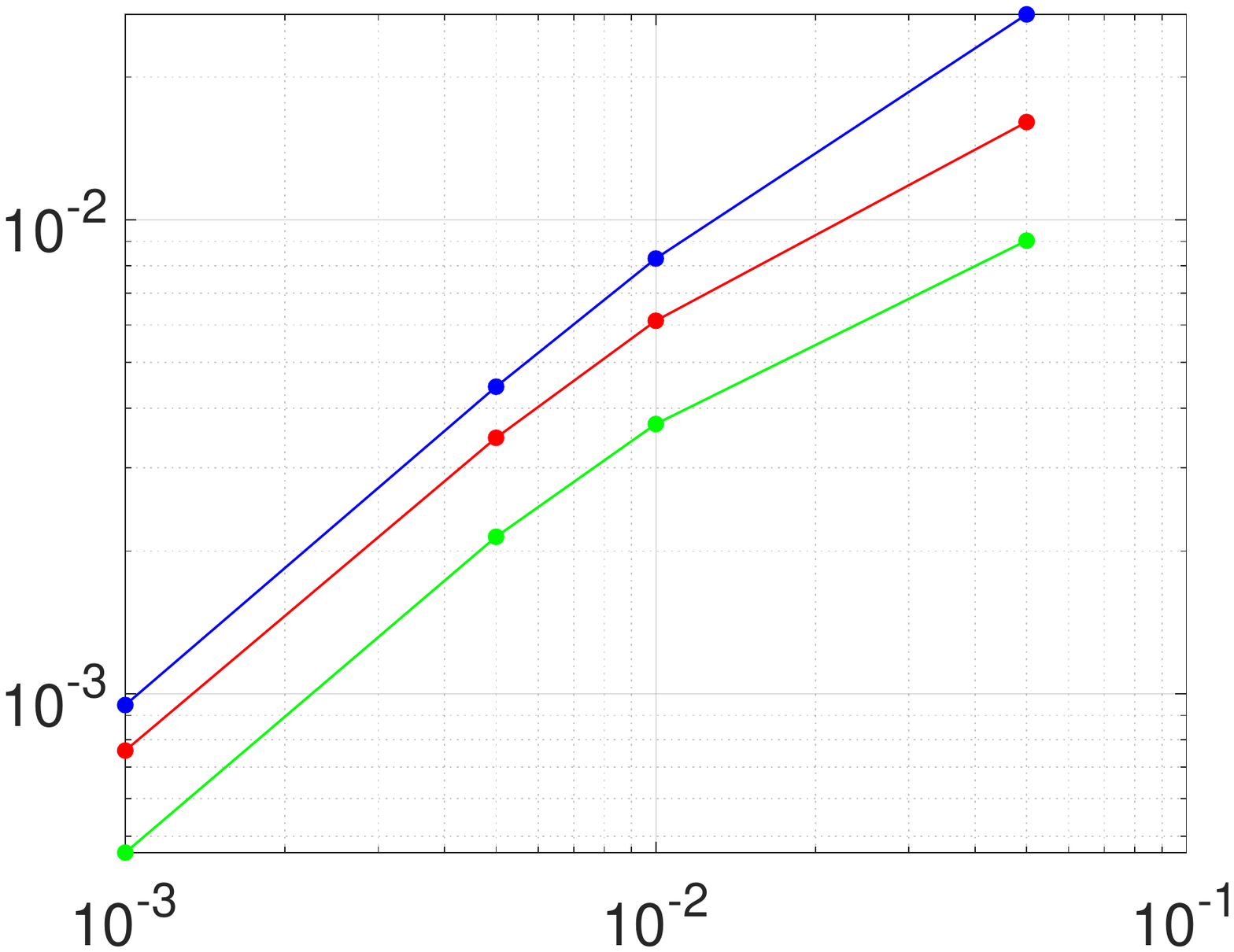}
\put(-10,30){\rotatebox{90}{$|d^0_{\mathtt{B_i}}-d^{\varepsilon}_{\mathtt{B_i}}|$}}
\put(80,0){$\varepsilon$}
\end{overpic}
\caption{Convergence of the first three bifurcation points in loglog scale. We report on the horizontal axis the values of $\varepsilon$, and on the vertical axis the difference between the bifurcation value of the fast--slow system $d^\varepsilon_{\mathtt{B_i}},$ and the corresponding one of the cross-diffusion system $d^0_{\mathtt{B_i}}\, i=1,2,3,$ (blue, red and green refer to  $\mathtt{B_1}$, $\mathtt{B_2}$, $\mathtt{B_3}$ respectively).}
\label{1D_fast_d_conv_loglog}
\end{figure}

In Figure \ref{1D_fast_d_conv_loglog} we quantitatively show the convergence of the first three bifurcation points on the homogeneous branch, namely $\mathtt{B_1}$, $\mathtt{B_2}$, $\mathtt{B_3}$, using the difference between the bifurcation values $d^0_\mathtt{B_i}$ and~$d^\varepsilon_\mathtt{B_i}$. The order of convergence is approximately one.

\subsection{Changing the bifurcation parameter}
\label{subsec:conv_r1}

After we have set up the system in \texttt{pde2path}, we can also change the bifurcation parameter. One possible choice is $r_1$, the growth rate of population $u$. This parameter appears in the reaction part of the system and the homogeneous coexistence state $(u_*,v_*)$ depends on its value. Note that, without diffusion, in the weak competition case (namely $a_1a_2-b_1b_2>0$) the homogeneous coexistence state is positive (and stable) when $b_1/a_2<r_1/r_2<a_1/b_2$, otherwise it is not meaningful and the non-coexistence states $(\bar u,0)$ and $(0,\bar v)$ are stable. With cross-diffusion, it has been shown in \cite{MBCKCS} that the homogeneous equilibrium can be destabilized and stable non-homogeneous solutions arise, which \emph{survive} in a region of the parameter space in which the homogeneous solution is no longer admissible, namely $r_1/r_2>a_1/b_2$. 

\begin{figure}
\hspace{5cm}
\subfloat[\label{sol_bl}]{
       \begin{overpic}[width=0.3\textwidth,tics=10,trim=90 240 100 250,clip]{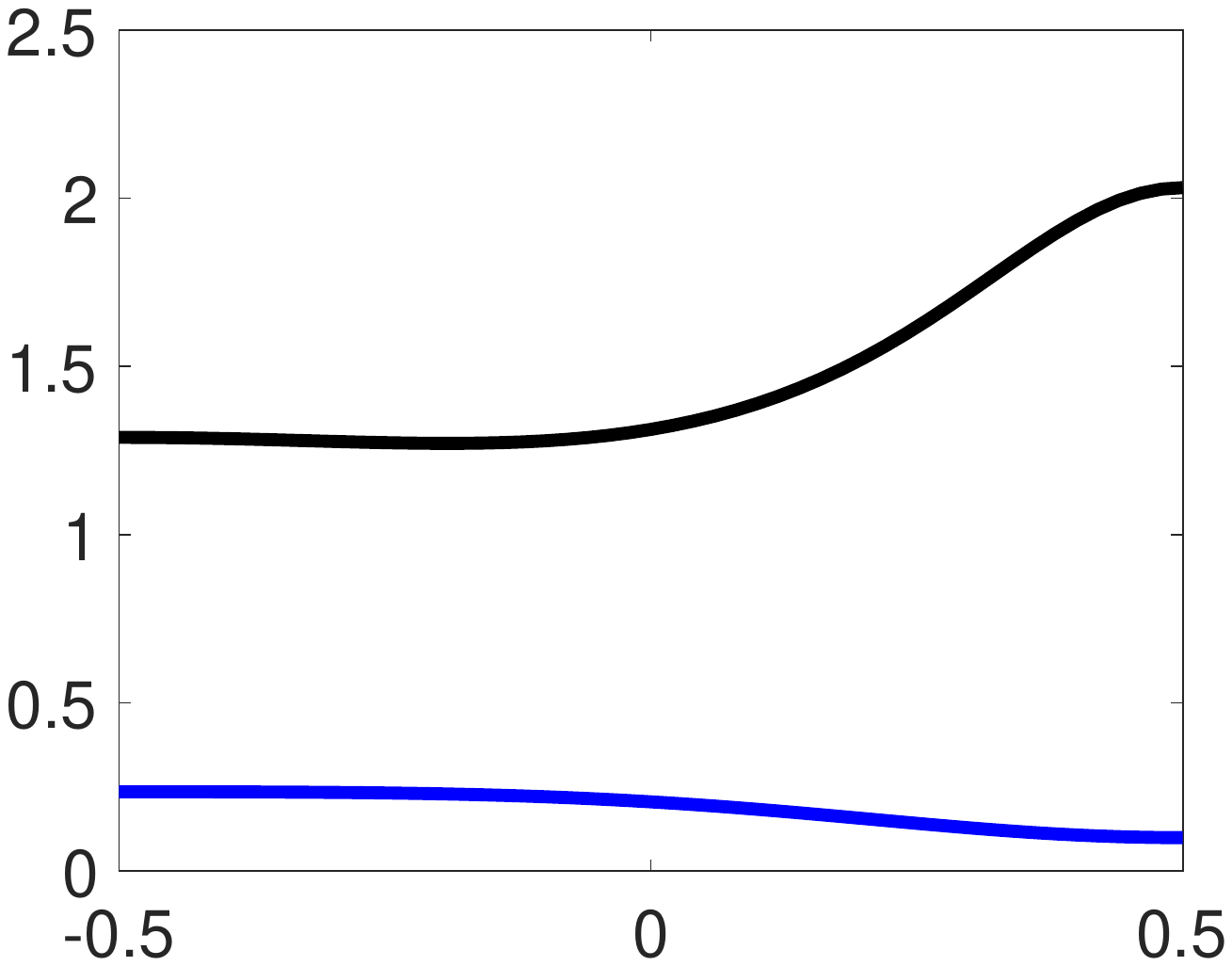}
         \put(-2,30){\rotatebox{90}{$u,\,v$}}
         \put(80,-2){$x$}
         \end{overpic}}
\subfloat[\label{sol_br}]{
       \begin{overpic}[width=0.3\textwidth,tics=10,trim=90 240 100 250,clip]{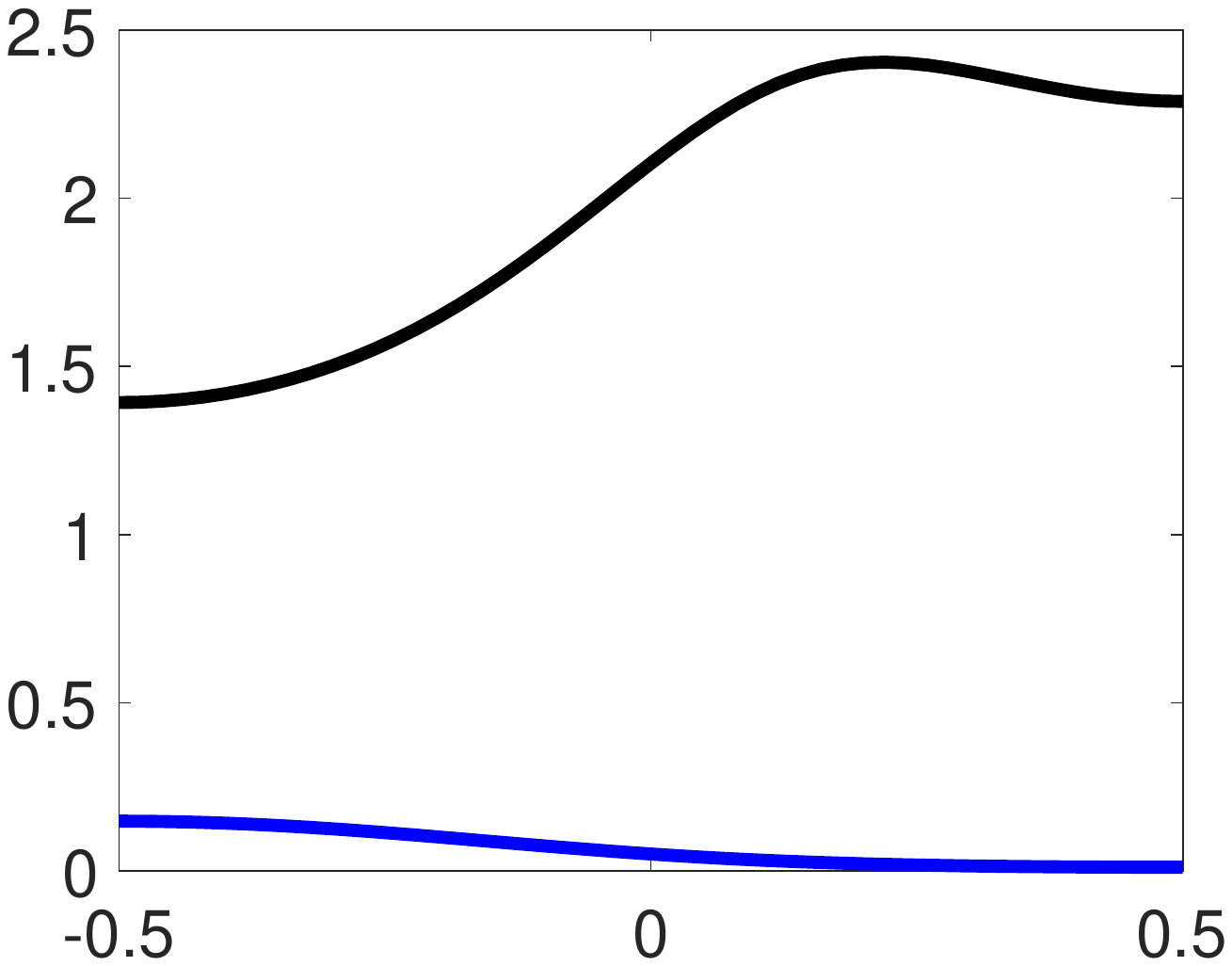}
         \put(-2,30){\rotatebox{90}{$u,\,v$}}
         \put(80,-2){$x$}
         \end{overpic}}\\[-1.5cm]
\begin{multicols}{2}
\subfloat[\label{sol_ru}]{
       \begin{overpic}[width=0.3\textwidth,tics=10,trim=90 240 100 250,clip]{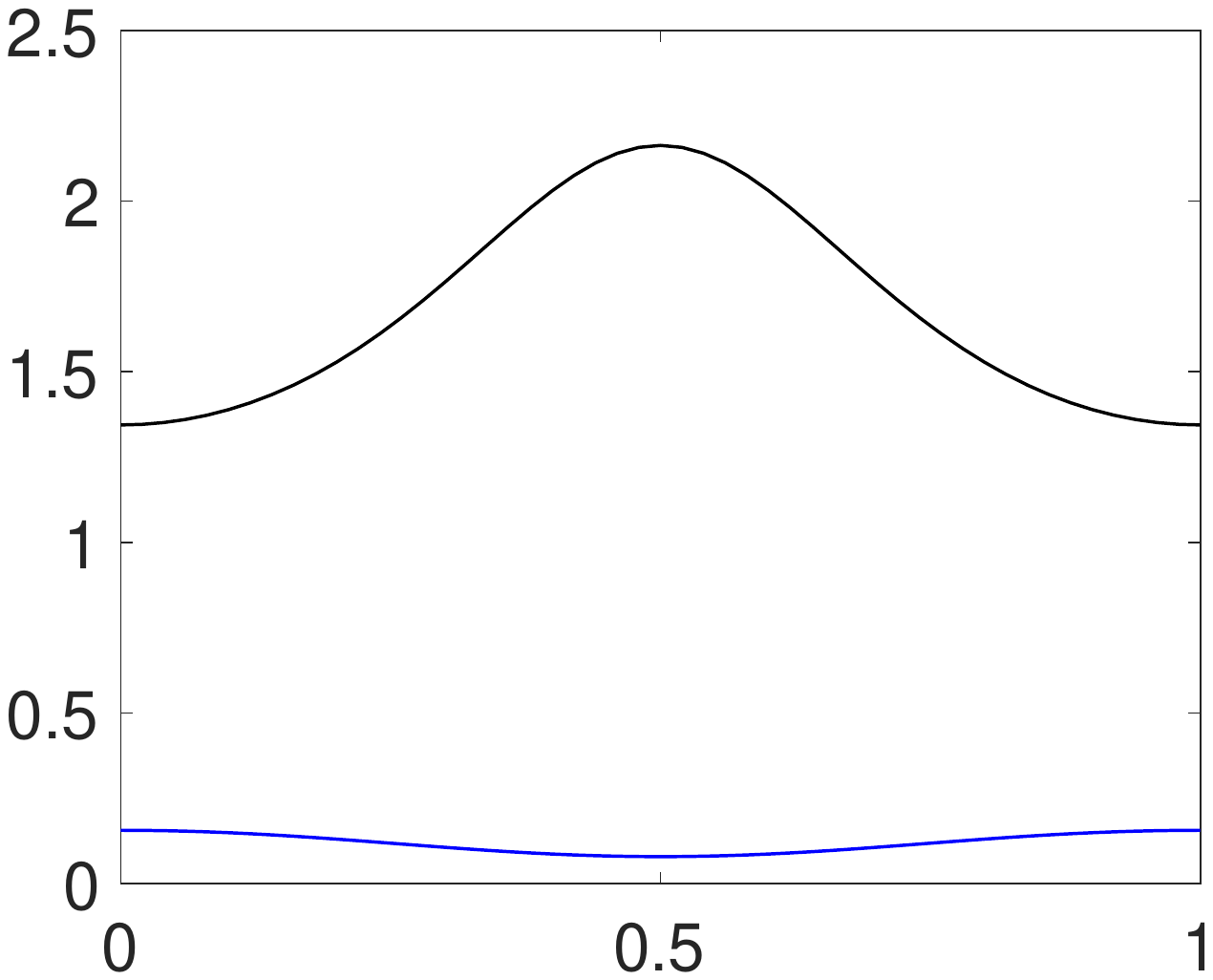}
         \put(-2,30){\rotatebox{90}{$u,\,v$}}
         \put(80,-2){$x$}
         \end{overpic}}\\
\subfloat[\label{sol_rd}]{
       \begin{overpic}[width=0.3\textwidth,tics=10,trim=90 240 100 250,clip]{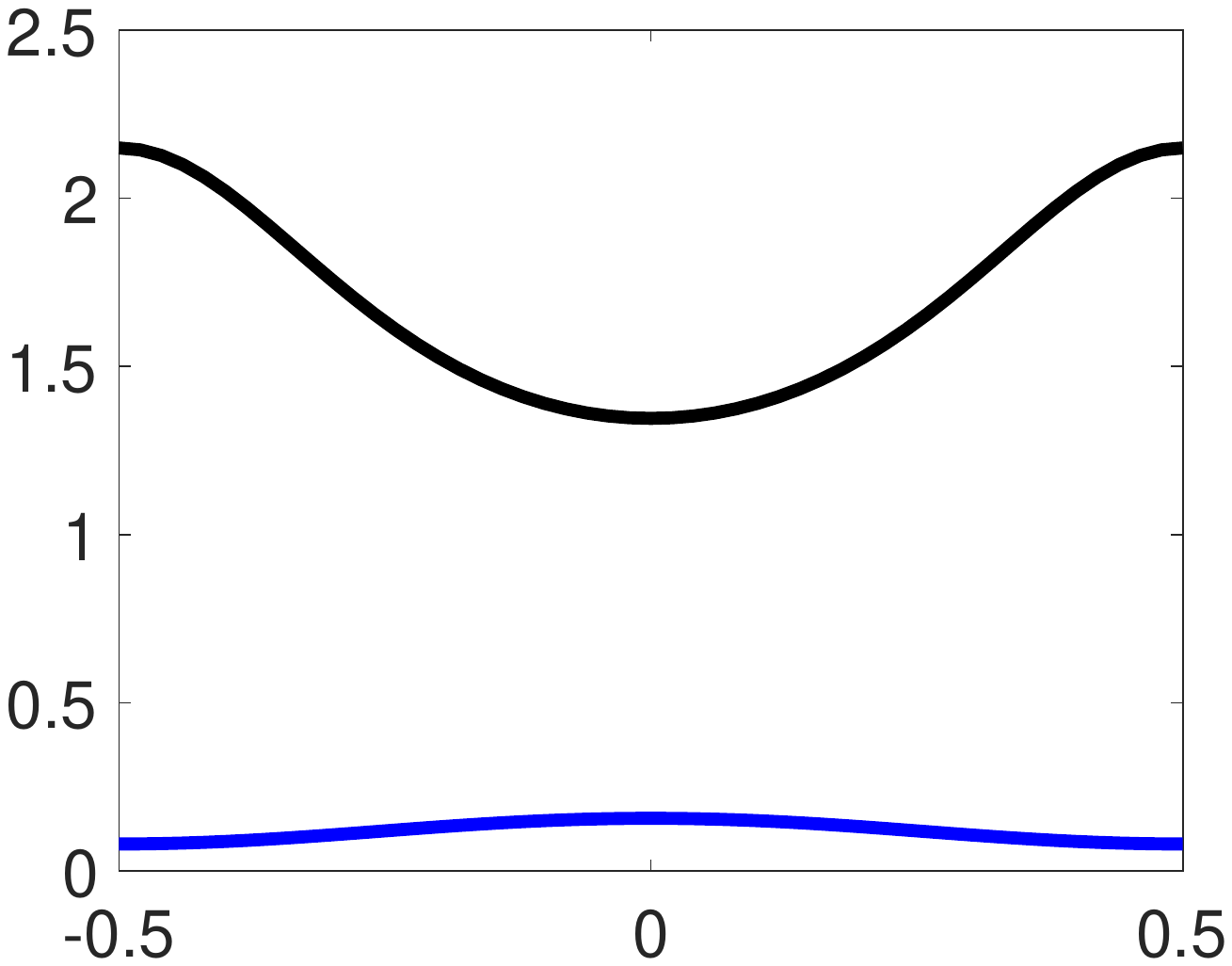}
         \put(-2,30){\rotatebox{90}{$u,\,v$}}
         \put(80,-2){$x$}
\end{overpic}}
\newpage\vspace*{\fill}        
\hspace{-3.3cm}
\subfloat[\label{bifdiagr1_p}]{
       \begin{overpic}[width=0.6\textwidth,tics=10,trim=90 240 100 250,clip]{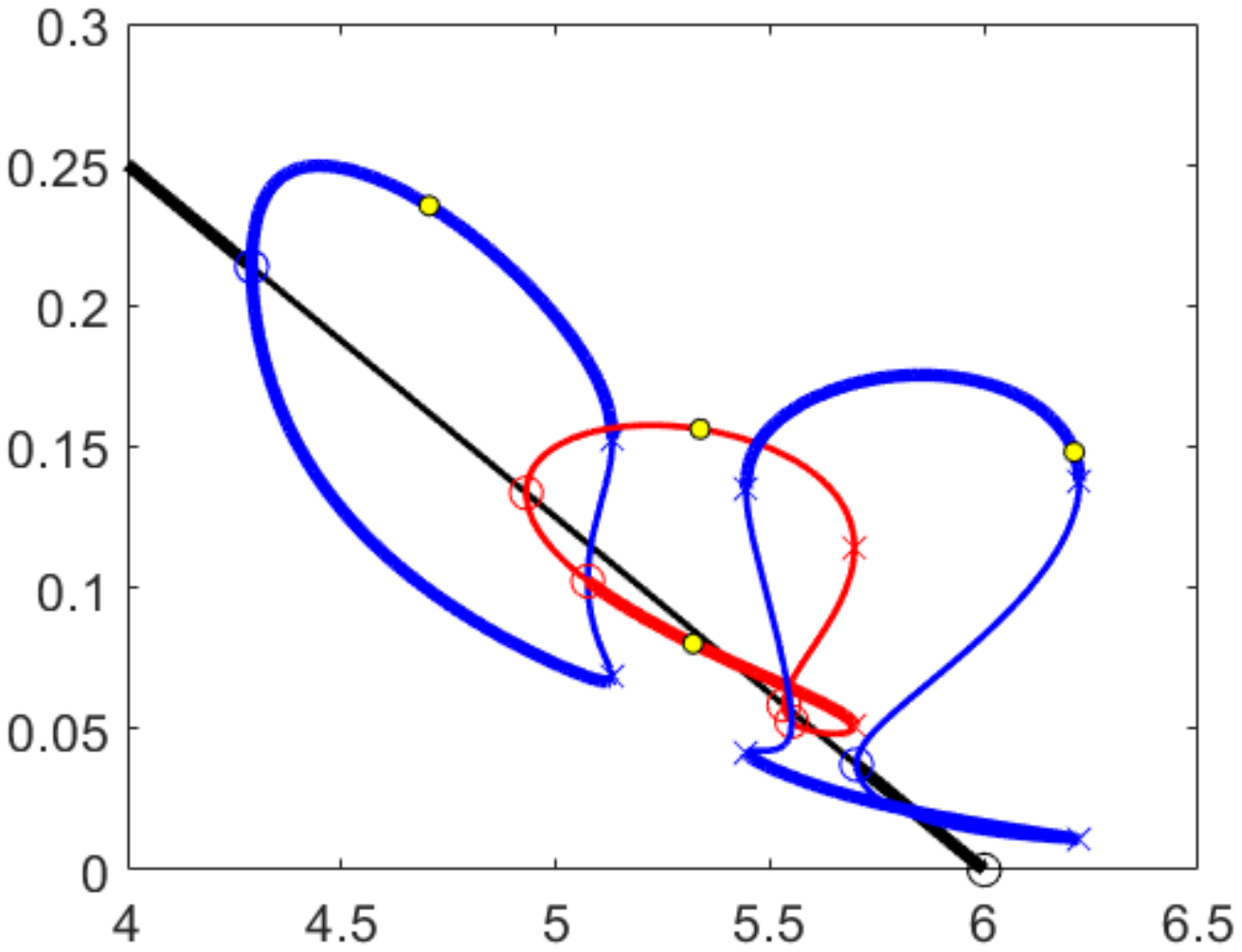}
         \put(0,30){\rotatebox{90}{$v(0)$}}
         \put(80,-2){$r_1$}
         \put(35,58){(a)}
         \put(83,40){(b)}
         \put(53,42){(c)}
         \put(53,18){(d)}
       \end{overpic} }      
\end{multicols}
\caption{\protect\subref{sol_bl}--\protect\subref{sol_rd} Different solution types along the branches in the bifurcation diagram \protect\subref{bifdiagr1_p} of the cross-diffusion system \eqref{cross} with bifurcation parameter $r_1$. Species $u$ and $v$ are denoted in black and blue, respectively. Thick lines correspond to stable solutions, thin lines to unstable ones.}
\label{1D_cross_r1_BifDiag_points}
\end{figure}

In Figure \ref{1D_cross_r1_BifDiag_points} we report the bifurcation diagram of the cross-diffusion system with respect to the bifurcation parameter $r_1$ and some solutions, obtained with the set of parameters of Table \ref{tab:param} and $d=0.02$. Note that it is possible to compare it with the one reported in Figure \ref{bifdiagv0}, since they are different cross-sections of a two-parameter bifurcation surface. The bifurcation diagram is composed by three \emph{rings} and the solutions profile is shown in Figures \ref{sol_bl}--\ref{sol_rd}. In particular, the two outer (blue) rings contain qualitatively similar solutions (Figures \ref{sol_bl} and \ref{sol_br}). Furthermore, Figure \ref{sol_br} shows a stable non-homogeneous steady state corresponding to a parameter value outside of the usual weak competition regime. Furthermore, the (red) branches originating from the second bifurcation point are non-symmetric regarding the stability properties even if they correspond to symmetric solutions (with respect to the homogeneous one) on the domain (Figures \ref{sol_ru} and \ref{sol_rd}). 

\begin{center}
\begin{figure}[!ht]
\vspace{-1.5cm}
\subfloat[$\varepsilon=0.01$\label{fast1D_r1_epsi_0p01}]{
       \begin{overpic}[width=0.5\textwidth,tics=10,trim=90 240 100 250,clip]{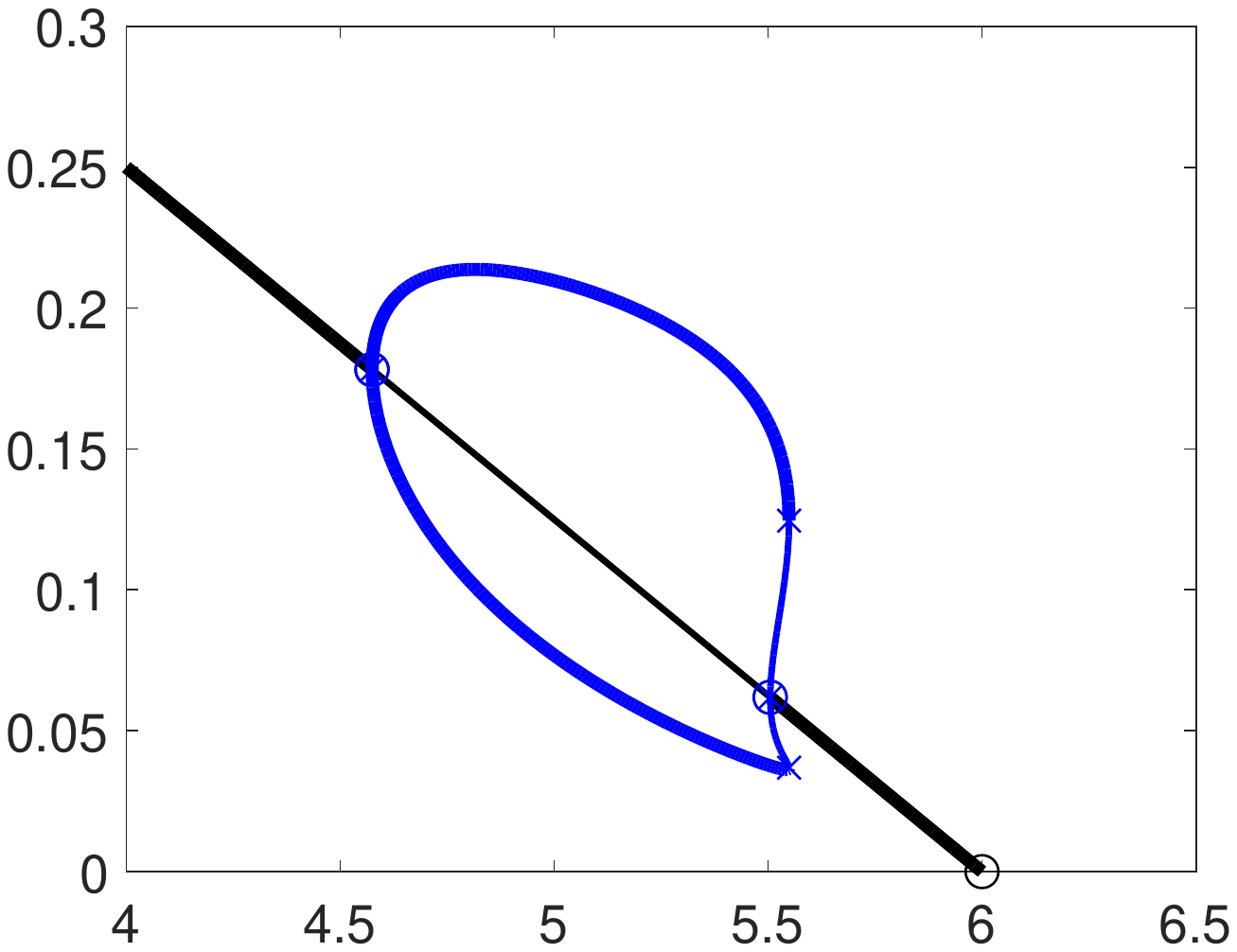}
         \put(0,30){\rotatebox{90}{$v(0)$}}
         \put(80,-2){$r_1$}
       \end{overpic}
}
\subfloat[$\varepsilon=0.005$\label{fast1D_r1_epsi_0p005}]{
       \begin{overpic}[width=0.5\textwidth,tics=10,trim=90 240 100 250,clip]{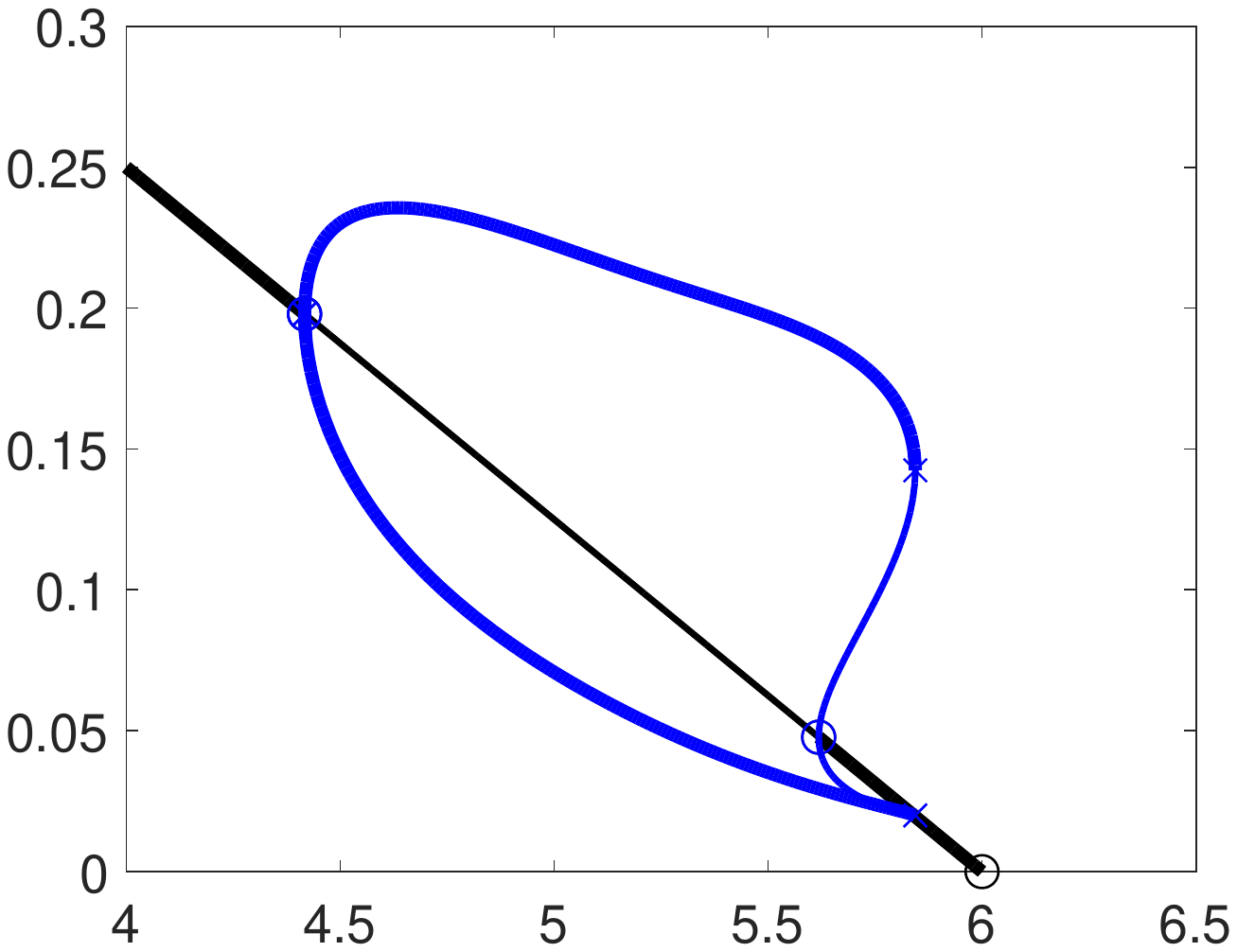}
         \put(0,30){\rotatebox{90}{$v(0)$}}
         \put(80,-2){$r_1$}
       \end{overpic}
}\\[-0.3cm]
\subfloat[$\varepsilon=0.001$\label{fast1D_r1_epsi_0p001}]{
       \begin{overpic}[width=0.5\textwidth,tics=10,trim=90 240 100 250,clip]{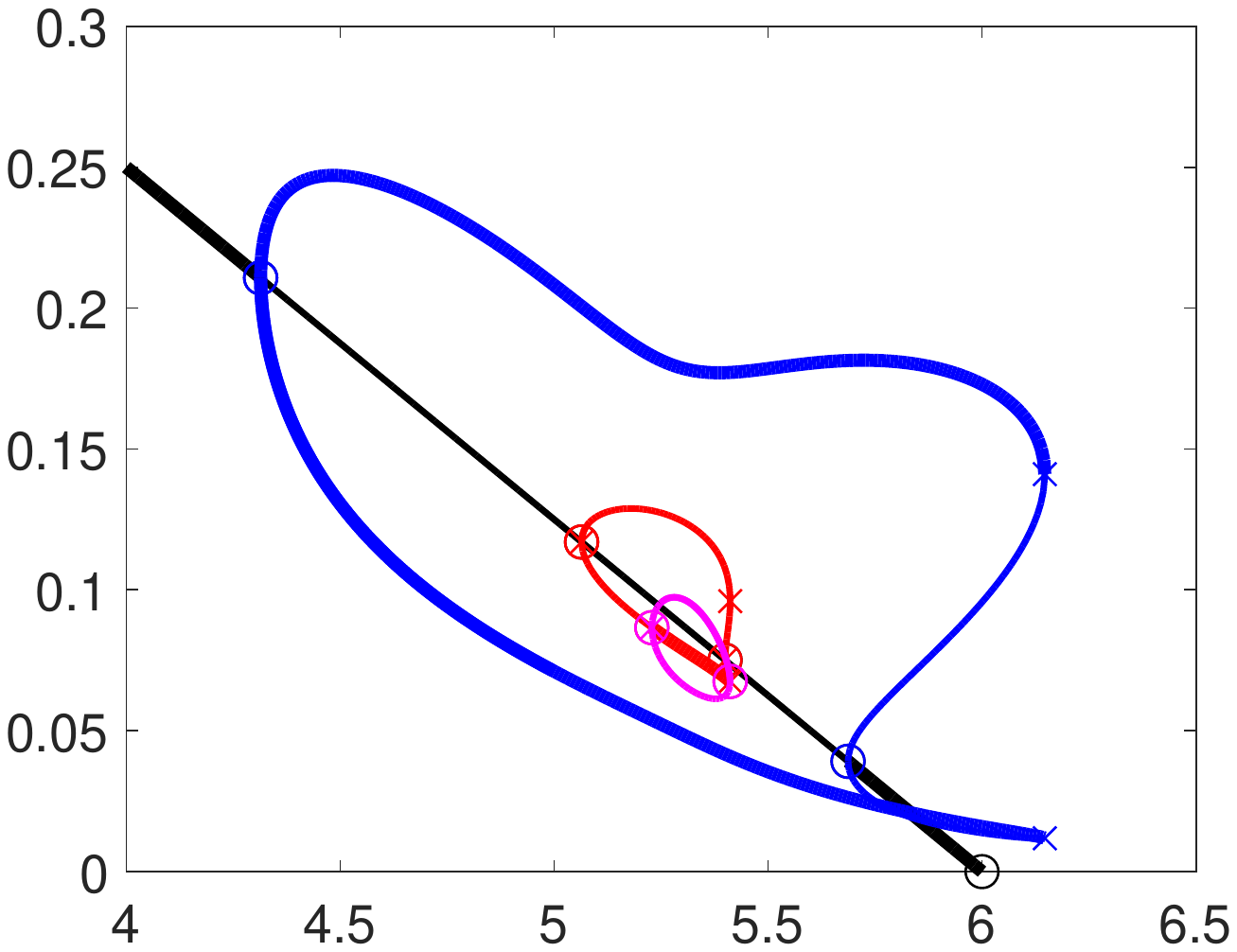}
         \put(0,30){\rotatebox{90}{$v(0)$}}
         \put(80,-2){$r_1$}
       \end{overpic}
}
\subfloat[$\varepsilon=0.0005$\label{fast1D_r1_epsi_0p0005}]{
       \begin{overpic}[width=0.5\textwidth,tics=10,trim=90 240 100 250,clip]{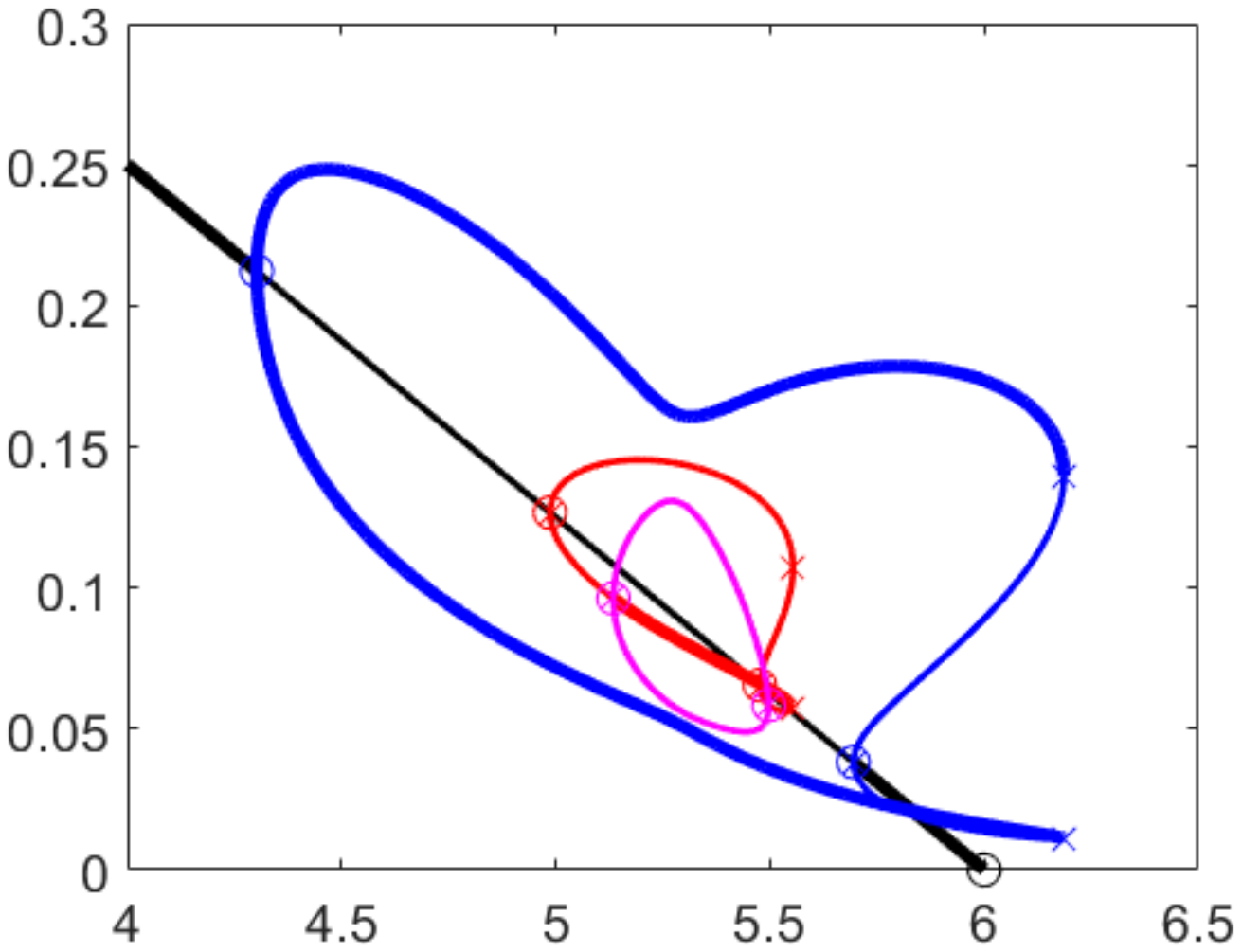}
         \put(0,30){\rotatebox{90}{$v(0)$}}
         \put(80,-2){$r_1$}
       \end{overpic}
}\\[-0.3cm]
\subfloat[$\varepsilon=0.0004$\label{fast1_Dr1_epsi_0p0004}]{
       \begin{overpic}[width=0.5\textwidth,tics=10,trim=90 240 100 250,clip]{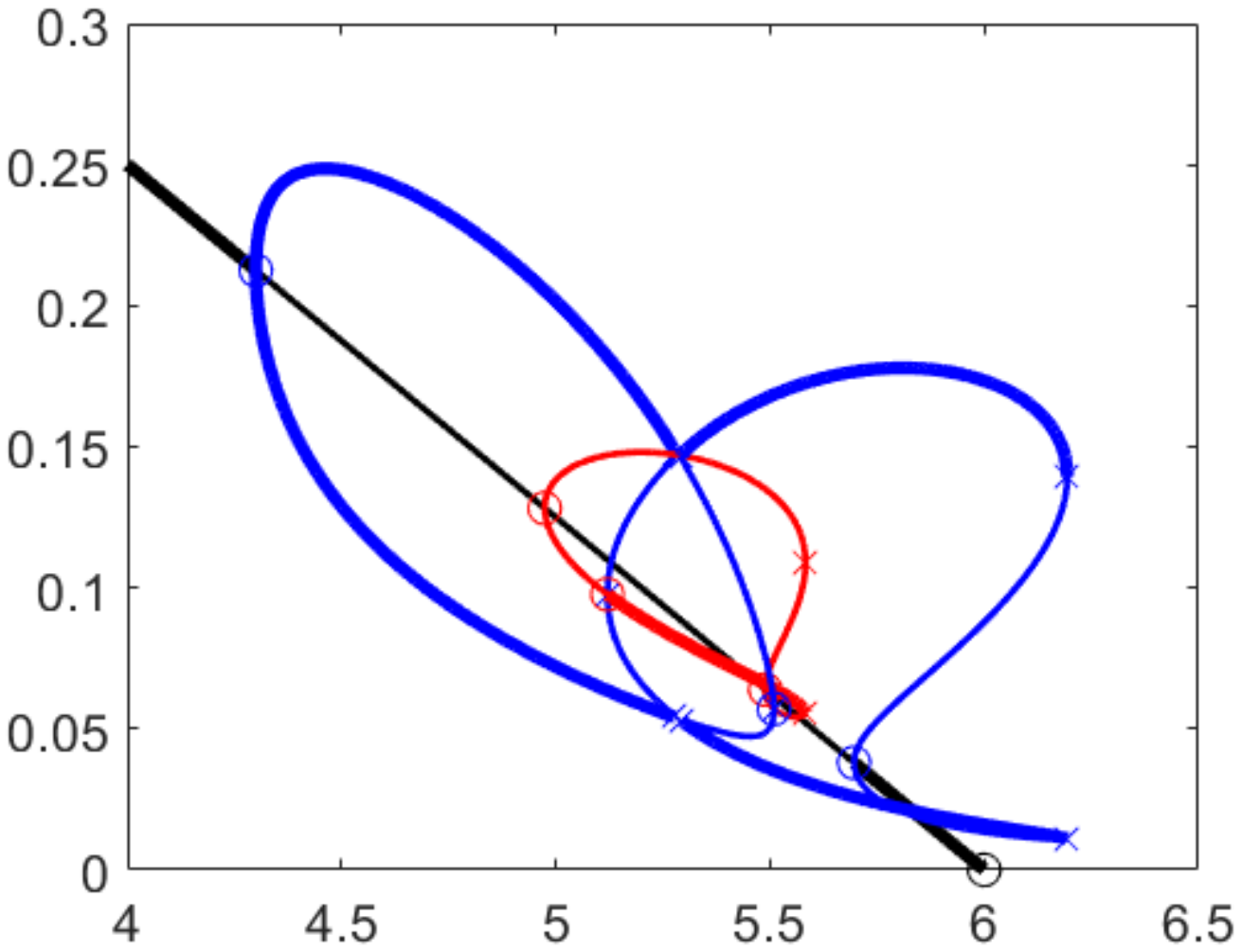}
         \put(0,30){\rotatebox{90}{$v(0)$}}
         \put(80,-2){$r_1$}
       \end{overpic}
}
\subfloat[$\varepsilon=0.0001$\label{fast1D_r1_epsi_0p0001}]{
       \begin{overpic}[width=0.5\textwidth,tics=10,trim=90 240 100 250,clip]{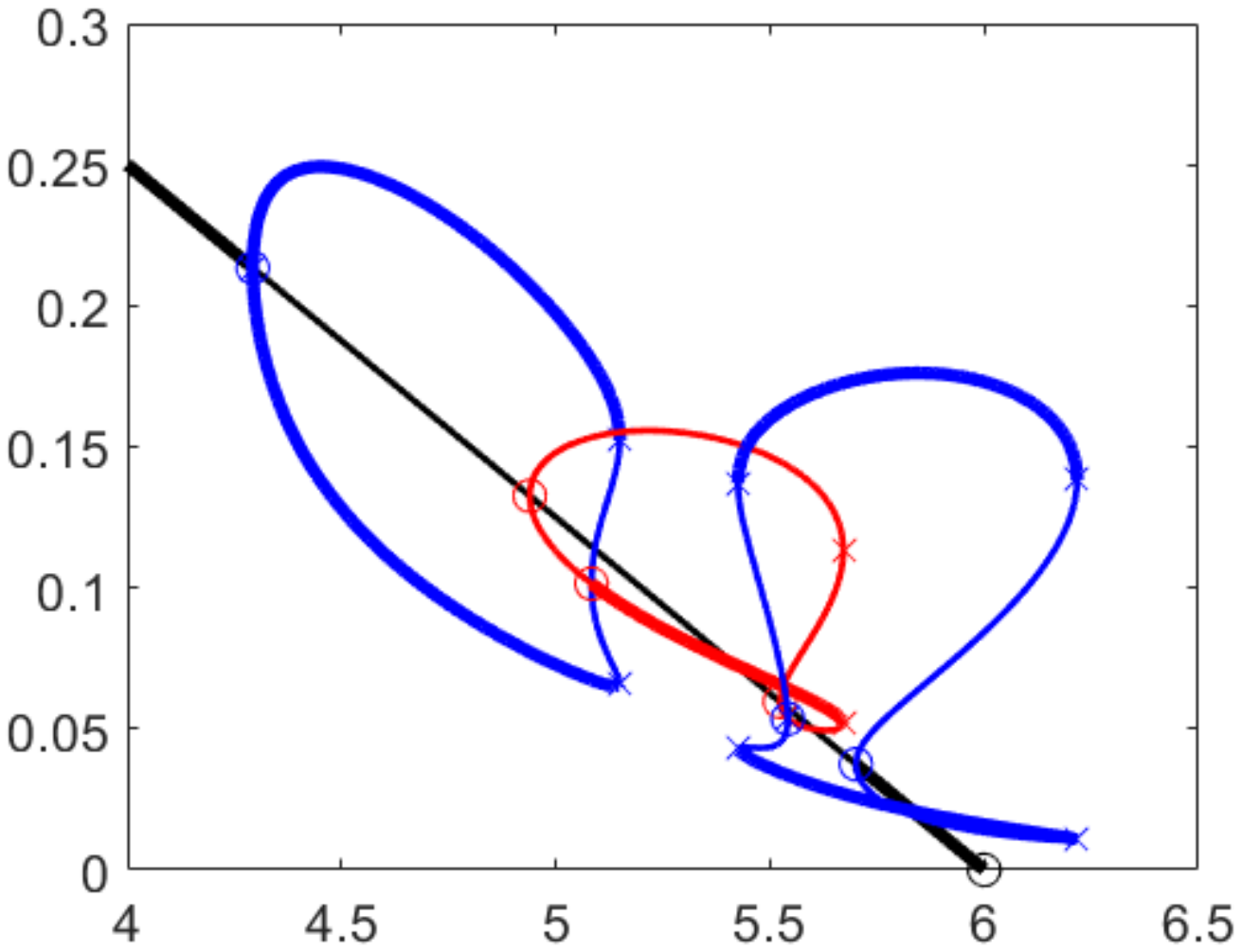}
         \put(0,30){\rotatebox{90}{$v(0)$}}
         \put(80,-2){$r_1$}
       \end{overpic}
}
\caption{Bifurcation diagrams with respect to the parameter $r_1$: \protect\subref{fast1D_r1_epsi_0p01}--\protect\subref{fast1D_r1_epsi_0p0001} correspond to the fast--slow system \eqref{fast} for different and smaller values of $\varepsilon$. We clearly observe a highly non-trivial deformation of the bifurcation diagram as $\varepsilon$ is decreased, starting from a ring and then heart-shape, we eventually have a ``broken-heart'' structure leading to two rings in the cross-diffusion limit.}
\label{diagbif_r1}
\end{figure}
\end{center}

In Figure \ref{diagbif_r1}, we show how the bifurcation diagram with respect to $r_1$ behaves when the fast time-scale parameter $\varepsilon$ decreases. On the vertical axis the value $v(0)$ is reported. Also in this case we observe convergence of the bifurcation structure of the fast-reaction system to the one of the cross-diffusion system of Figure \ref{bifdiagr1_p} (well approximated by taking $\varepsilon=10^{-4}$, corresponding to Figure \ref{fast1D_r1_epsi_0p0001}). However, we observe that the bifurcation structure is not only expanding as in the previous section, but it is qualitatively changing as $\varepsilon$ decreases. In particular, in Figure \ref{fast1D_r1_epsi_0p01} and \ref{fast1D_r1_epsi_0p005} there are only two non-homogeneous branches connecting two primary branch points on the homogeneous branch and forming a \emph{bifurcation ring}. For smaller values of $\varepsilon$ other bifurcation points (and consequently other rings) appear inside the main bifurcation ring; the inner rings expand while the outer ring folds in the middle part forming a heart shape. Then the outer ring interacts with the inner ring and separates into two rings, giving rise to the bifurcation structure of the cross-diffusion system, i.e., the heart-shaped structure breaks.

\begin{figure}
\centering
\begin{overpic}[width=0.5\textwidth,tics=10,trim=112 240 120 250,clip]{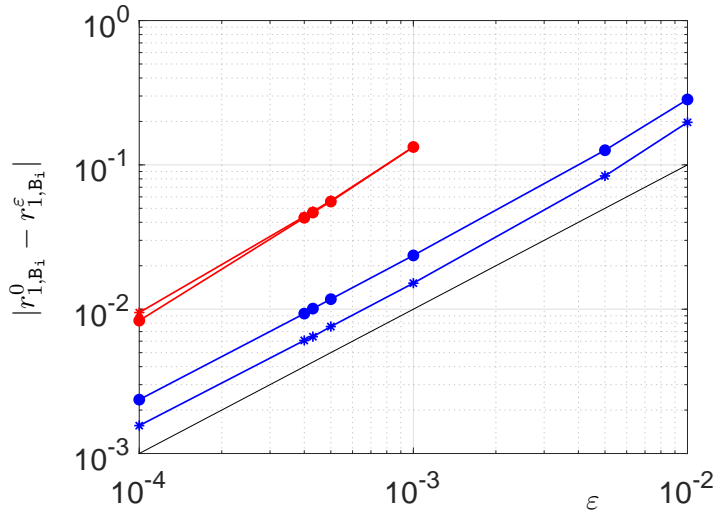}
\put(-10,30){\rotatebox{90}{$|r^0_{1,\mathtt{B_i}}-r^{\varepsilon}_{1,\mathtt{B_i}}|$}}
\put(80,0){$\varepsilon$}
\end{overpic}
\caption{Convergence of the bifurcation points in loglog scale. We report on the horizontal axis the values of $\varepsilon$, and on the vertical axis the difference between the bifurcation value of the fast--slow system $r^\varepsilon_{1,\mathtt{B_i}},$ and the corresponding one of the cross-diffusion system $r^0_{1,\mathtt{B_i}}\, i=1,2,3,4$ (blue and red refer to Figure \ref{diagbif_r1}, dots and stars denote the first and the second bifurcation point, respectively, for each color). The black line corresponds the order of convergence one.}
\label{1D_fast_r1_conv_loglog}
\end{figure}

Also in this case we quantitatively show the convergence of the first bifurcation points on the homogeneous branch in Figure \ref{1D_fast_r1_conv_loglog}. Although the bifurcation structure is qualitatively changing, the convergence rate is comparable to Figure \ref{1D_fast_d_conv_loglog}. This clearly shows that locally one can expect a convergence rate near the homogeneous branch. In fact, this leads one to the conjecture that even the global diagram could be captured asymptotically by a convergence rate towards the singular limit but since we have only captured part of the full diagram, this is hard to validate completely numerically. Furthermore, many different convergence metrics are conceivable if one moves beyond single points so we leave this conjecture as a topic for future work. 

\FloatBarrier


\section{Numerical continuation on a 2D rectangular domain}
\label{sec:nc2D}

In this section we consider a 2D rectangular domain with edges of length $L_x=1$ and $L_y=4$. This choice reduces the presence of multiple branch points. As in the previous section, we set the standard diffusion coefficient $d$ as our main bifurcation parameter. The bifurcation analysis close to the homogeneous branch can be performed (see Appendix~\ref{app:TI}). It allows to compute the values of the parameter $d$ at the bifurcation points on the homogeneous branch, which can be compared to the values numerically obtained. It also predicts the shape of the steady state for each bifurcation point close to the homogeneous branch, by looking at the eigenvalues of the Laplacian at the bifurcation point and the associated eigenfunctions~\cite{Kielhoefer,KuehnBook1}. 

\begin{figure}[!ht]
\centering
    \subfloat{
       \begin{overpic}[width=0.7\textwidth,tics=10,trim=80 280 80 250,clip]{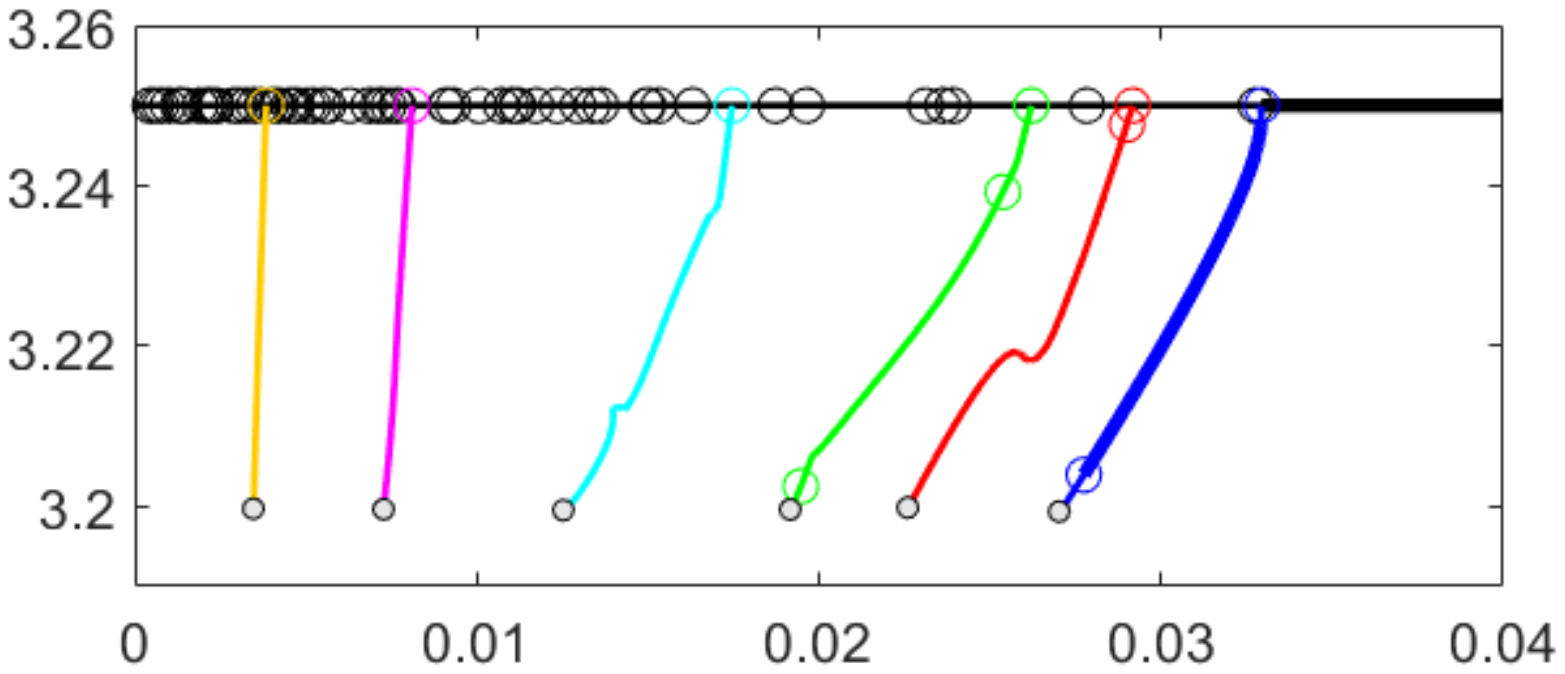}%
         \put(-7,25){\rotatebox{90}{$||u||_{L^2}$}}
         \put(85,0){$d$}
         \put(63,10){(a)}
         \put(53,10){(b)}
         \put(46,10){(c)}
         \put(31,10){(d)}
         \put(20,10){(e)}
         \put(12,10){(f)}
       \end{overpic} }\\
       \setcounter{subfigure}{0}
     \subfloat[\label{sol2Da}]{
       \begin{overpic}[width=0.3\textwidth,tics=10,trim=90 240 90 250,clip]{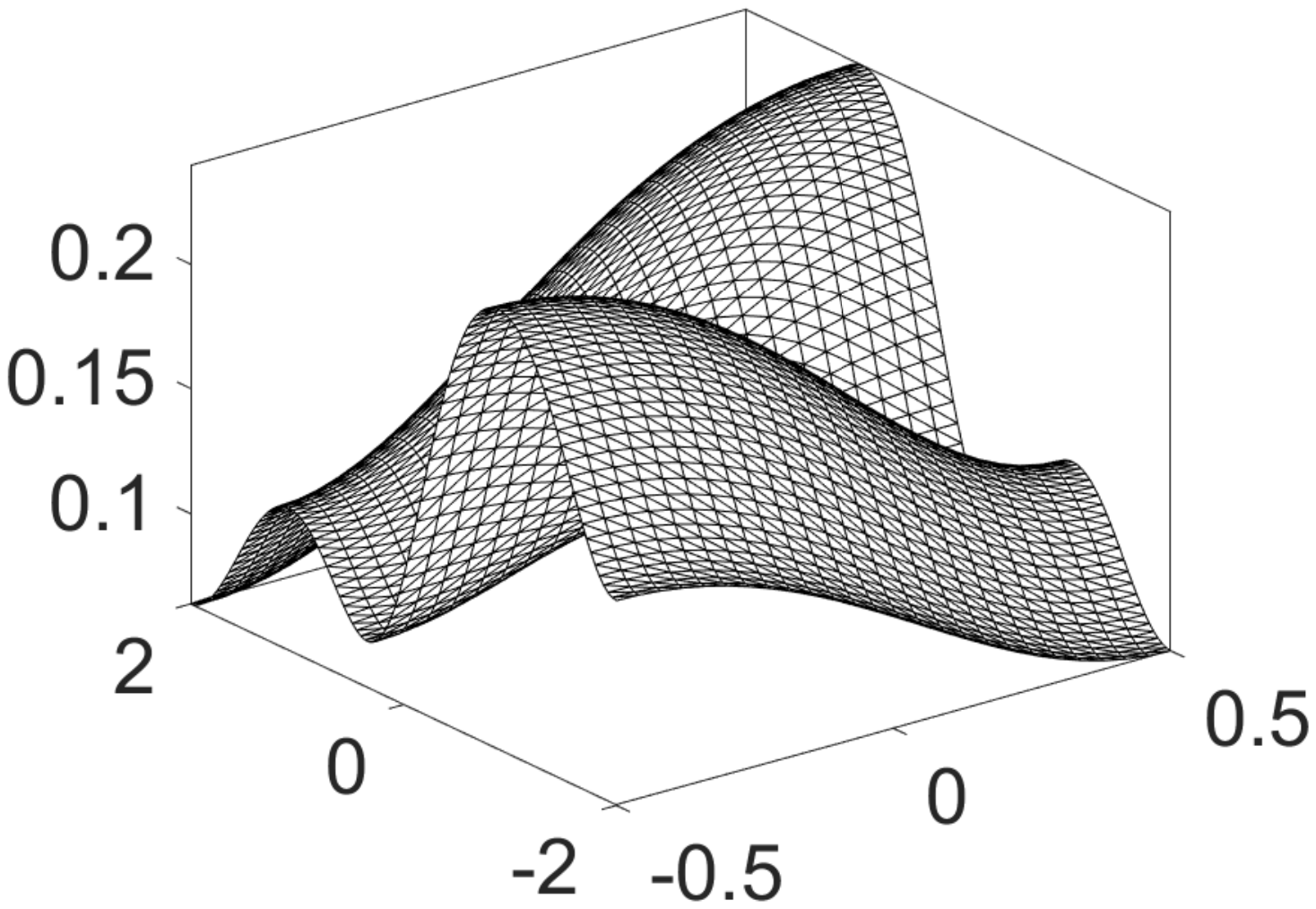}
         \put(80,1){$x$}
         \put(15,5){$y$}
         \end{overpic}}
	\subfloat[\label{sol2Db}]{
       \begin{overpic}[width=0.3\textwidth,tics=10,trim=90 240 90 250,clip]{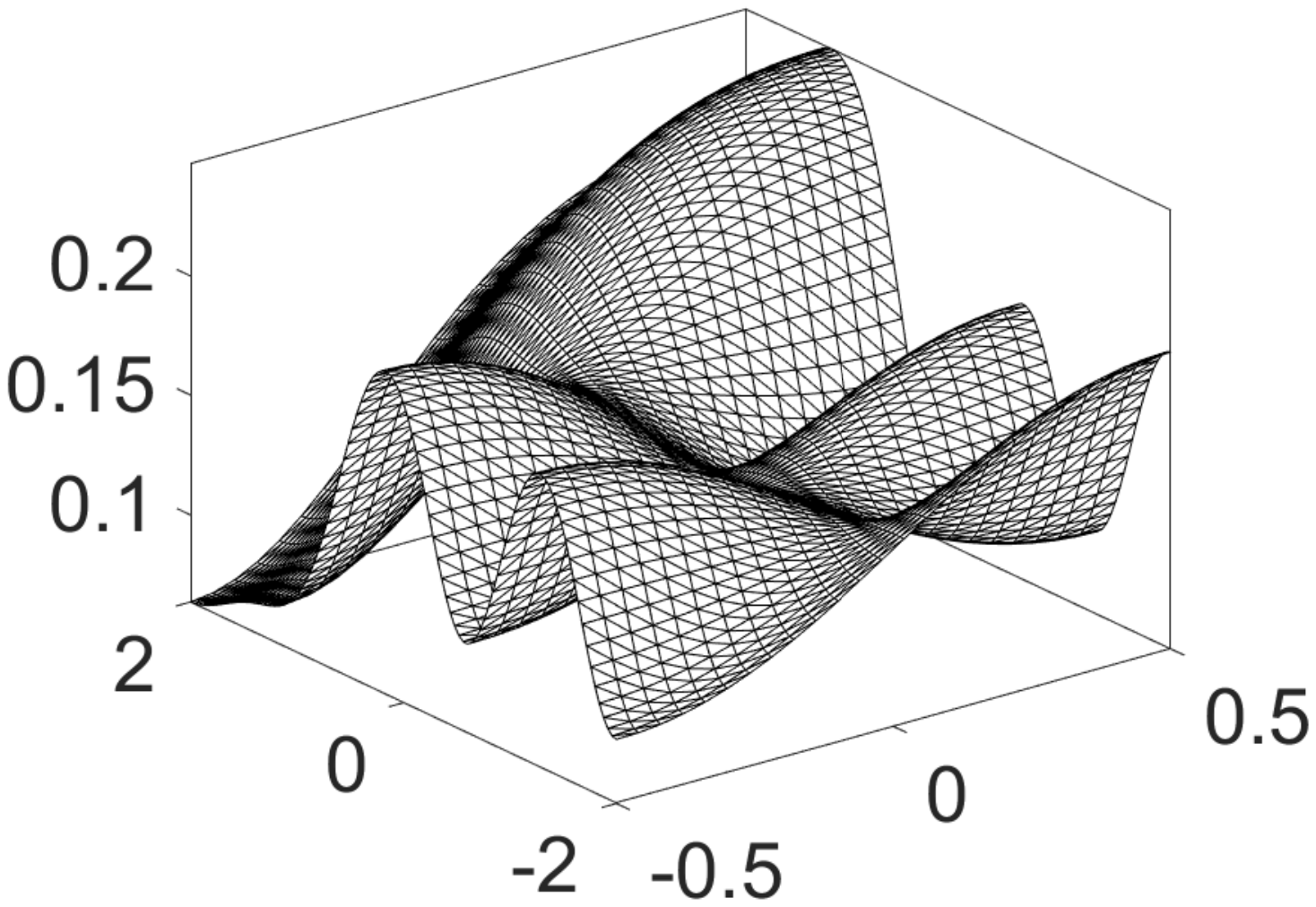}
                \put(80,1){$x$}
         \put(15,5){$y$}
         \end{overpic}}
	\subfloat[\label{sol2Dc}]{
       \begin{overpic}[width=0.3\textwidth,tics=10,trim=90 240 90 250,clip]{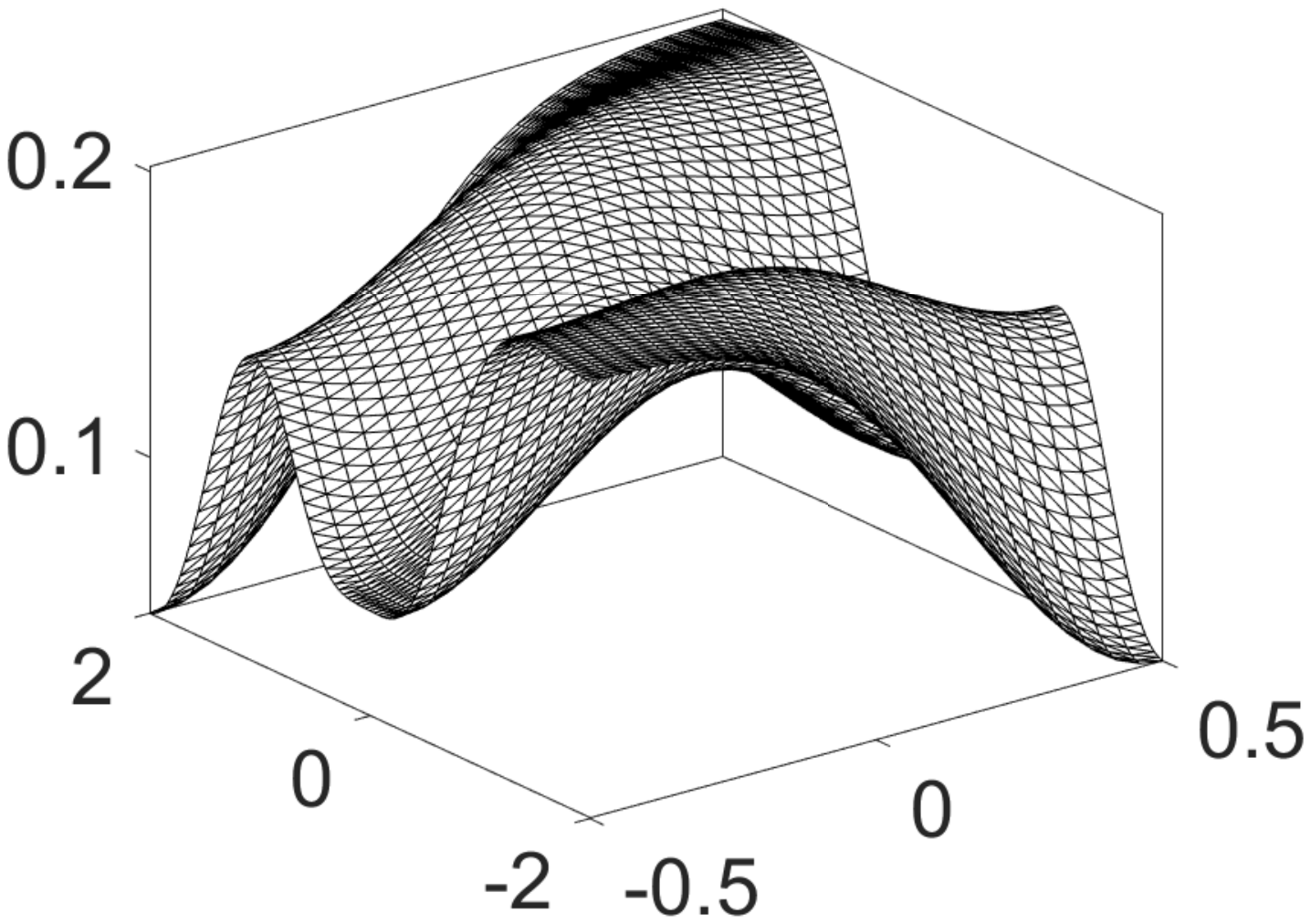}
                \put(80,1){$x$}
         \put(15,5){$y$}
         \end{overpic}}\\
	\subfloat[\label{sol2Dd}]{
       \begin{overpic}[width=0.3\textwidth,tics=10,trim=90 240 90 250,clip]{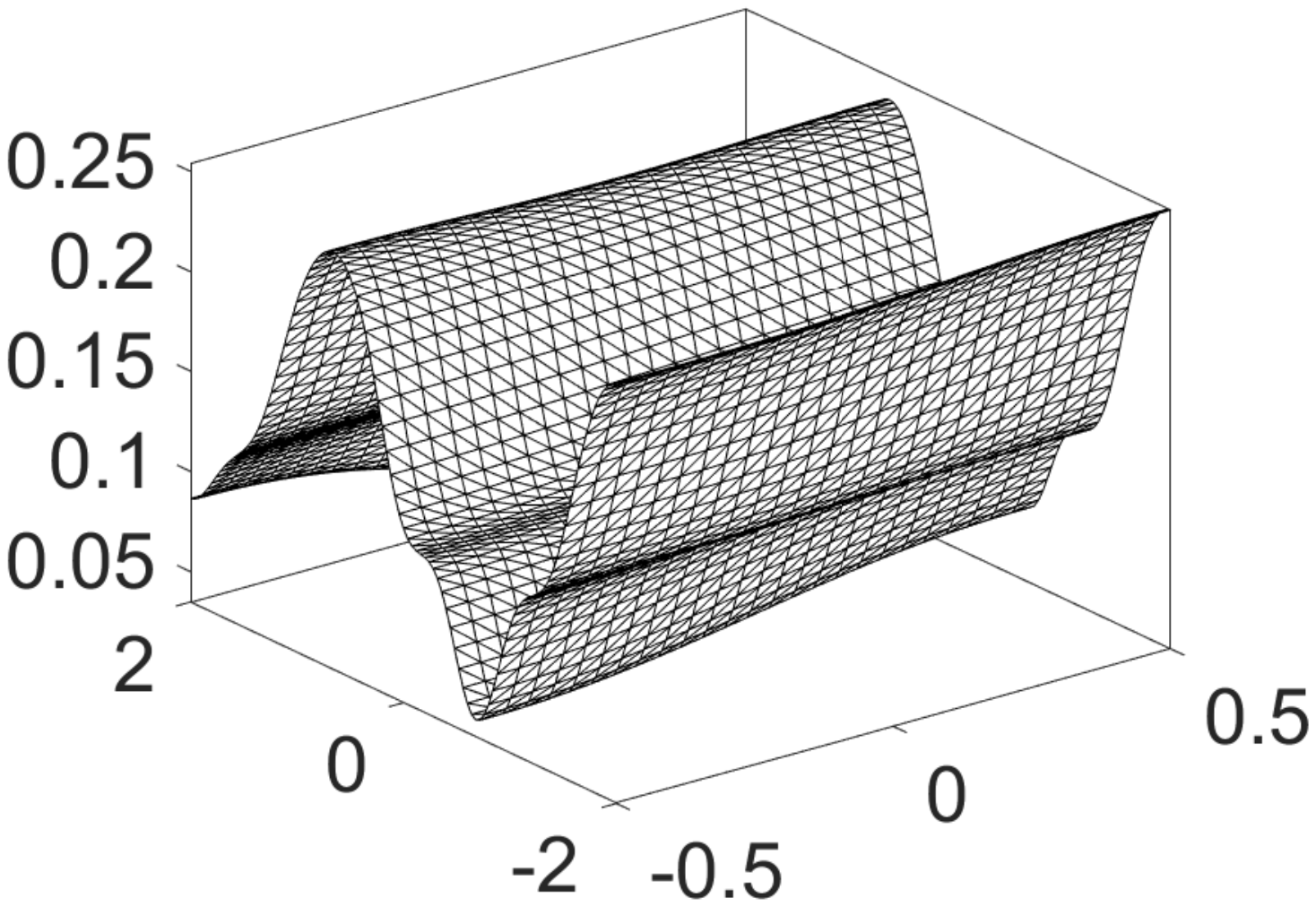}
                \put(80,1){$x$}
         \put(15,5){$y$}
         \end{overpic}}
	\subfloat[\label{sol2De}]{
       \begin{overpic}[width=0.3\textwidth,tics=10,trim=90 240 90 250,clip]{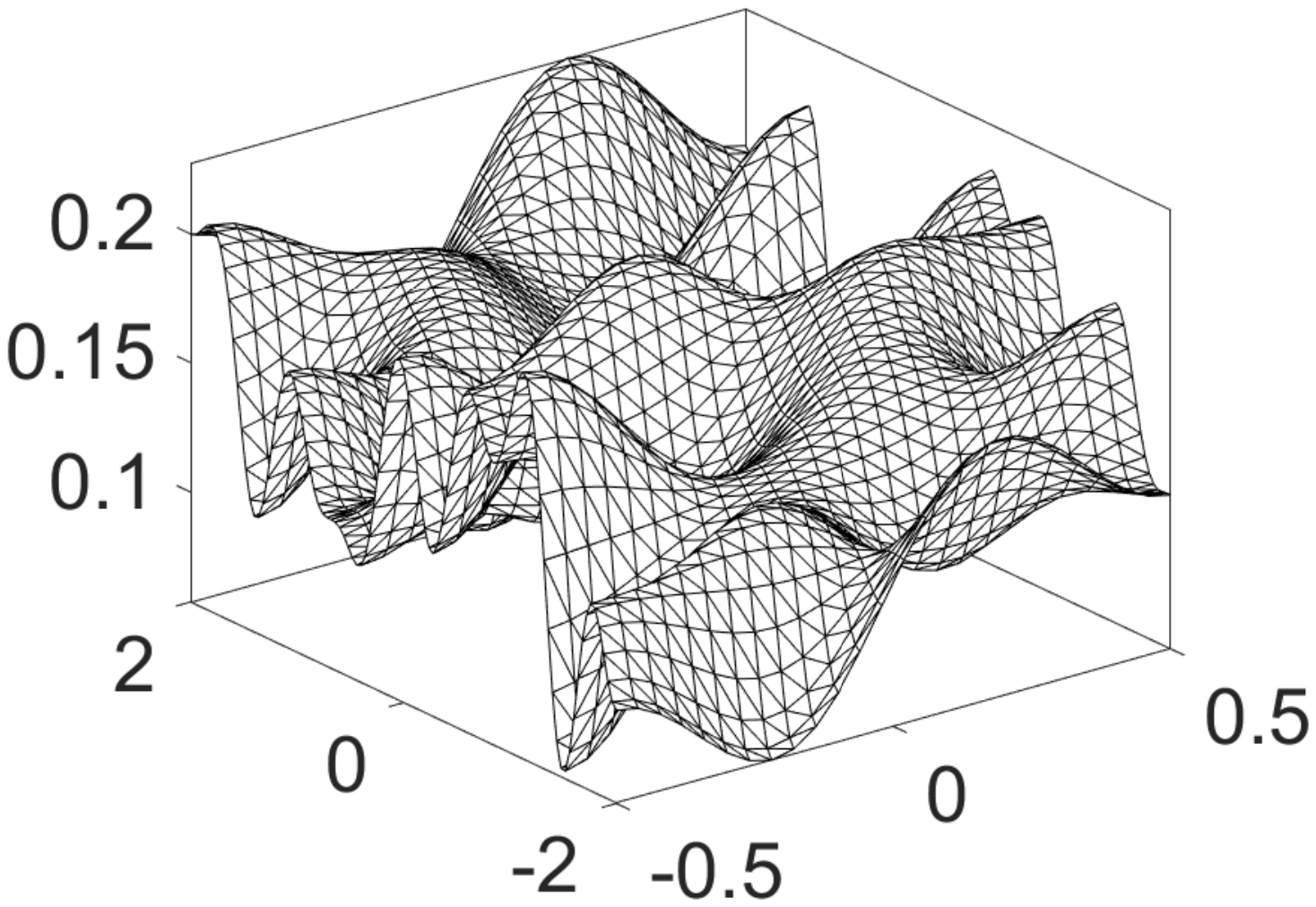}
                \put(80,1){$x$}
         \put(15,5){$y$}
         \end{overpic}}
         \subfloat[\label{sol2Df}]{
       \begin{overpic}[width=0.3\textwidth,tics=10,trim=90 240 90 250,clip]{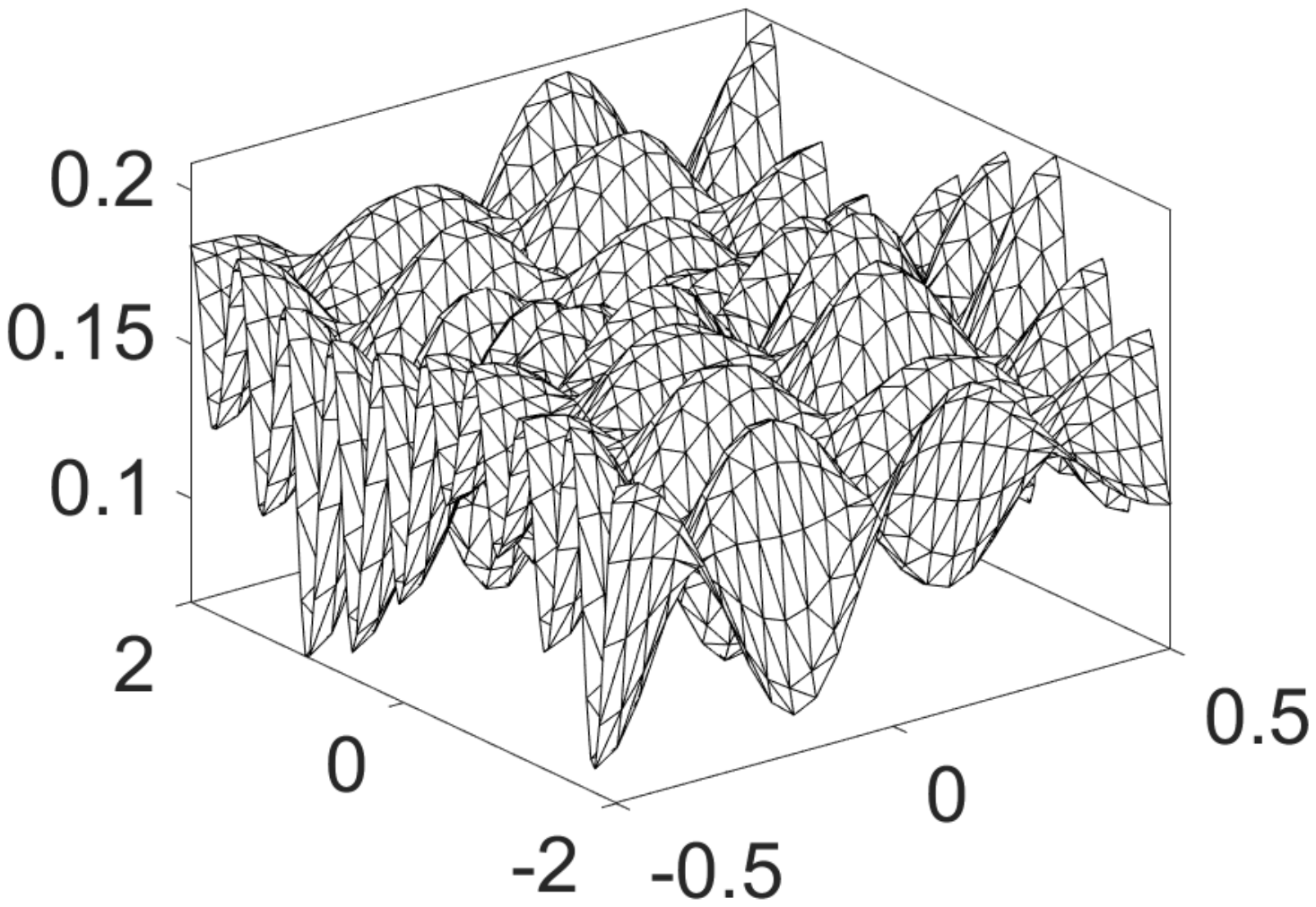}
                \put(80,1){$x$}
         \put(15,5){$y$}
         \end{overpic}}
\caption{Bifurcation diagram and different solution types along the branches. Upper panel: partial bifurcation diagram relative to a rectangular domain $L_x=1,\; L_y=4$. Lower panel: \protect\subref{sol2Da}--\protect\subref{sol2Df} solutions on different branches (species $u$).}
\label{2D_DiagBifCross_sol}
\end{figure}

In Figure~\ref{2D_DiagBifCross_sol} we show part of the bifurcation diagram close to the homogeneous branch with respect to the $L^2$-norm of the species $u$. As in the 1D case, for decreasing values of $d$ the homogeneous solution destabilizes at the first bifurcation point, and stable non-homogeneous steady states appear (blue branch). For smaller values of $d$ other bifurcation points occur. We show in the figure just some of the successive branches and we report the corresponding solution at the gray points (Figures~\ref{sol2Da}--\ref{sol2Df}), which turn out to be unstable. Note that certain cross-sections of the 2D-solutions have a similar shape as the steady states of 1D steady states.

\begin{figure}[!ht]
\centering
\begin{multicols}{2}
\vspace*{\fill}   
\hspace{-0.4cm}  
\subfloat{
\begin{overpic}[width=0.43\textwidth,tics=10,trim=100 240 110 260,clip]{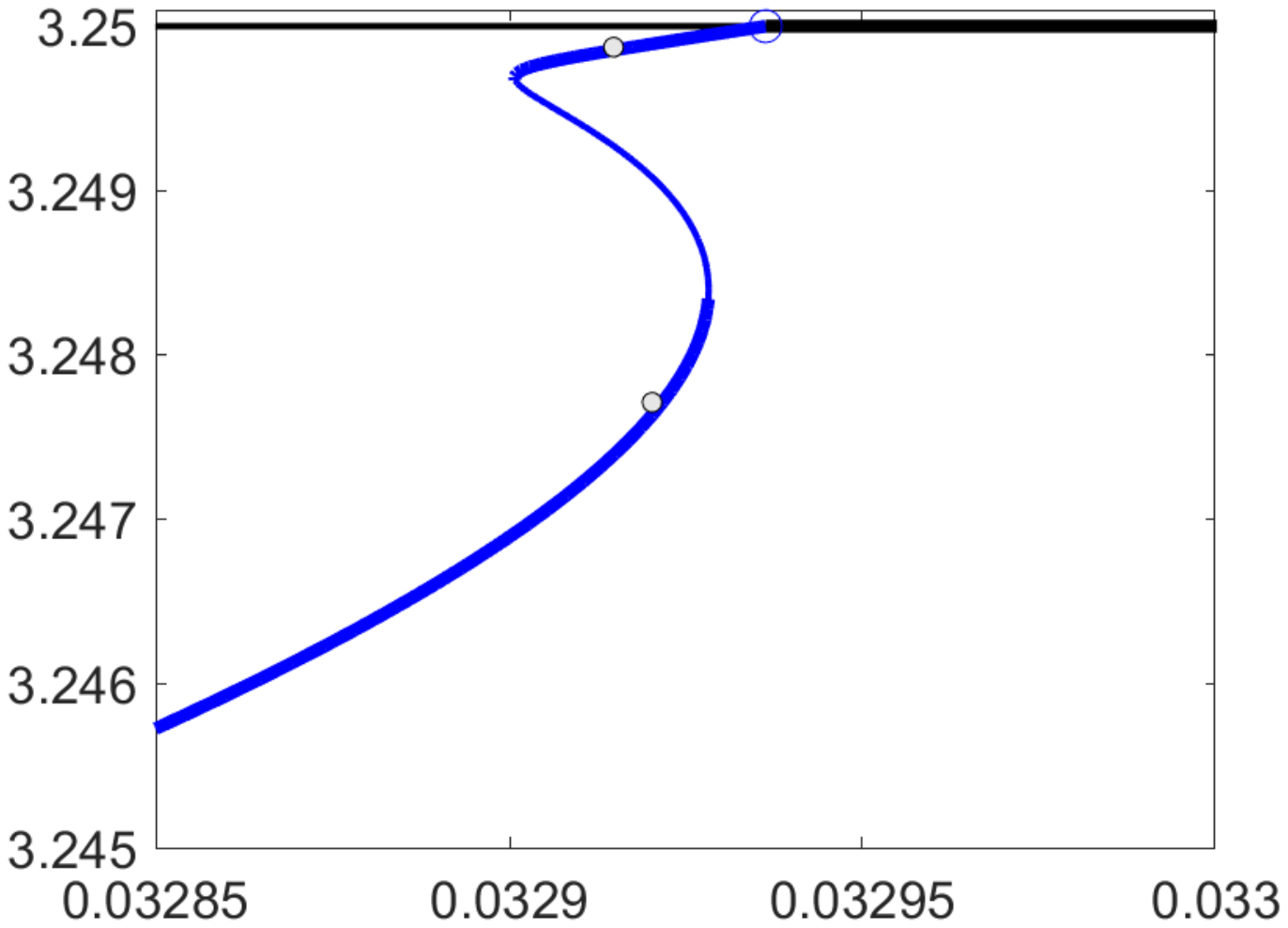}%
         \put(-5,30){\rotatebox{90}{$||u||_{L^2}$}}    
         \put(50,63){(a),(b)}
         \put(50,35){(c),(d)} 
\end{overpic} }
\vspace*{\fill}  
\newpage    
        
\setcounter{subfigure}{0}
\centering
\hspace{-1.3cm} 
\subfloat[\label{sol2Da_3}]{
\begin{overpic}[width=0.28\textwidth,tics=10,trim=90 240 90 250,clip]{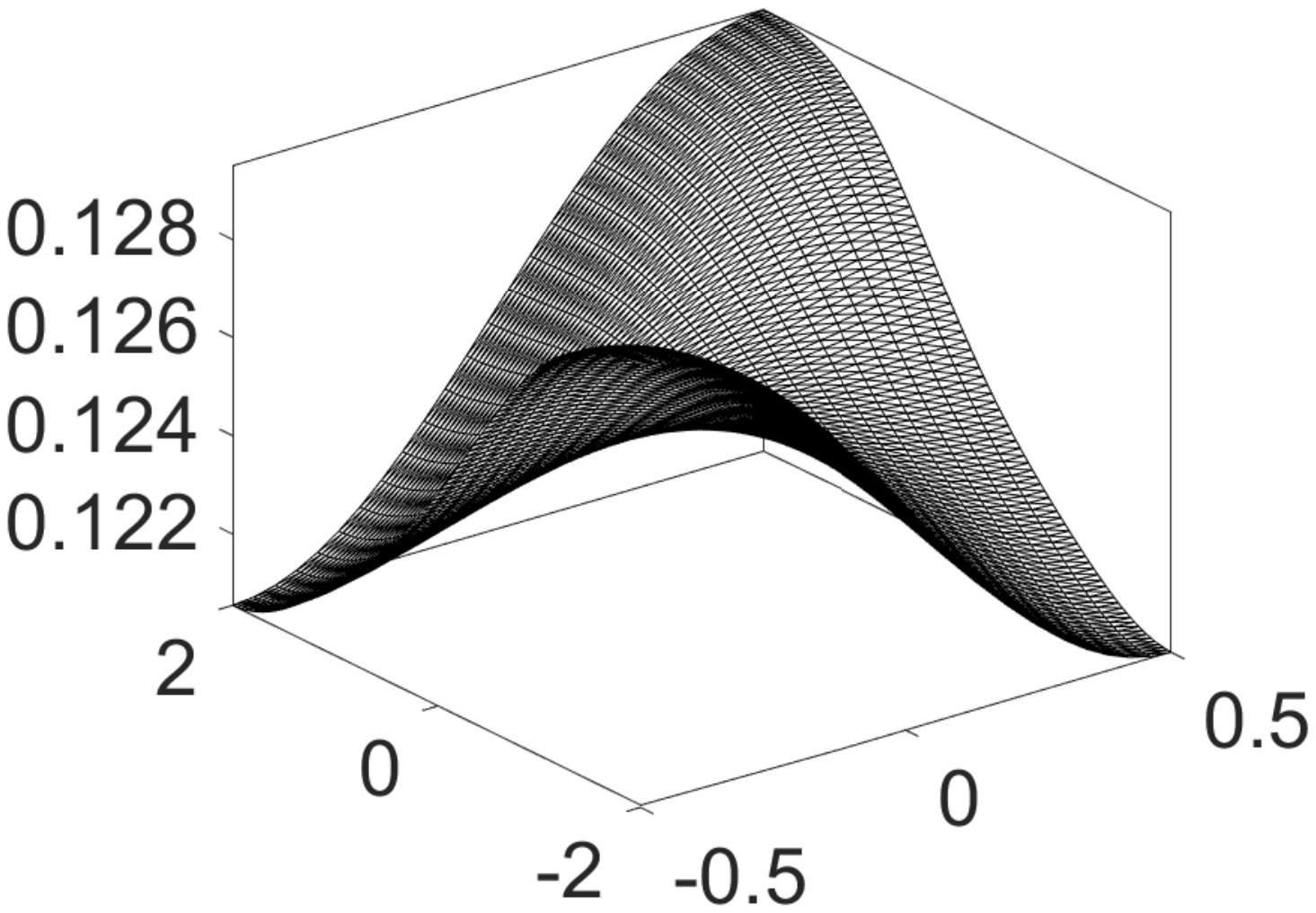}
         \put(80,1){$x$}
         \put(15,5){$y$}
\end{overpic}}
\hspace{-0.4cm} 
\subfloat[\label{sol2Da_3u}]{
\begin{overpic}[width=0.28\textwidth,tics=10,trim=90 240 90 250,clip]{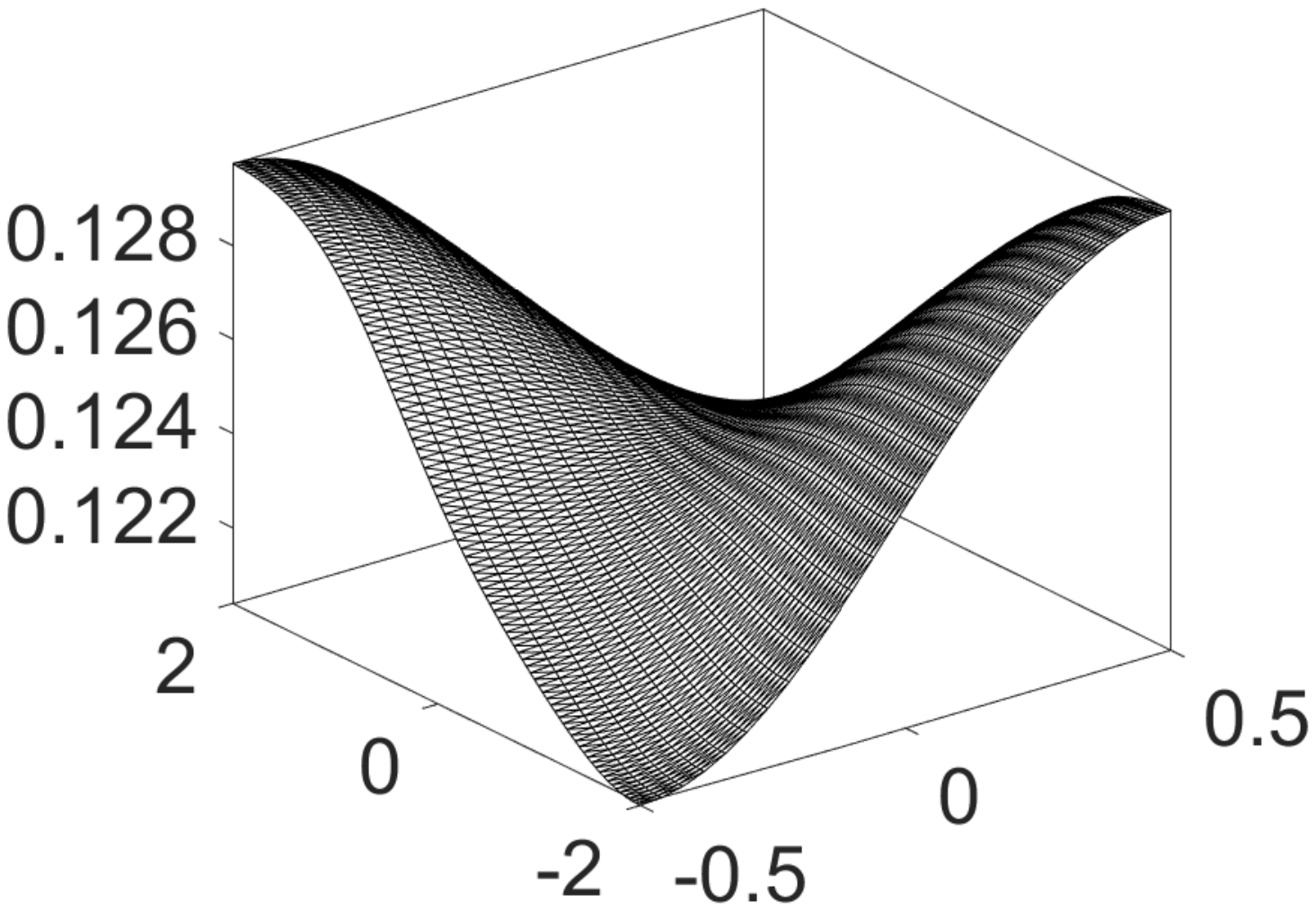}
     \put(80,1){$x$}
	\put(15,5){$y$}
\end{overpic}}\\    
\hspace{-1.3cm} 
\subfloat[\label{sol2Db_3}]{
\begin{overpic}[width=0.28\textwidth,tics=10,trim=90 240 90 250,clip]{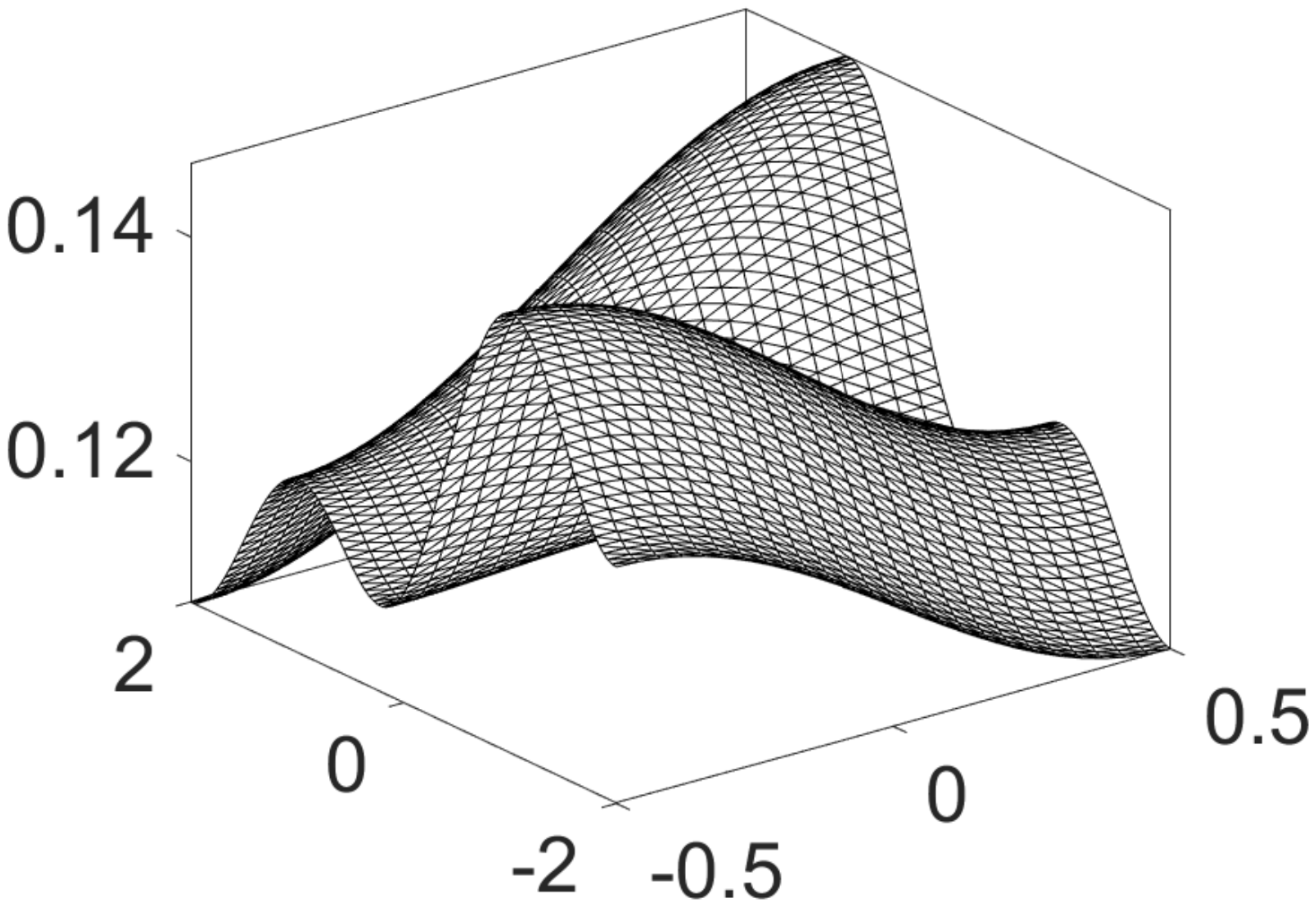}
         \put(80,1){$x$}
         \put(15,5){$y$}
\end{overpic}}
\hspace{-0.4cm} 
\subfloat[\label{sol2Db_3u}]{
\begin{overpic}[width=0.28\textwidth,tics=10,trim=90 240 90 250,clip]{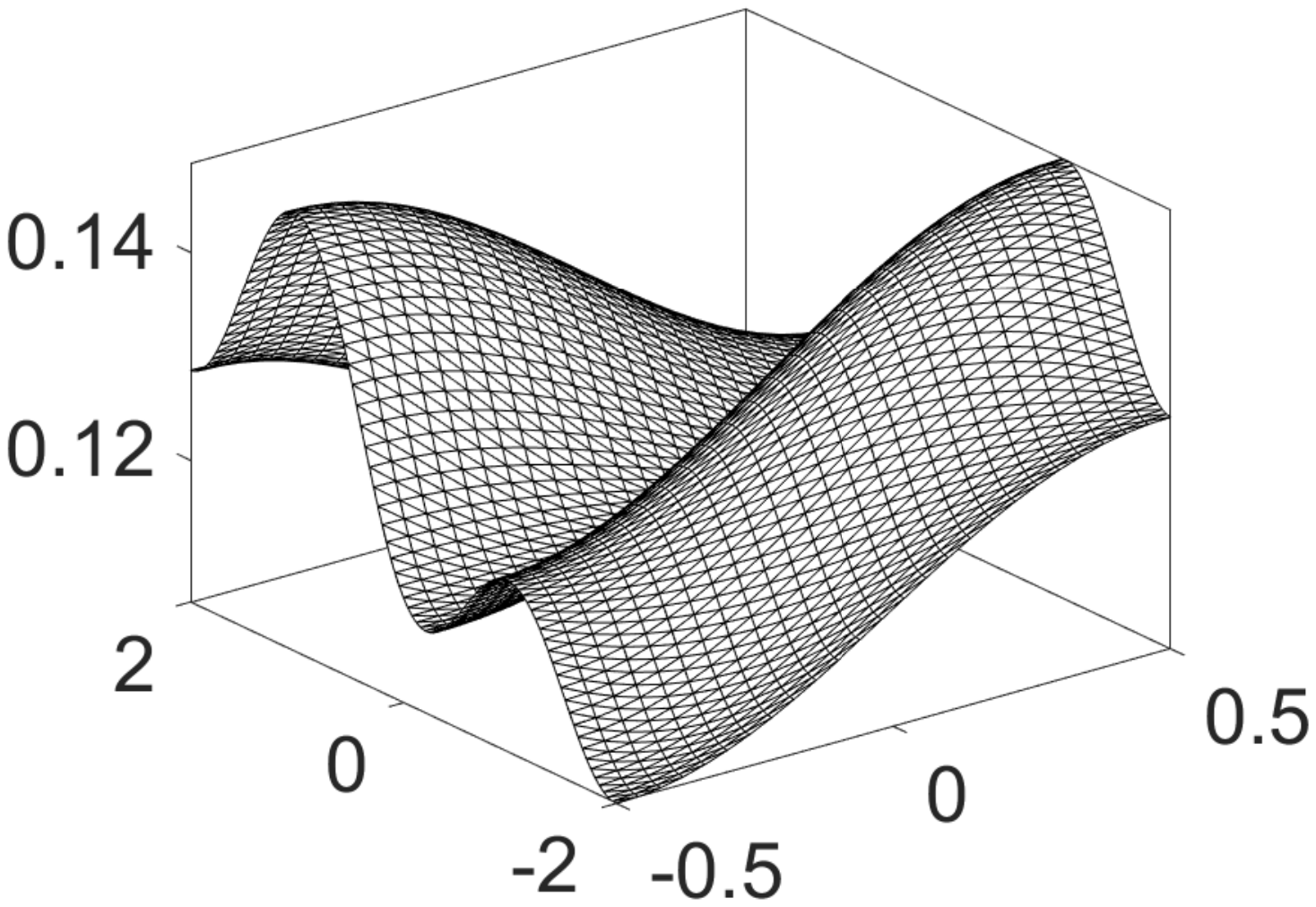}
         \put(80,1){$x$}
         \put(15,5){$y$}
\end{overpic}}
\end{multicols}
\caption{Bifurcation diagram and multi-stable solutions along the first branch. Left panel: zoom in of Figure~\ref{2D_DiagBifCross_sol} close to the first bifurcation point. Right panel: \protect\subref{sol2Da_3}--\protect\subref{sol2Db_3u} solutions (species $u$) corresponding to points on the first bifurcating branch marked with the gray dots.}
\label{2D_DiagBifCross_sol3}
\end{figure}

Different from the 1D case, in Figure~\ref{2D_DiagBifCross_sol3} the enlargement of the initial part of the first (blue) branch is reported. It shows that the first bifurcation point is supercritical, but then two successive fold bifurcations (indicated with a cross) appear leading to multi-stability of solutions: in a small range of the bifurcation parameter the system admits four (two symmetric) stable non-homogeneous steady states. Their shape is reported in Figures~\ref{sol2Da_3}--\ref{sol2Db_3u}. Close to the homogeneous branch (Figures \ref{sol2Da_3} and~\ref{sol2Da_3u}) the steady states have half a bump on both edges (in agreement with the eigenvalues corresponding to the first bifurcation point, see Appendix~\ref{app:TI}), while going further along the branch the shape modifies. 

\begin{figure}[!ht]
\centering
    \subfloat{
       \begin{overpic}[width=0.7\textwidth,tics=10,trim=50 250 50 250,clip]{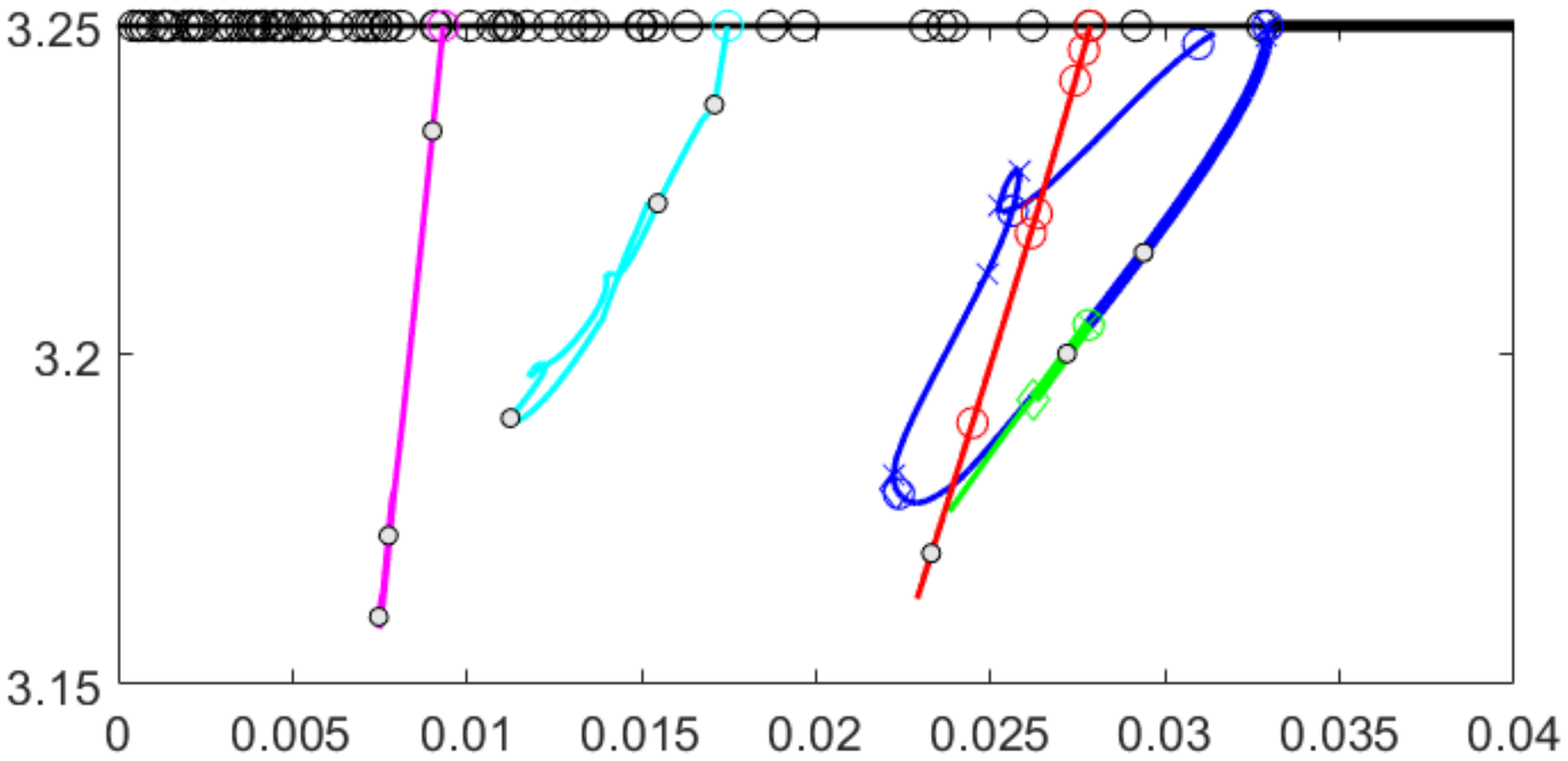}%
         \put(-2,25){\rotatebox{90}{$||u||_{L^2}$}}
         \put(98,9){$d$}
         \put(75,35){(a)}
         \put(70,28){(b)}
         \put(60,12){(c)}
         \put(47,45){(d)}
         \put(42,35){(e)}
         \put(33,20){(f)}
         \put(22,40){(g)}
         \put(19,17){(h)}
         \put(19,12){(i)}         
       \end{overpic} }\\
       \setcounter{subfigure}{0}
     \subfloat[\label{sol2Da_2}]{
       \begin{overpic}[width=0.3\textwidth,tics=10,trim=90 240 90 250,clip]{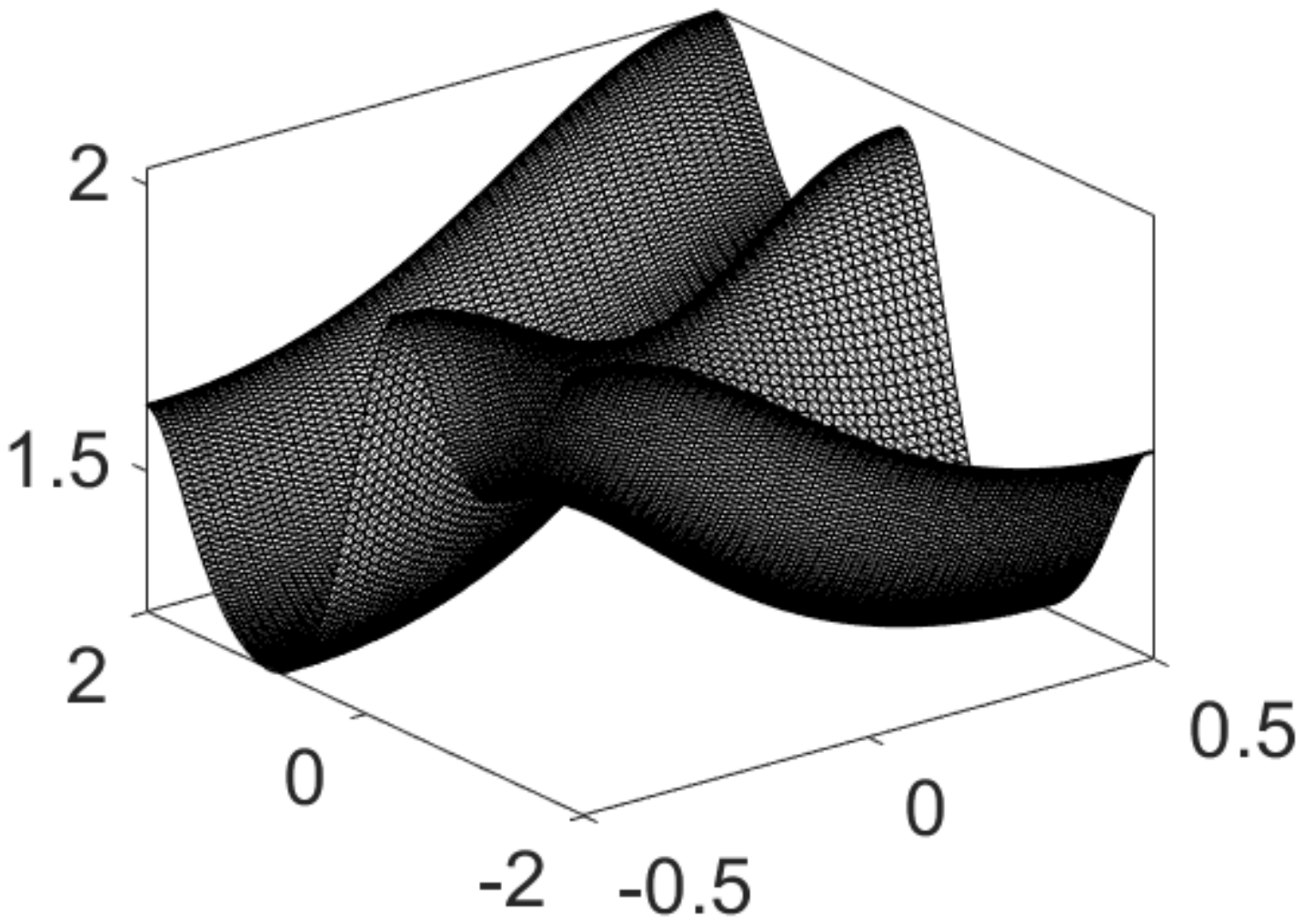}
         \put(80,1){$x$}
         \put(15,5){$y$}
         \end{overpic}}
	\subfloat[\label{sol2Db_2}]{
       \begin{overpic}[width=0.3\textwidth,tics=10,trim=90 240 90 250,clip]{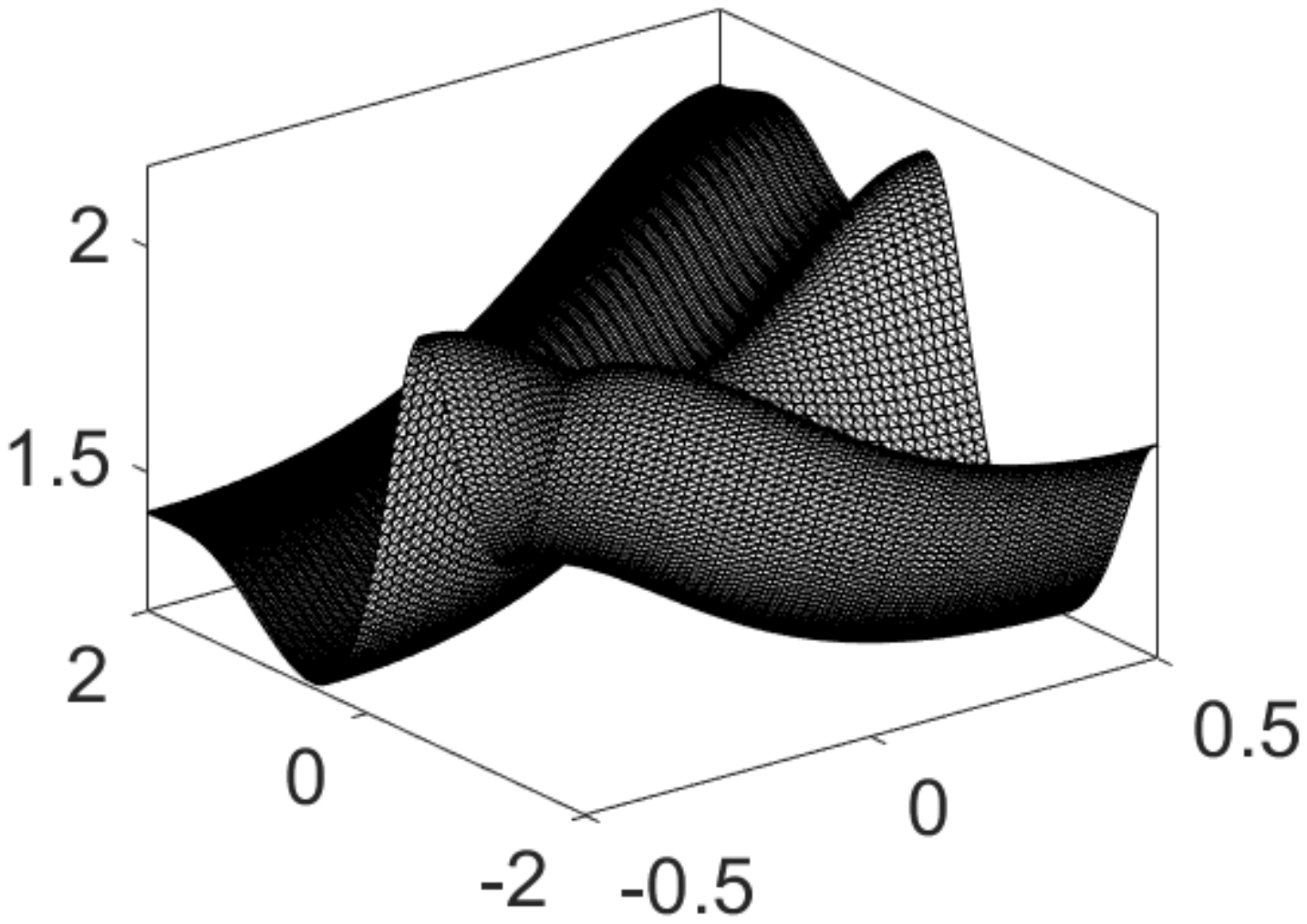}
                \put(80,1){$x$}
         \put(15,5){$y$}
         \end{overpic}}
	\subfloat[\label{sol2Dc_2}]{
       \begin{overpic}[width=0.3\textwidth,tics=10,trim=90 240 90 250,clip]{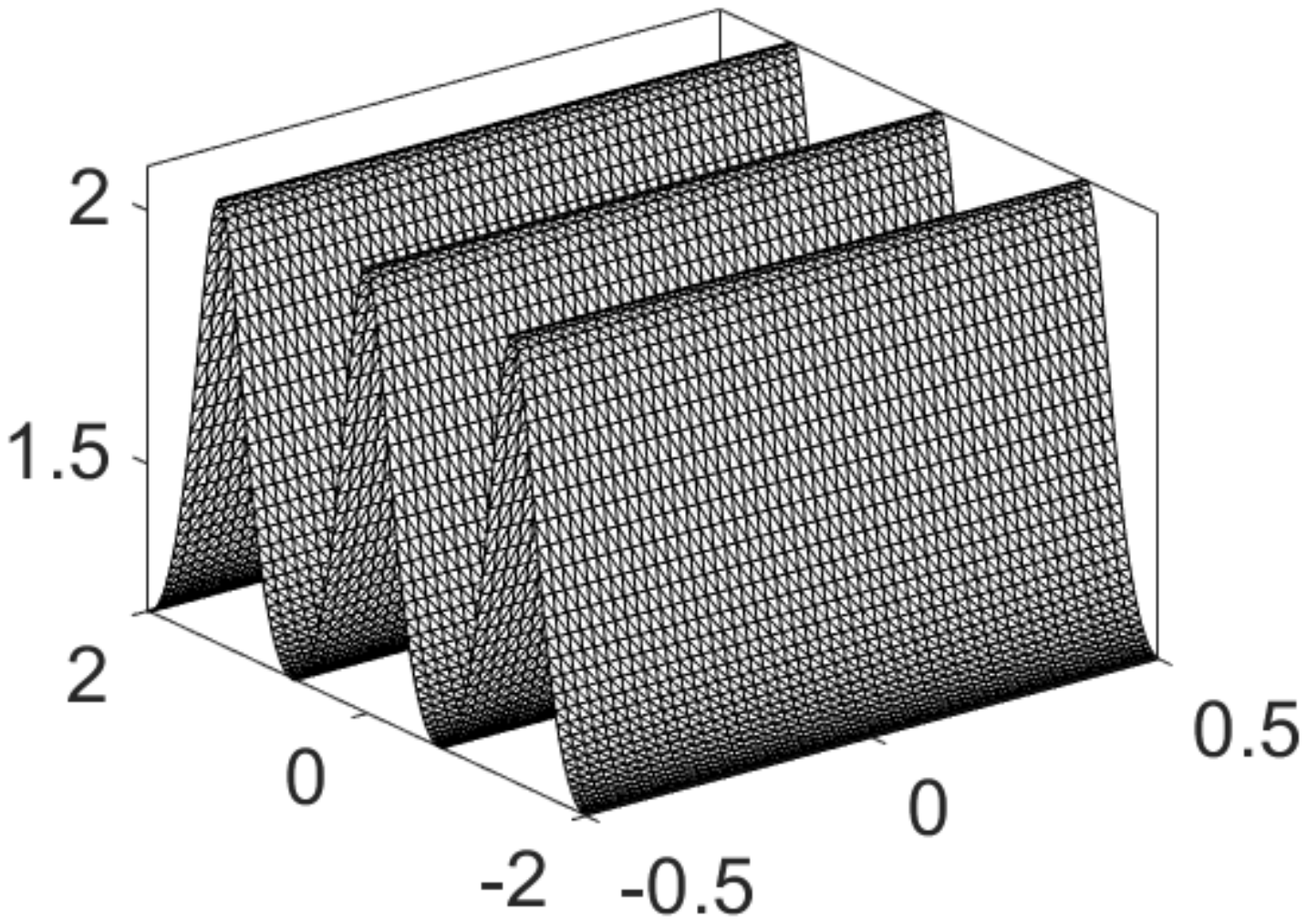}
                \put(80,1){$x$}
         \put(15,5){$y$}
         \end{overpic}}\\
	\subfloat[\label{sol2Dd_2}]{
       \begin{overpic}[width=0.3\textwidth,tics=10,trim=90 240 90 250,clip]{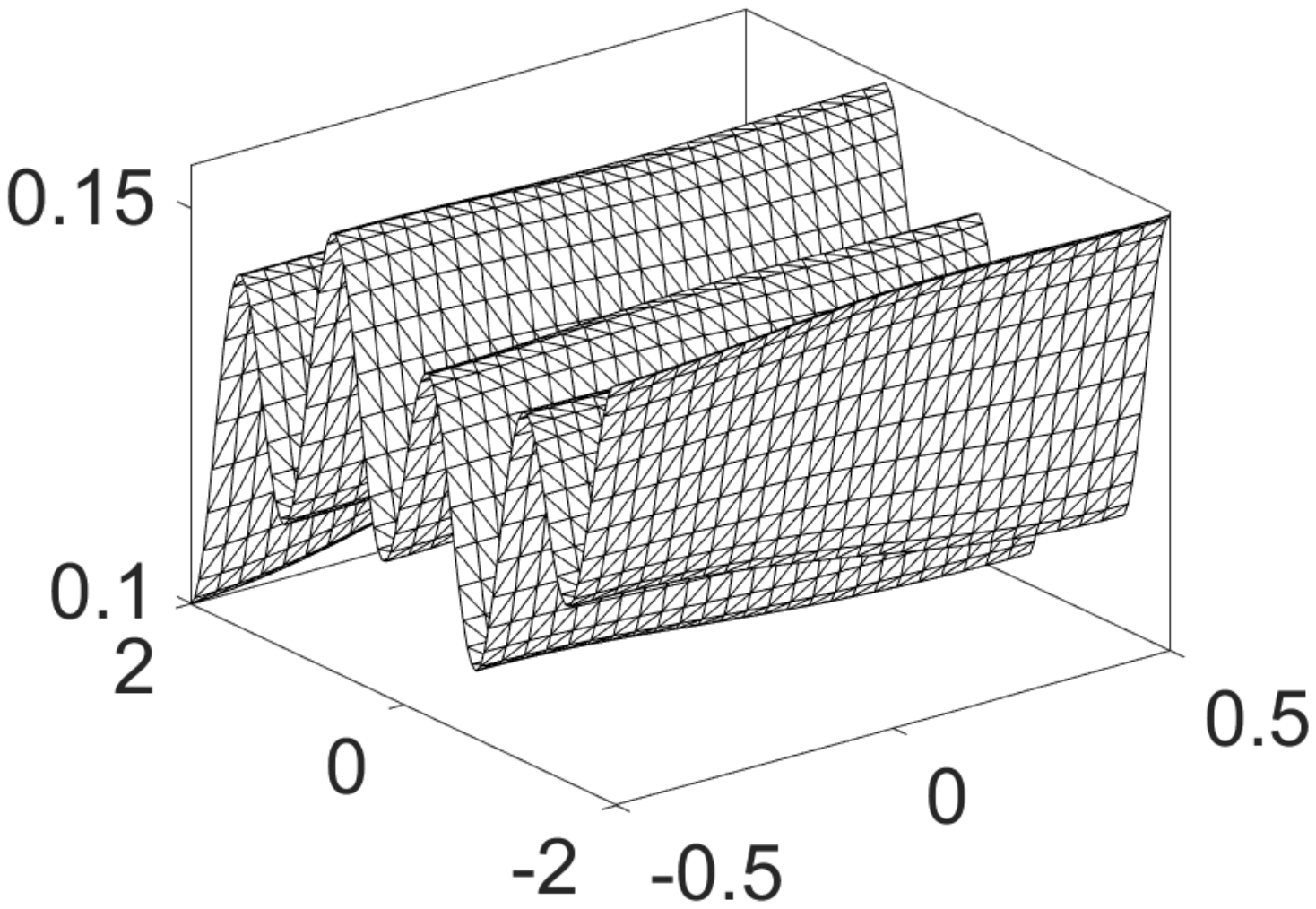}
                \put(80,1){$x$}
         \put(15,5){$y$}
         \end{overpic}}
	\subfloat[\label{sol2De_2}]{
       \begin{overpic}[width=0.3\textwidth,tics=10,trim=90 240 90 250,clip]{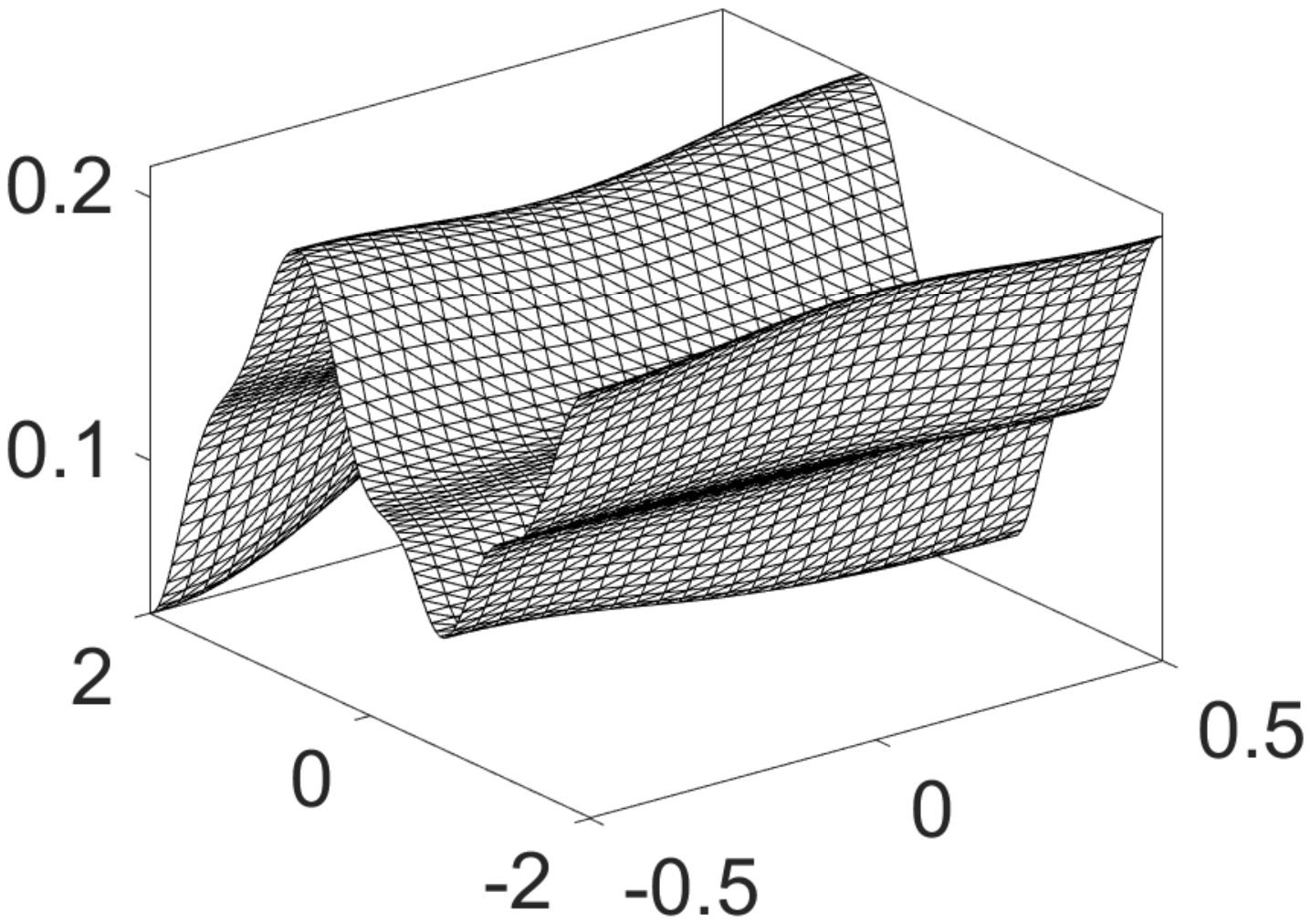}
                \put(80,1){$x$}
         \put(15,5){$y$}
         \end{overpic}}
         \subfloat[\label{sol2Df_2}]{
       \begin{overpic}[width=0.3\textwidth,tics=10,trim=90 240 90 250,clip]{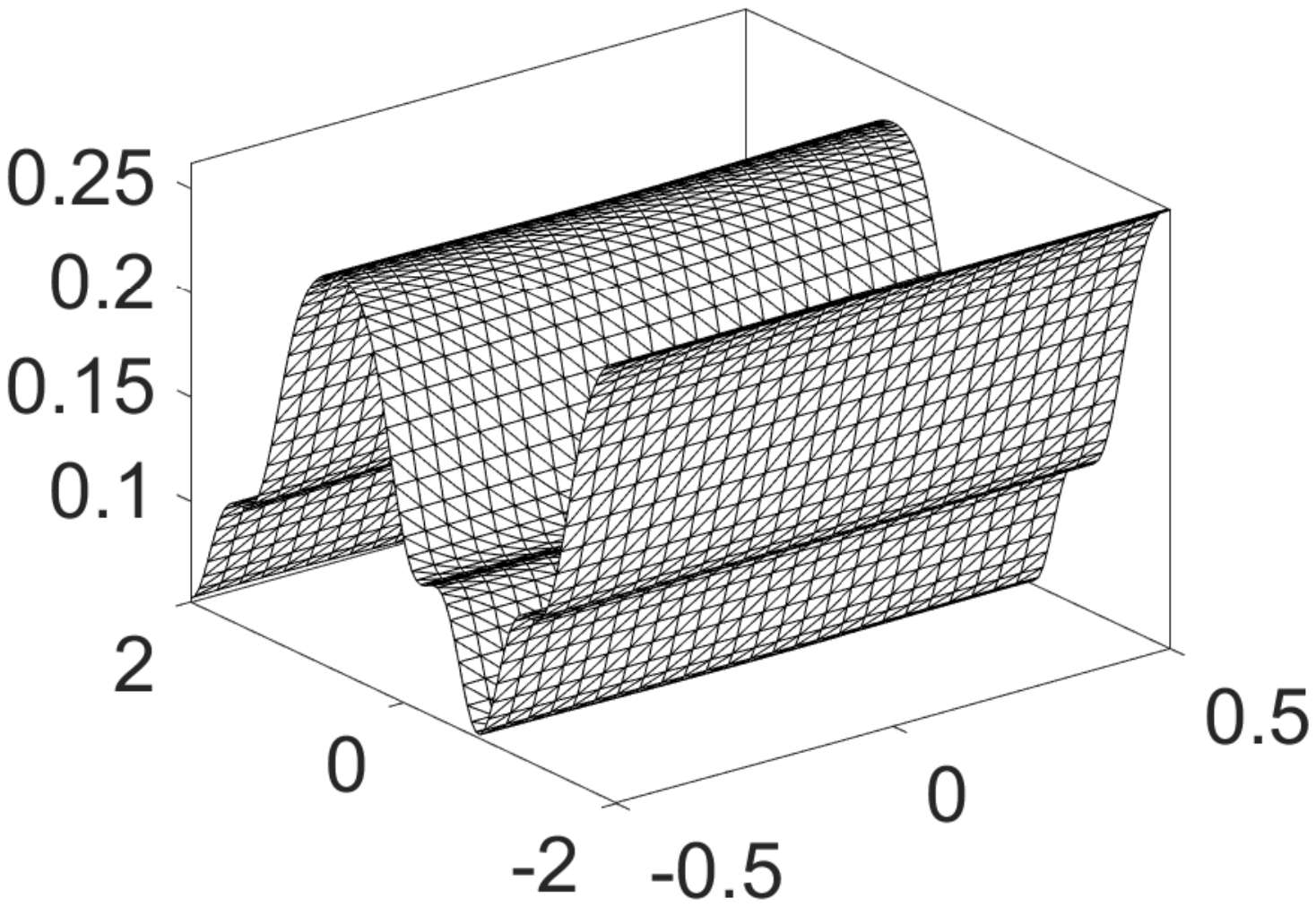}
                \put(80,1){$x$}
         \put(15,5){$y$}
         \end{overpic}}\\
          \subfloat[\label{sol2Dg_2}]{
        \begin{overpic}[width=0.3\textwidth,tics=10,trim=90 240 90 250,clip]{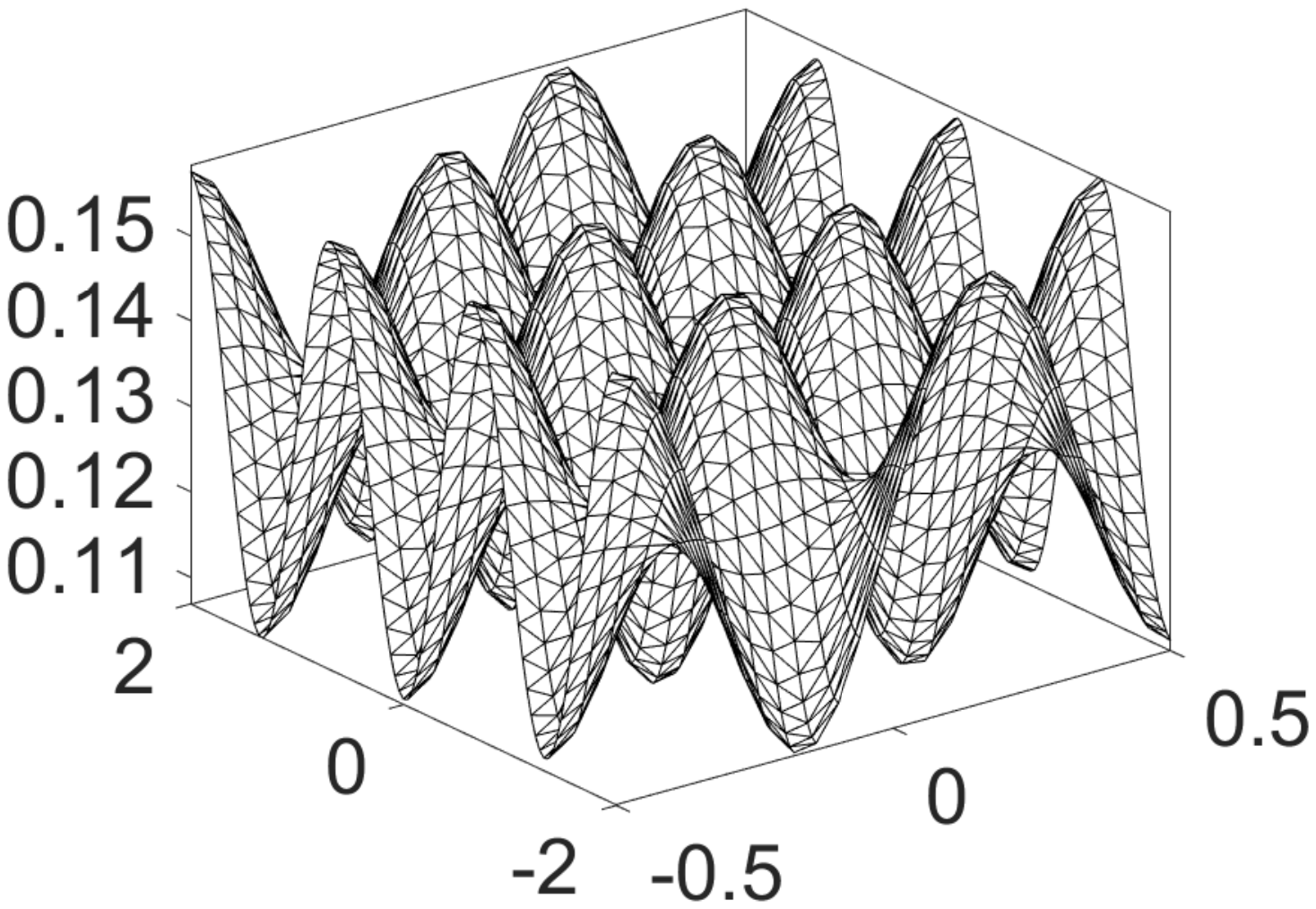}
                \put(80,1){$x$}
         \put(15,5){$y$}
         \end{overpic}}
	\subfloat[\label{sol2Dh_2}]{
       \begin{overpic}[width=0.3\textwidth,tics=10,trim=90 240 90 250,clip]{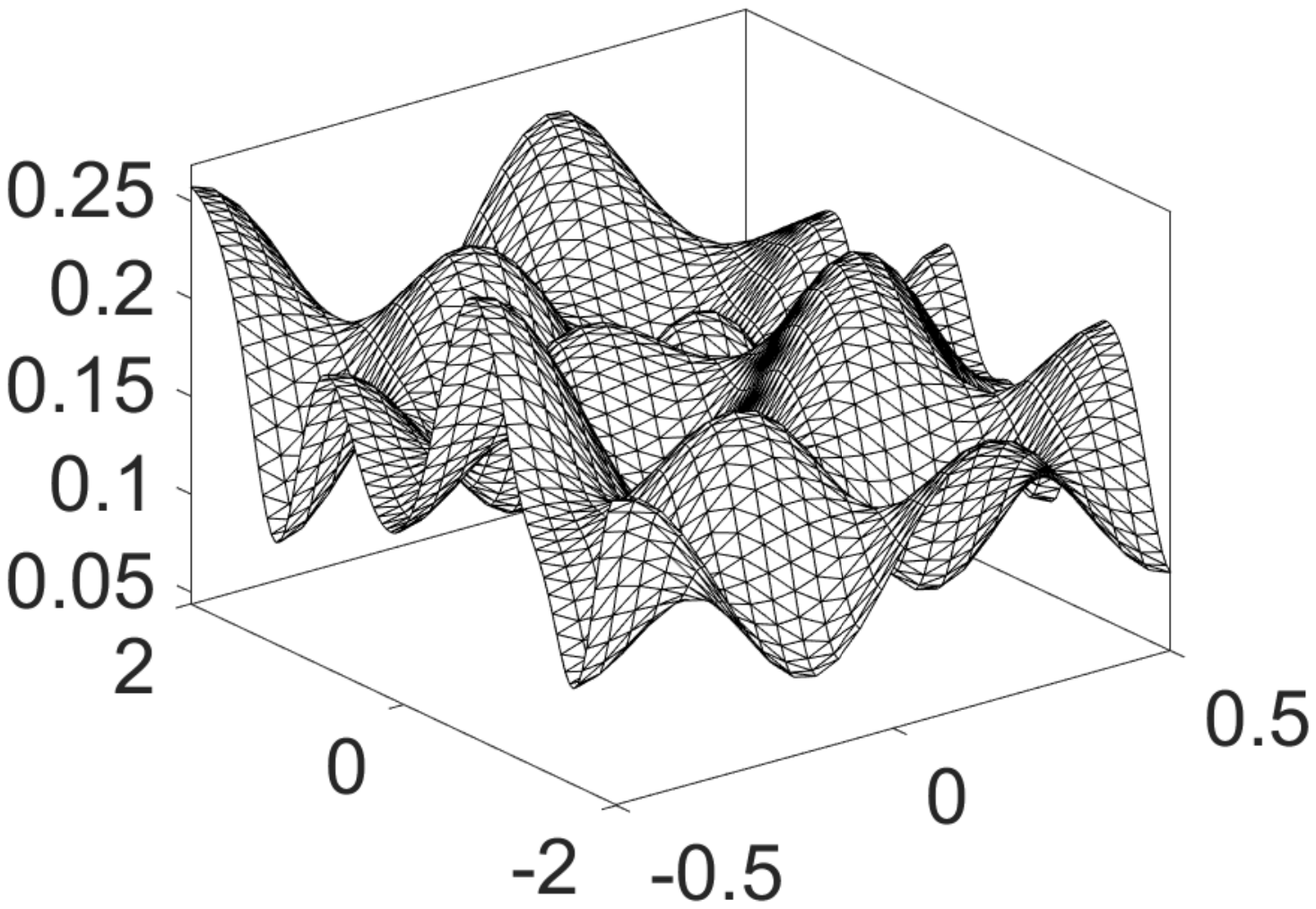}
                \put(80,1){$x$}
         \put(15,5){$y$}
         \end{overpic}}
         \subfloat[\label{sol2Di_2}]{
       \begin{overpic}[width=0.3\textwidth,tics=10,trim=90 240 90 250,clip]{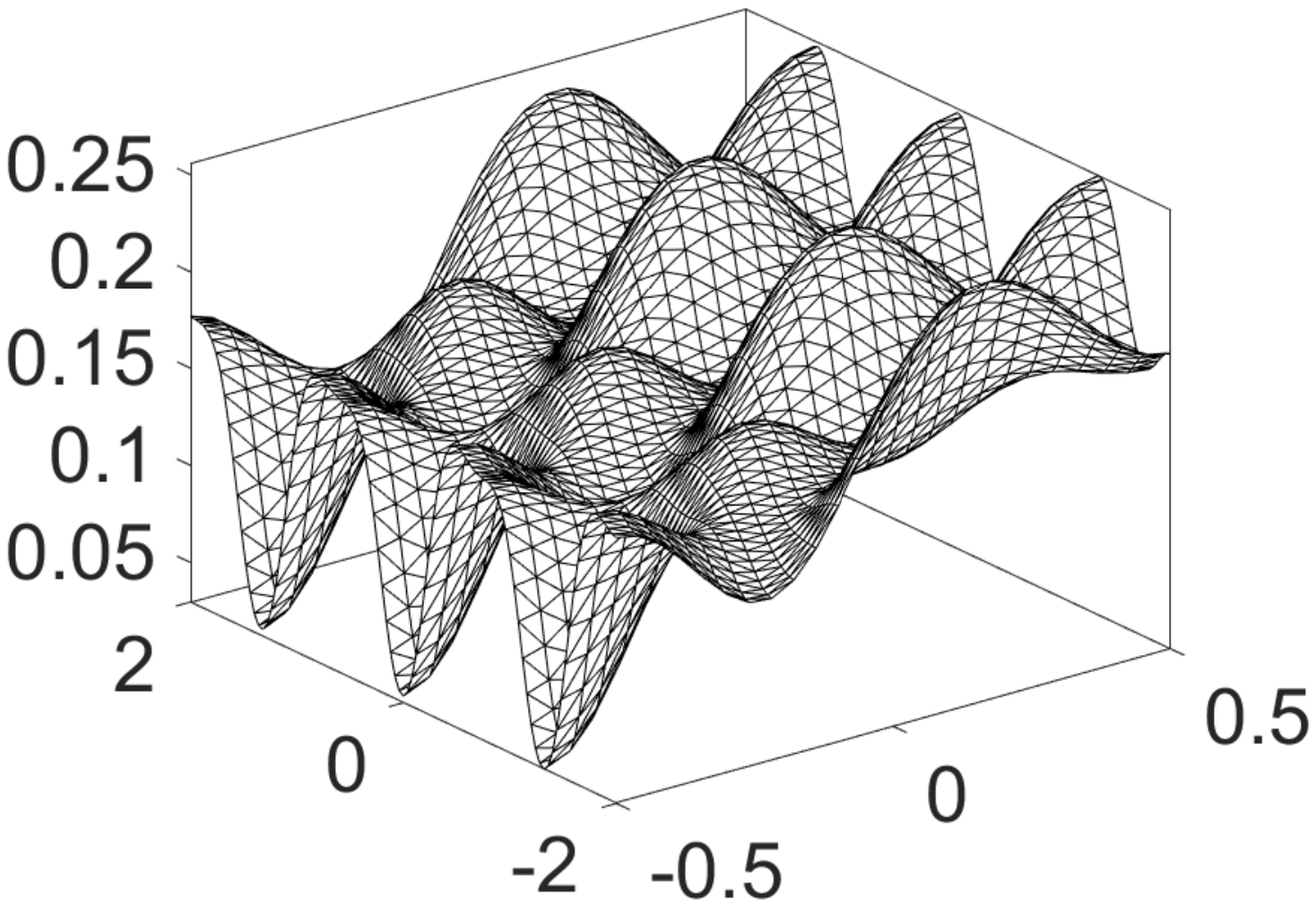}
                \put(80,1){$x$}
         \put(15,5){$y$}
         \end{overpic}}
\caption{Bifurcation diagram and different solution types along the branches. Upper panel: partial bifurcation diagram relative to a rectangular domain $L_x=1,\; L_y=4$. Lower panel: \protect\subref{sol2Da_2}--\protect\subref{sol2Di_2} solutions (species $u$) on different branches.}
\label{2D_DiagBifCross_sol2}
\end{figure}

Finally, Figure \ref{2D_DiagBifCross_sol} shows three different branches, computed far away from the homogeneous one. The first (blue) branch undergoes a further bifurcation and a secondary stable branch arises (marked in green). Then, on this branch a Hopf bifurcation point has been detected. Note that this type of bifurcation is not present in the 1D case in the weak competition regime and yields the existence of time-periodic solutions. Furthermore, the bifurcation diagram in the 2D case seems to be even more intricate compared to the 1D case, many branches are curly and they swirl back and forth. The solution along the branches deforms (see Figures \ref{sol2Da_2}--\ref{sol2Di_2}), and various patterns are evident. Stripes occur in~\ref{sol2Db_2},~\ref{sol2Df_2} and spots in~\ref{sol2Dg_2}, \ref{sol2Di_2}, although they are unstable.

\FloatBarrier


\newpage
\subsection{Convergence of the bifurcation structure}
\label{subsec:2Dconv_d}

\begin{figure}
\subfloat[$\varepsilon=0.05$\label{fast2D_e0p05}]{
       \begin{overpic}[width=0.5\textwidth,tics=10,trim=90 240 100 250,clip]{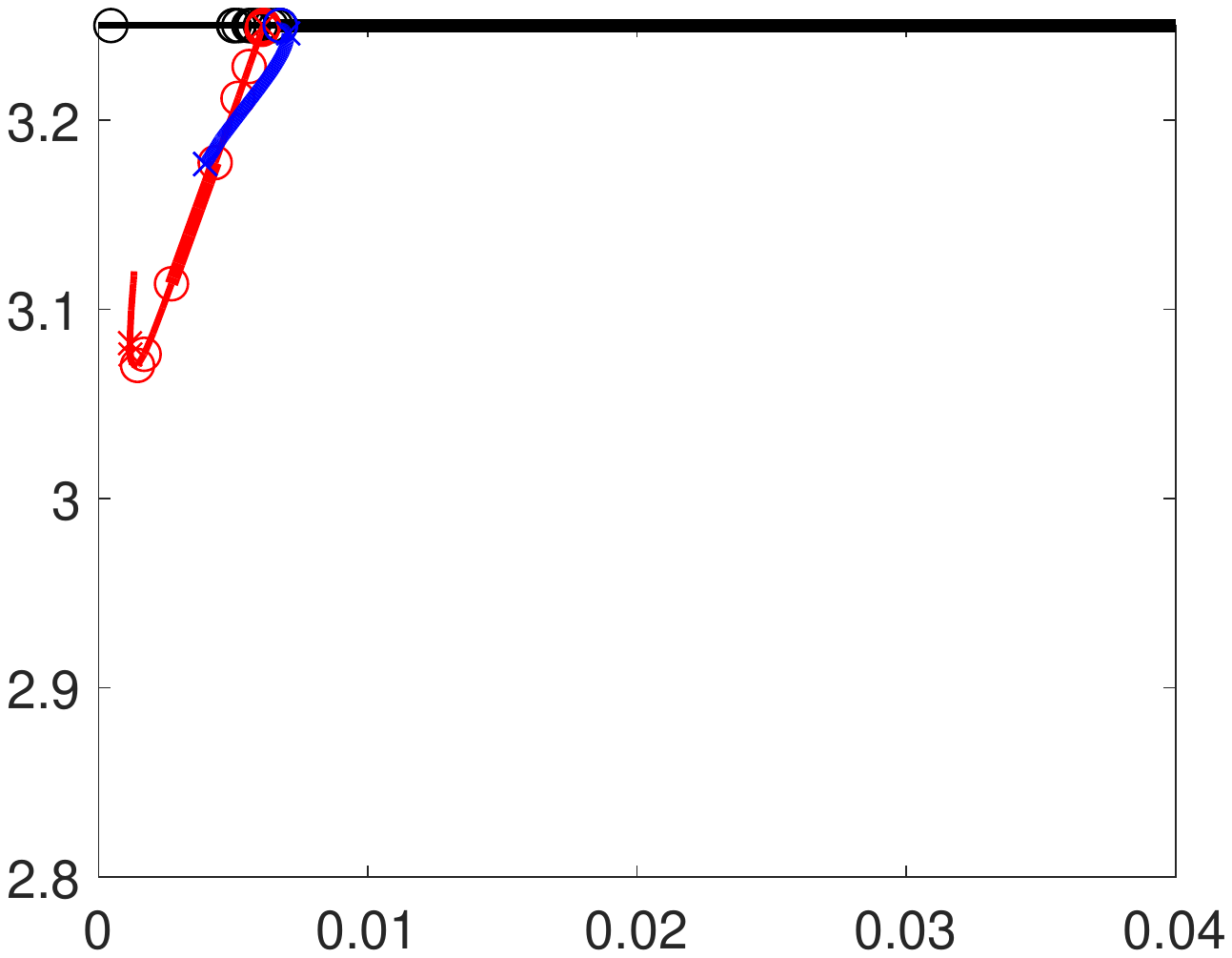}
         \put(0,30){\rotatebox{90}{$||u||_{L^2}$}}
         \put(80,-2){$d$}
       \end{overpic}
}
\subfloat[$\varepsilon=0.01$\label{fast2D_e0p01}]{
       \begin{overpic}[width=0.5\textwidth,tics=10,trim=90 240 100 250,clip]{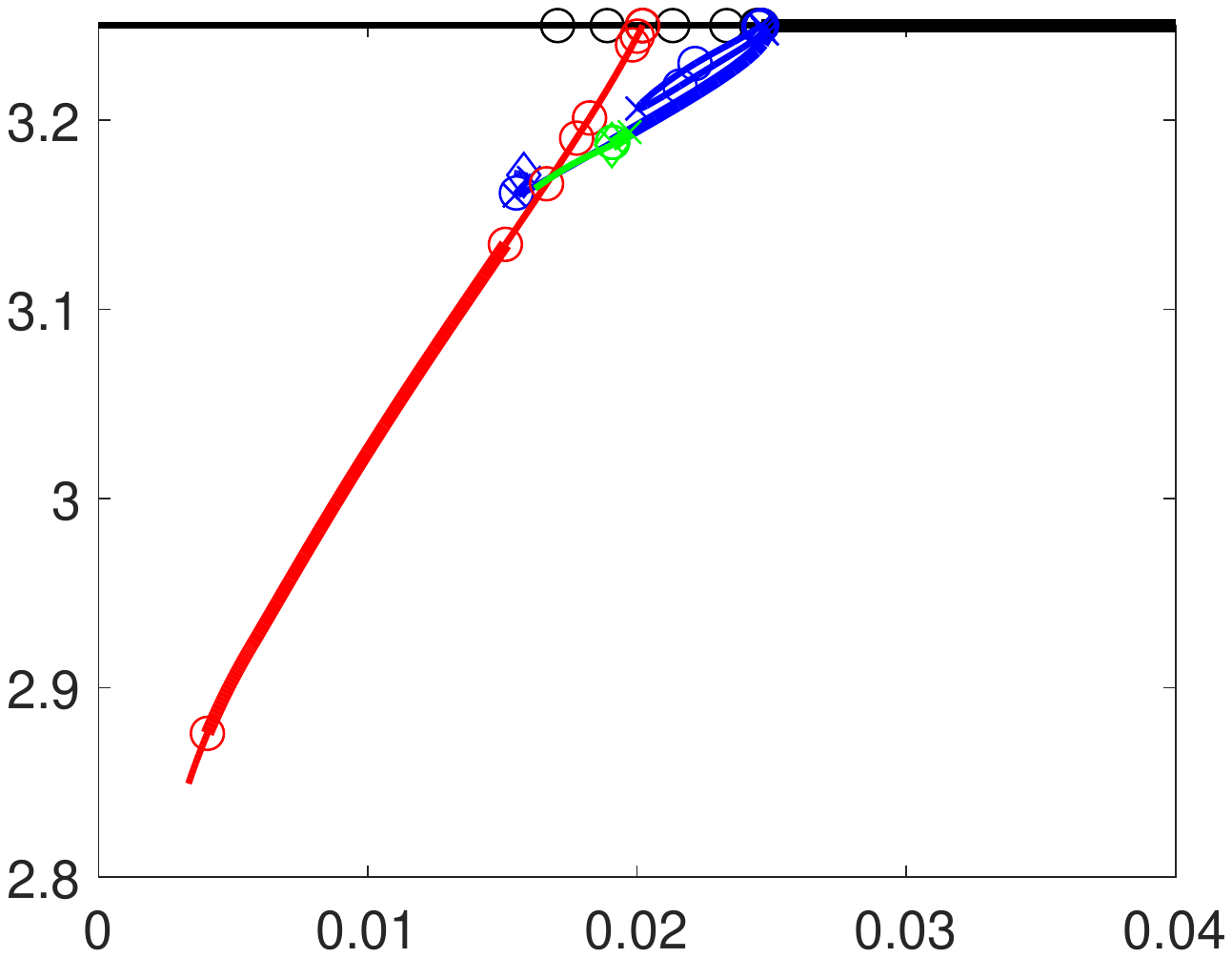}
         \put(0,30){\rotatebox{90}{$||u||_{L^2}$}}
         \put(80,-2){$d$}
       \end{overpic}
}\\
\subfloat[$\varepsilon=0.001$\label{fast2D_e0p001}]{
       \begin{overpic}[width=0.5\textwidth,tics=10,trim=90 240 100 250,clip]{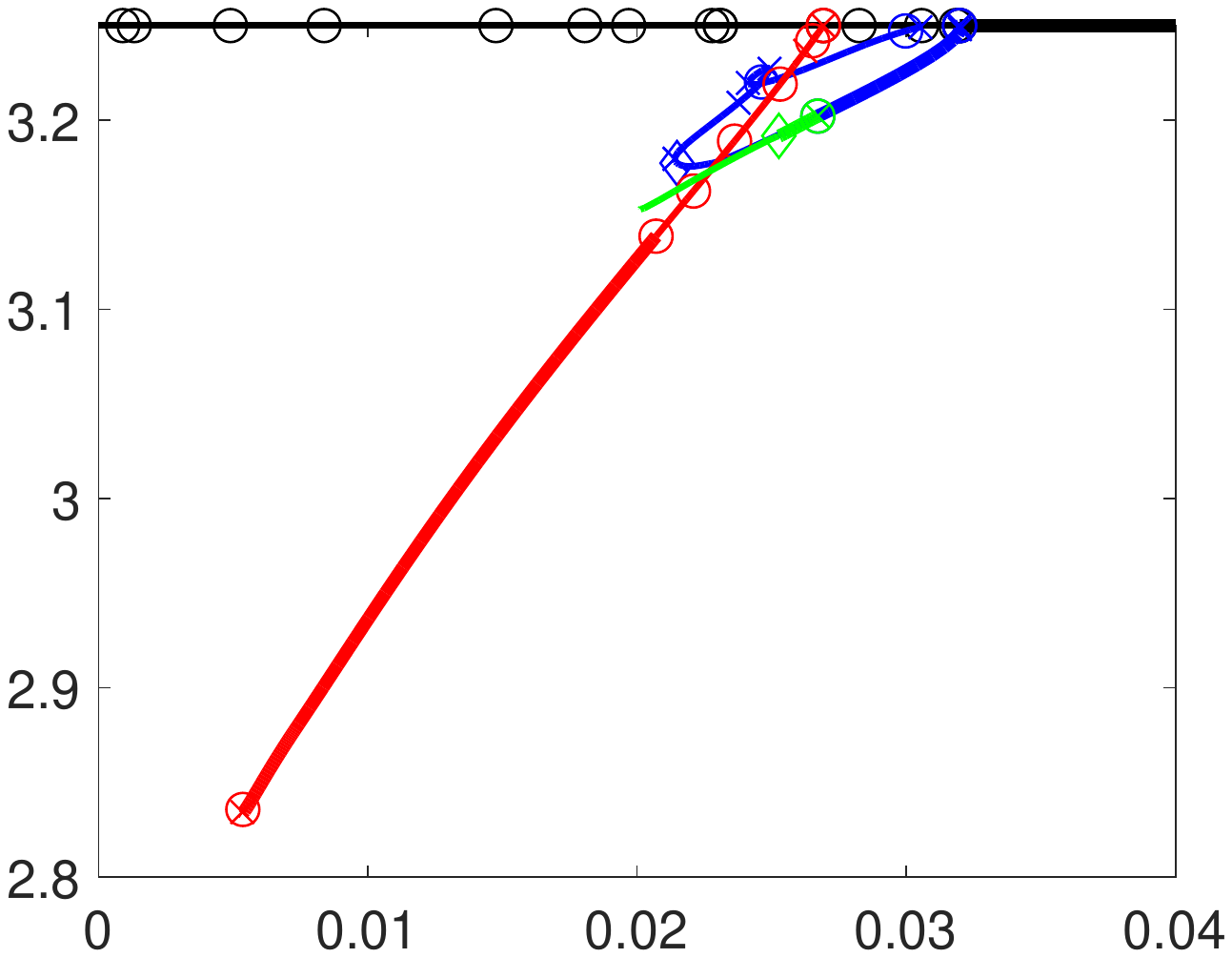}
         \put(0,30){\rotatebox{90}{$||u||_{L^2}$}}
         \put(80,-2){$d$}
       \end{overpic}
}
\subfloat[$\varepsilon=0.0001$\label{fast2D_e0p0001}]{
       \begin{overpic}[width=0.5\textwidth,tics=10,trim=90 240 100 250,clip]{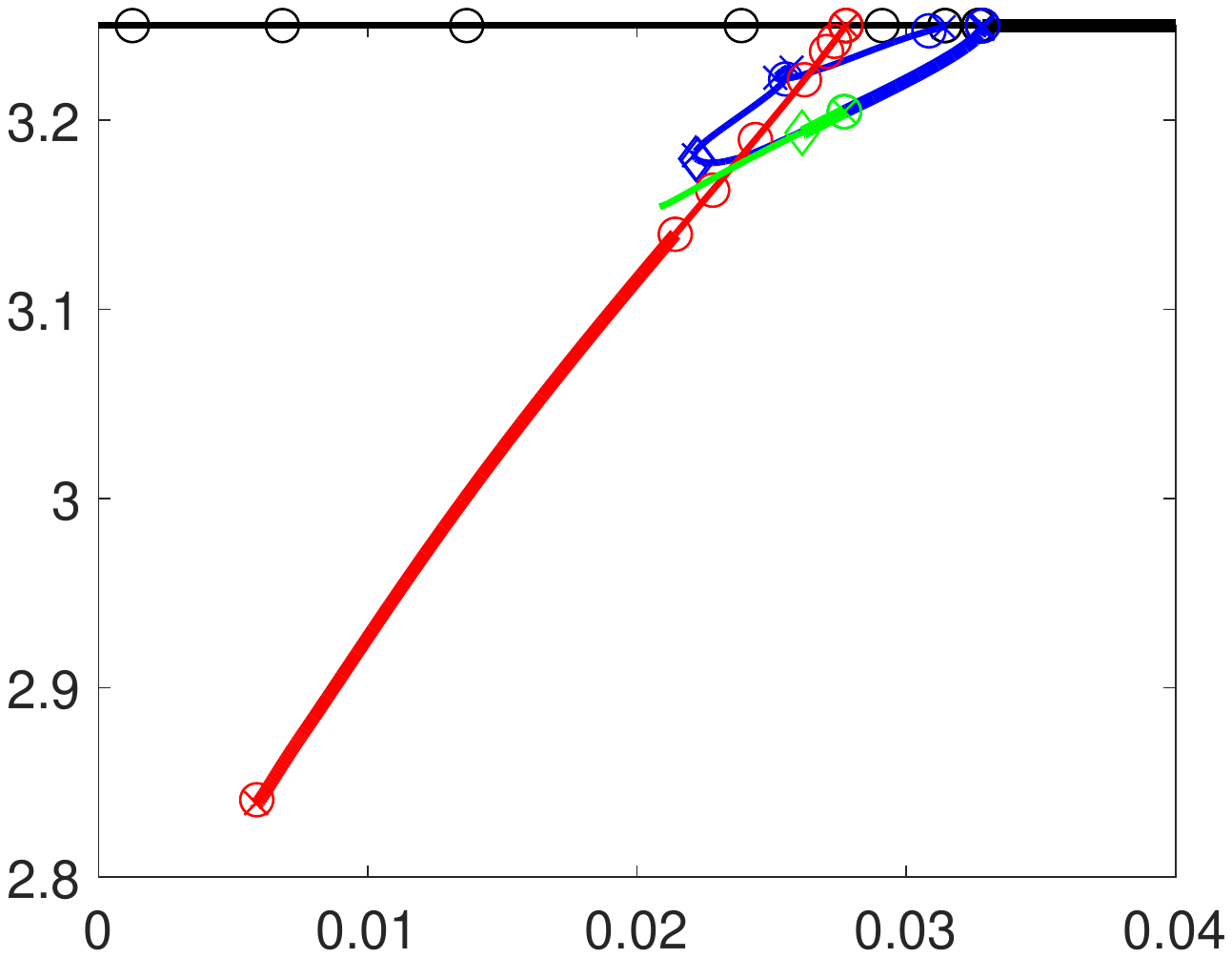}
         \put(0,30){\rotatebox{90}{$||u||_{L^2}$}}
         \put(80,-2){$d$}
       \end{overpic}
}\\
\centering
\subfloat[$\varepsilon=0$ (CD)\label{cross2D}]{
       \begin{overpic}[width=0.5\textwidth,tics=10,trim=90 240 100 250,clip]{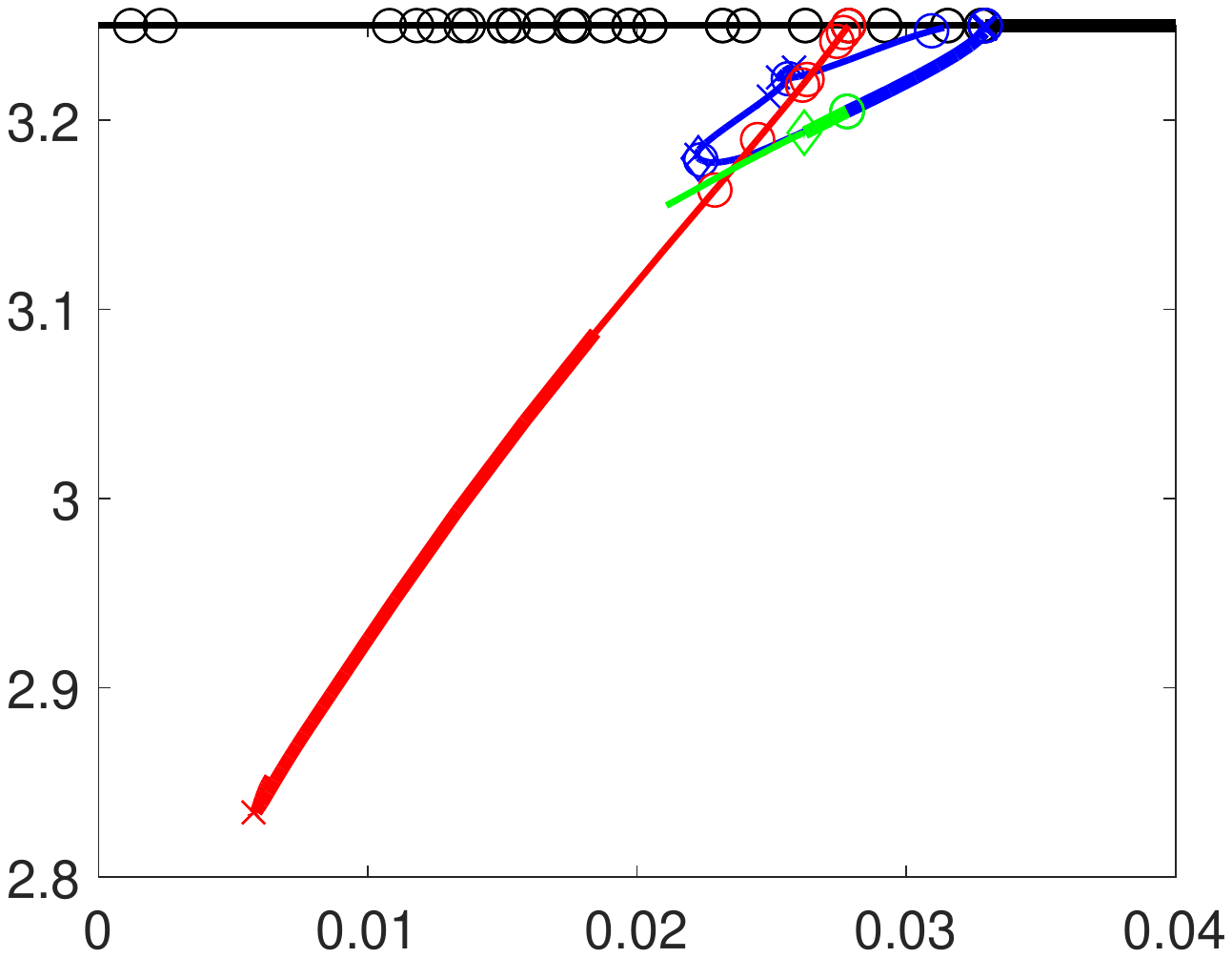}
         \put(0,30){\rotatebox{90}{$||u||_{L^2}$}}
         \put(80,-2){$d$}
       \end{overpic}
}
\caption{Bifurcation diagrams obtained for the 2D domain $[-0.5,0.5]\times[-2,2]$: \protect\subref{fast2D_e0p05}--\protect\subref{fast2D_e0p0001} correspond to the fast-reaction system \eqref{fast} for different and smaller values of $\varepsilon$, while \protect\subref{cross2D} corresponds to the cross-diffusion system \eqref{cross}.}
\label{DiagBifConv2D}
\end{figure}

\begin{figure}
\centering
\subfloat[$u_1$\label{2D_fast_sol_epsi0p0001_blue_u1}]{
       \begin{overpic}[width=0.25\textwidth,tics=10,trim=200 240 200 250,clip]{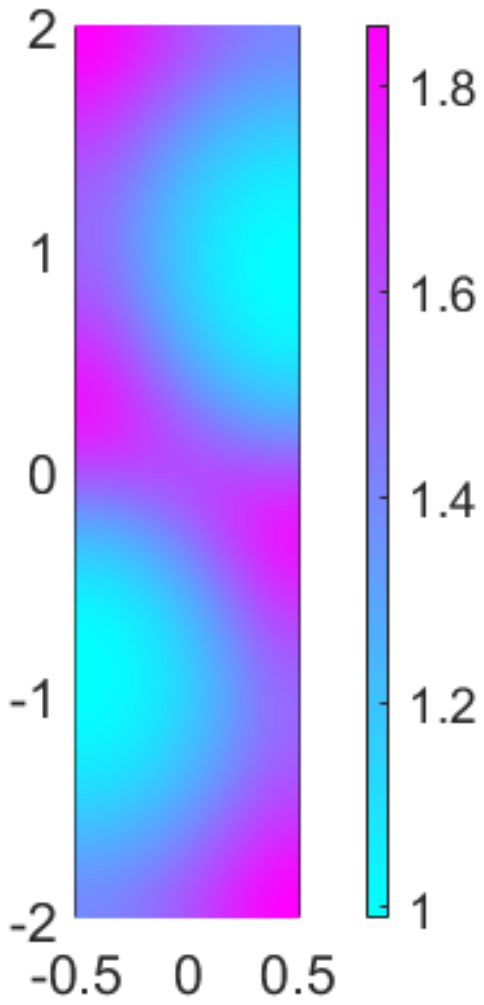}
       \end{overpic}
}
\subfloat[$u_2$\label{2D_fast_sol_epsi0p0001_blue_u2}]{
       \begin{overpic}[width=0.25\textwidth,tics=10,trim=200 240 200 250,clip]{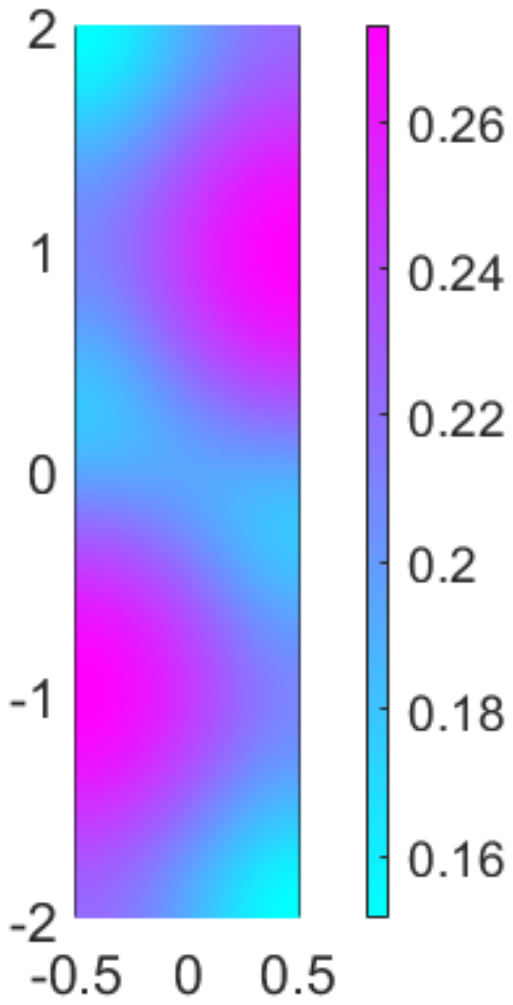}
       \end{overpic}
}
\subfloat[$v$\label{2D_fast_sol_epsi0p0001_blue_v}]{
       \begin{overpic}[width=0.25\textwidth,tics=10,trim=200 240 200 250,clip]{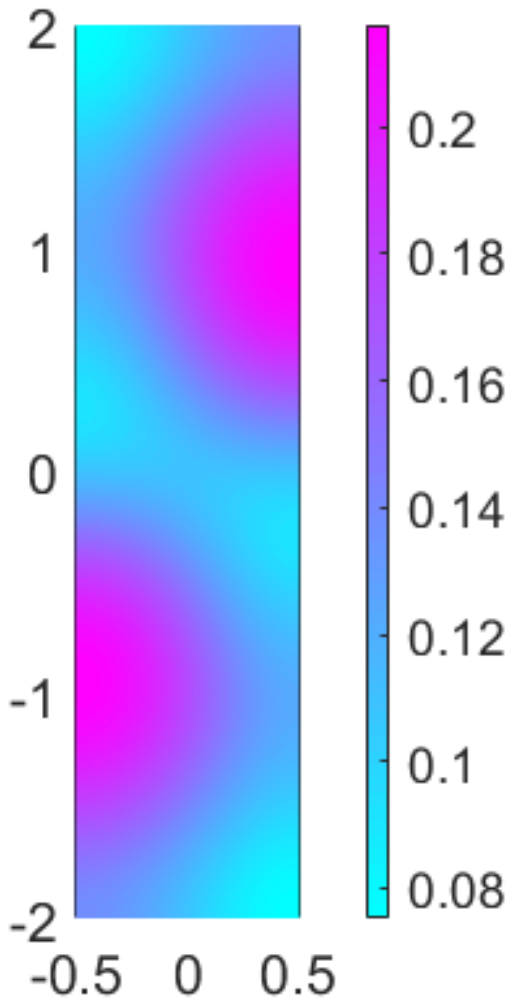}
       \end{overpic}
}\\
\subfloat[$u_1$\label{2D_fast_sol_epsi0p0001_green_u1}]{
       \begin{overpic}[width=0.25\textwidth,tics=10,trim=200 240 200 250,clip]{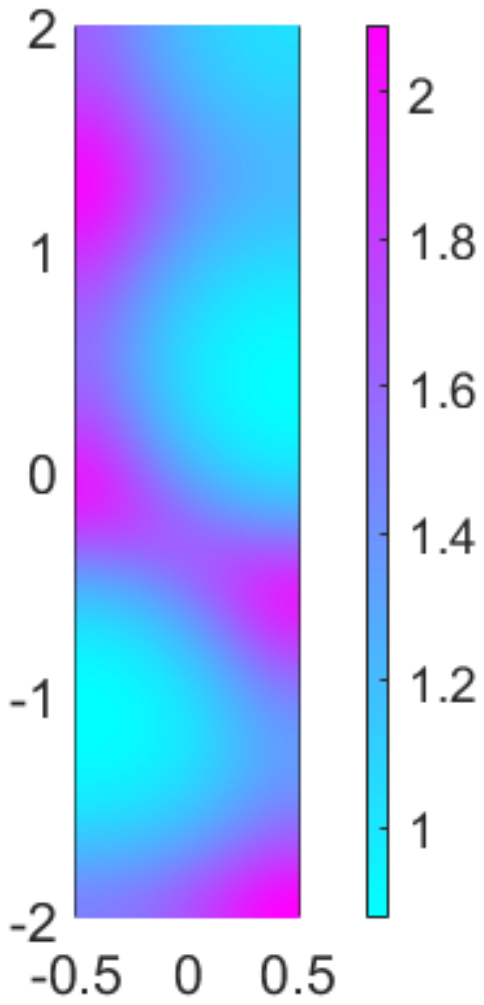}
       \end{overpic}
}
\subfloat[$u_2$\label{2D_fast_sol_epsi0p0001_green_u2}]{
       \begin{overpic}[width=0.25\textwidth,tics=10,trim=200 240 200 250,clip]{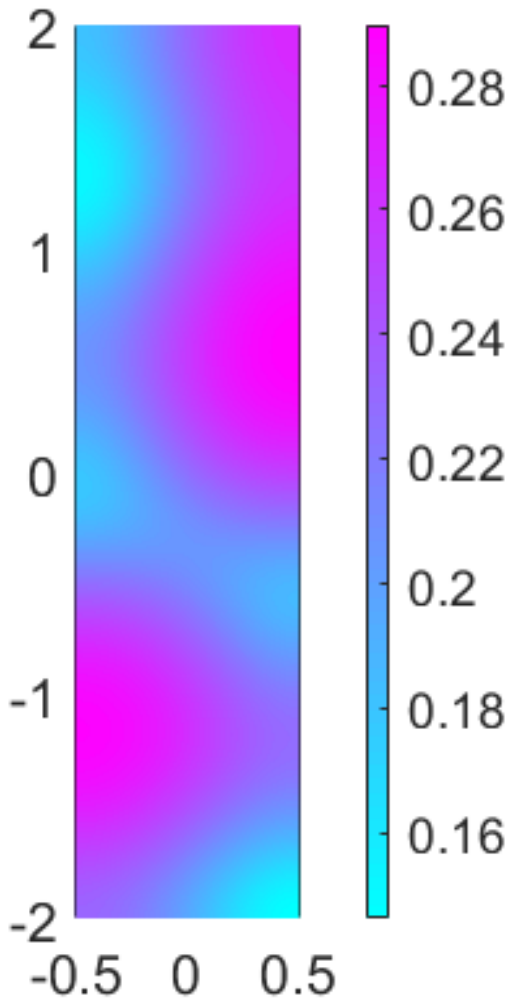}
       \end{overpic}
}
\subfloat[$v$\label{2D_fast_sol_epsi0p0001_green_v}]{
       \begin{overpic}[width=0.25\textwidth,tics=10,trim=200 240 200 250,clip]{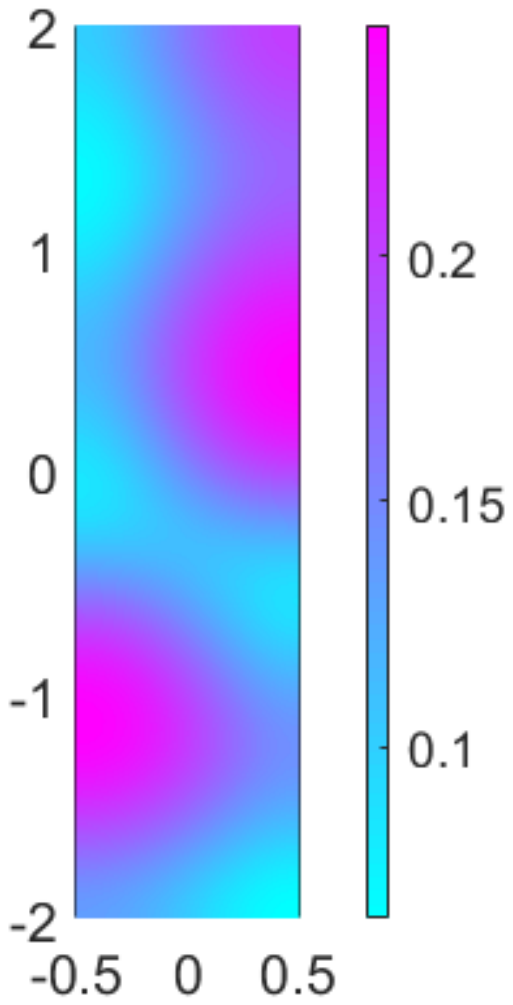}
       \end{overpic}
}\\
\subfloat[$u_1$\label{2D_fast_sol_epsi0p0001_red_u1}]{
       \begin{overpic}[width=0.25\textwidth,tics=10,trim=200 240 200 250,clip]{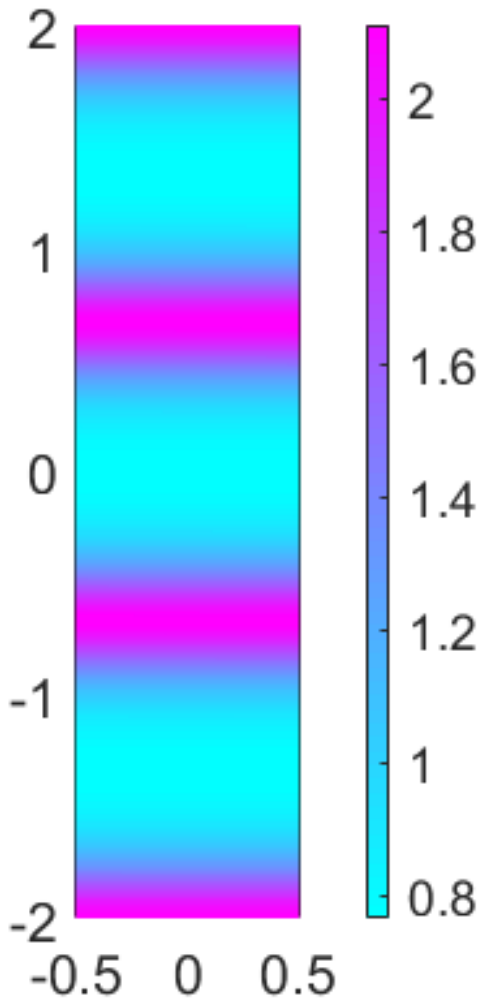}
       \end{overpic}
}
\subfloat[$u_2$\label{2D_fast_sol_epsi0p0001_red_u2}]{
       \begin{overpic}[width=0.25\textwidth,tics=10,trim=200 240 200 250,clip]{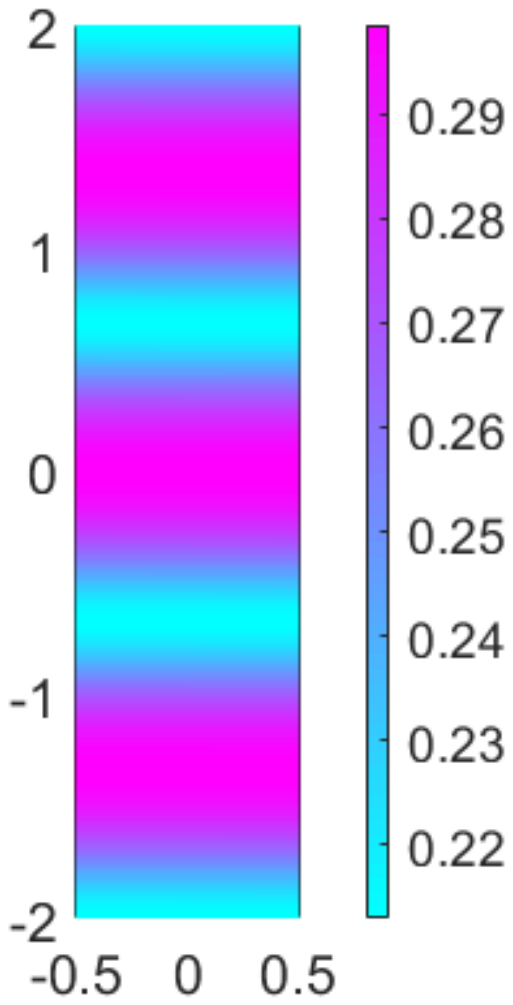}
       \end{overpic}
}
\subfloat[$v$\label{2D_fast_sol_epsi0p0001_red_v}]{
       \begin{overpic}[width=0.25\textwidth,tics=10,trim=200 240 200 250,clip]{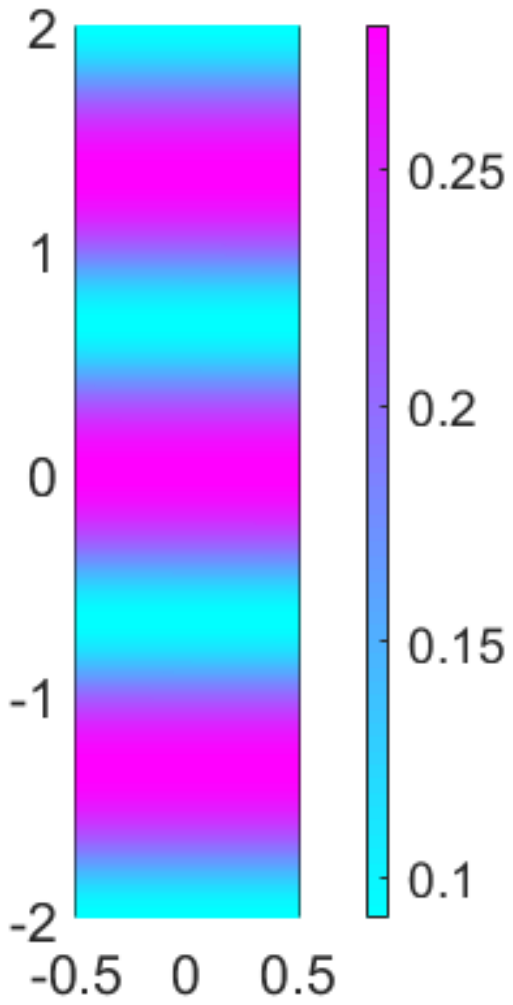}
       \end{overpic}
}
\caption{Three different stable solutions to the fast-reaction system \eqref{fast} (for each one the three species are reported) with $\varepsilon=0.0001$ on the 2D domain $[-0.5,0.5]\times[-2,2]$. With respect to the bifurcation diagram \ref{fast2D_e0p0001}, Figures \protect\subref{2D_fast_sol_epsi0p0001_blue_u1}--\protect\subref{2D_fast_sol_epsi0p0001_blue_v} correspond to the blue branch at $d=0.0297$, Figures \protect\subref{2D_fast_sol_epsi0p0001_green_u1}--\protect\subref{2D_fast_sol_epsi0p0001_green_v} to the green (secondary) branch at $d=0.0262$, while Figures \protect\subref{2D_fast_sol_epsi0p0001_red_u1}--\protect\subref{2D_fast_sol_epsi0p0001_red_v} correspond to the red branch at $d=0.0187$.}
\label{SolFast2D}
\end{figure}

As in the 1D case, we investigate how the bifurcation structure deforms  on a 2D domain with respect to the time-scale separation parameter $\varepsilon$. In Figure \ref{DiagBifConv2D} the bifurcation diagrams for different values of~$\varepsilon$ are reported showing the convergence of the bifurcation structure of the fast-reaction system \eqref{fast} to the cross-diffusion one \eqref{cross}, shown in Figure \ref{2D_fast_sol_epsi0p0001_red_v}. For the sake of clarity of the visualization we only show stable branches. For $\varepsilon=0.05$ the instability region of the homogeneous branch is reduced, and the two primary branches roughly sketch the cross-diffusion one: the first bifurcation point is supercritical, the first branch is stable and the other (red) represents a stable part. However, there is no a secondary bifurcation branch giving rise to a Hopf point. For decreasing values of $\varepsilon$, the cross-diffusion bifurcation structure is well recognizable. 

Finally, in Figure \ref{SolFast2D} we show three different stable stationary solutions of the fast-reaction system with $\varepsilon=0.0001$: they corresponds to points on the blue, green and red branches of Figure \ref{fast2D_e0p0001}, showing different (stable) patterns. It is interesting to observe the distribution of the two different states of species $u$ on the domain. It turns out that the ``excited individuals'' and the competing species $v$ are more abundant in the same region of the domain. Furthermore, the steady state solutions satisfy the algebraic equation derived by a quasi-steady state ansatz from the three-species system \eqref{fast}
$$u_2\left(1-\dfrac v M\right)-u_1 \dfrac v M=0.$$ 

In addition, we have found other stable non-homogeneous solutions of the cross-diffusion system, and the corresponding steady states are stripe patterns.

\clearpage
\section{Conclusion and Outlook}
\label{sec:concl}

In this paper we have investigated the bifurcation structure of the triangular SKT model and of the corresponding fast-reaction system in 1D and 2D domains in the weak competition case via numerical continuation. Despite the fact that part of the bifurcation structure of this system in 1D has been already computed \cite{breden2018existence, breden2013global,iida2006diffusion, izuhara2008reaction}, the key points of our work can be here summarized.
\begin{itemize}[leftmargin=0.3cm]
\item[-] We have adapted the continuation software \texttt{pde2path} to treat cross-diffusion terms. Providing the code in the appendix, this work can be used as a guide to implement such class of problems. 
\item[-] We have reproduced the already computed bifurcation diagrams for the triangular SKT model on a 1D domain with respect to $d$. Even though this is not per se a new result, we have quantitatively checked, how accurate the computation of bifurcation points is. It is worthwhile to note that we have also provided information about the stability of non-homogeneous steady states, which are in agreement with previous results \cite{breden2018existence,iida2006diffusion} and confirms the reliability of our new software setup for cross-diffusion systems.
\item[-] Once the software setup has been established, we have changed the bifurcation parameter to obtain novel structures. We have selected as new bifurcation parameter the growth rate of the first species, which appears in the reaction part. We have shown that the bifurcation structure qualitatively changes as the small parameter tends to zero leading to a ``broken-heart'' bifurcation structure.
\item[-] With respect to both of the considered bifurcation parameters, we have also provided a novel precise quantification of the convergence of the bifurcation points on the homogeneous branch of the fast-reaction system to the cross-diffusion ones.
\item[-] A major new contribution is that we have computed the bifurcation diagram and the non-homogeneous steady states of the triangular SKT model on a 2D domain (rectangular). This case is intricate; the resulting bifurcation diagram is not as clear as in 1D. We have highlighted the main characteristics of the diagrams, and we have presented different types of non-homogeneous steady states, with a focus on stability and pattern formation. We have seen that solutions can exhibit spots or stripes, depending on the parameters, but such solutions turn out to be unstable (at least as far as we have computed the bifurcation branches).
\item[-] We have provided a link between the computed bifurcation diagrams in 1D and 2D domains, and the Turing instability analysis as a tool to fully understand and validate the software output. It is worthwhile to note that the results obtained with the Turing instability analysis only provide insights close to the homogeneous branch. However, the global bifurcation structure has to be numerically computed to achieve a full comprehension of the possible (stationary) outcomes of the system. 
\item[-] Our numerical calculations can now provide also guidelines and conjectures for analytical approaches to cross-diffusion systems, e.g., where in parameter space one can expect entropy structures to behave differently, or where multi-stability and deformation of global branches has to be taken into account. 
\end{itemize} 

Several research directions arise at this point. On the one hand, it is natural to further combine analytical methods and numerical continuation techniques to the study of the full cross-diffusion systems in order to to understand the role of cross-diffusion terms on pattern formation and their influence on the bifurcation diagram \cite{MBCKCS}. On the other hand, continuing our work to numerically investigate of the convergence of the bifurcation structure of a four-equation fast-reaction system leading to the full cross-diffusion SKT system, not only the triangular one, is a straightforward continuation of this work. In this context, analytical proofs of convergence of the solutions of the four species system with two different time-scales to the non-triangular cross-diffusion one is an open problem \cite{desvillettes2018non}. 

Furthermore, the 2D case we have computed has shown that Hopf bifurcations can appear, which can give rise to time-periodic solutions. Their presence could then also be investigated further using the continuation software \texttt{pde2path}, in order to provide a clearer picture on possible long-time asymptotic dynamics of the model. Another direction is to explore other cross-diffusion systems which have been derived by time-scale arguments \cite{conforto2018about,desvillettes2018non,desvillettes20XXnew}, or simply proposed to describe different processes \cite{jungel2017meeting, lacitignola2018cross}. To this end, we have shown that the continuation software is easily applicable to different models and it can be adapted to treat systems that do not exactly belong to the class of problems for which it has been developed. For instance, another direction could be its extension to other non-standard diffusion processes such as fractional reaction--diffusion equations~\cite{NECKCS}.

\bigskip
\textbf{Acknowledgements: }CK have been supported by a Lichtenberg Professorship of the VolkswagenStiftung. CS has received funding from the European Union's Horizon 2020 research and innovation programme under the Marie Sk\l odowska--Curie grant agreement No.~754462. Support by INdAM-GNFM is gratefully acknowledged by CS. CK also acknowledges partial support of the EU within the TiPES project funded the European Unions Horizon 2020 research and innovation programme under grant agreement No.~820970.
\bibliographystyle{plain}
\bibliography{bibliography}
\begin{appendices}
\section{Software setup}
\label{app:setup}

We provide here an explanation how to adapt the \texttt{pde2path} software setup to our problem in the \texttt{OOPDE} setting \cite{prufert2014oopde}, in particular how to treat the cross-diffusion term. Since the purpose of this work is not to give a complete overview of the software for beginner users, we do not explain in detail the basic setup; see \cite{deWitt2017pde2path, uecker2014pde2path,  rademacher2018oopde} for complete guides on the continuation software and for the notation adopted in the following.

In \cite{rademacher2018oopde} the software setup for the quasilinear Allen--Cahn equation is explained, and only recently  this approach was extended to treat a chemotaxis reaction--diffusion system involving a quasi-linear cross-diffusion term of the form $\nabla \cdot (u \nabla v)$ \cite{uecker2019pattern}. However, the previous approach is not directly applicable to the cross-diffusion system \eqref{cross}, since the cross-diffusion term is not written in divergence form. Then we need to rewrite it as
\begin{equation}\label{A:split}
\Lap((d_1+d_{12}v)u)=\nabla \cdot((\underbrace{d_1+d_{12}v}_{c(v)})\nabla u+ \underbrace{d_{12}u}_{\tilde c(u)} \nabla v)
\end{equation}
The steady state problem reads 
$$
0=G(u):=
-\begin{pmatrix} \nabla \cdot (c(v) \nabla u)+\nabla \cdot (\tilde c(u) \nabla v)\\ d_2 \Lap v\end{pmatrix}
-\begin{pmatrix} (r_1-a_1 u-b_1 v)u \\(r_2-b_2 u-a_2 v)v \end{pmatrix},
$$
and on the FEM level it becomes
\begin{equation}\label{A:sm}
G(u)= 
\begin{pmatrix} K_{21}(v)u+K_{12}(u)v\\ Kv\end{pmatrix}
\begin{pmatrix} F_1(u,v)\\ F_2(u,v),\end{pmatrix}
\end{equation}
where $K$ is the standard one-component Neumann-Laplacian (stiffness matrix), $K_{21}(v)u$ and $K_{12}(u)v$ implement $\nabla\cdot(c(v) \nabla u)$ and $\nabla \cdot (\tilde c(u) \nabla v)$ respectively, while $F_1,\,F_2$ belong to the reaction part. In detail,
\begin{equation}\label{A:KK}
(K_{21}(v))_{ij}=\int_\Omega c(v) \nabla\phi_i \cdot \nabla \phi_j~\txtd x,\qquad (K_{12}(u))_{ij}=\int_\Omega \tilde c(u) \nabla\phi_i \cdot \nabla \phi_j~\txtd x
\end{equation}
depend on $v$ and on $u$. Hence, they have to be computed at each step. In the \texttt{OOPDE} setting, we can employ the routine \texttt{assema} to compute those matrices, but this needs $c(v)$ and $\tilde c(u)$ on each element center, which is obtained interpolating $u$ and $v$ from the nodes to the element centers, as it can be seen in Listings \ref{lst:sG} and \ref{lst:oosetfemops}, which show the main files implementing the triangular cross-diffusion system \eqref{cross} in \texttt{pde2path}.

\lstset{
numbers=left,
firstnumber=1,
stepnumber=5,
numbersep=5pt,
numberstyle=\small\color{black},
basicstyle=\ttfamily,
keywordstyle=\color{black},
commentstyle=\color{green},
stringstyle=\color{black},
tabsize=2,
backgroundcolor=\color{gray!30!white},
captionpos=b}
\begin{lstlisting}[
caption={\texttt{sG.m}. In particular in lines 18, 19 we define the functions $c,\,\tilde{c}$ defined in \eqref{A:split}, and in line 20, 21 we use them to build the matrices $K_{21}$ and $K_{12}$ appearing in \eqref{A:KK} using the \texttt{OOPDE} routine \texttt{assema}. The system matrix is then assembled in line 26, following \eqref{A:sm}.},
label={lst:sG}]
function r=sG(p,u) 
% compute pde-part of residual
	u1=u(1:p.np); % extract the first component
	u2=u(p.np+1:2*p.np); % extract the second component
	par=u(p.nu+1:end); % extract parameters
	d=par(1); d12=par(2);
	r1=par(3); a1=par(5); b1=par(7); 
	r2=par(4); a2=par(6); b2=par(8);

	f1=(r1-a1*u1-b1*u2).*u1;
	f2=(r2-b2*u1-a2*u2).*u2;

	gr=p.pdeo.grid;
	% interpolate to element centers
	u1t=gr.point2Center(u1);
	u2t=gr.point2Center(u2);

	c=d+d12*u2t;
	cc=d12*u1t;
	[Kc,~,F1]=p.pdeo.fem.assema(gr,c,0,f1);
	[Kcc,~,~]=p.pdeo.fem.assema(gr,cc,0,f1);
	[K,~,F2]=p.pdeo.fem.assema(gr,1,0,f2);
	N=sparse(p.pdeo.grid.nPoints,p.pdeo.grid.nPoints);
	p.mat.K=[[Kc Kcc];[N d*K]];
	F=[F1;F2];
	r=p.mat.K*[u1;u2]-F;
end
\end{lstlisting}

\begin{lstlisting}[caption={\texttt{oosetfemops.m}. The stiffness matrix needs to be build at each step, while the mass matrix $M$ can be assembled here.},label={lst:oosetfemops}]
function p=oosetfemops(p)
    gr=p.pdeo.grid;
    [~,M,~]=p.pdeo.fem.assema(gr,0,1,1); % assemble 'scalar' M
    p.mat.M=kron([[1,0];[0,1]],M); %build 2-component system M
end
\end{lstlisting}

%
%

\section{Turing instability analysis}
\label{app:TI}

In this section the detailed Turing instability analysis of system \eqref{cross} can be found. Although it is a simpler case than the one treated in \cite{MBCKCS}, this computation allows us to 
identify the bifurcation points along the homogeneous branch. These values can be compared to the ones numerically obtained using \texttt{pde2path}, so we present the analysis for reference here.


We recall that in the weak competition case ($a_1a_2-b_1b_2>0$) there exists a homogeneous equilibrium $(u_*,v_*)$, which is stable for the reaction under the additional condition on the growth rates $b_1/a_2<r_1/r_2<a_1/b_2$. Then, the Jacobian matrix of the reaction part and the linearization of the diffusion part of the cross-diffusion system \eqref{cross}, evaluated at the equilibrium $(u_*,v_*)$ are
$$J_*= \begin{pmatrix}-a_1u_* & -b_1u_*\\-b_2v_*& -a_2v_*\end{pmatrix}, \quad
J_\Delta^*=\begin{pmatrix}d+d_{12}v_*& d_{12}u_*\\0 & d\end{pmatrix},$$
where $\tr J_*<0, \, \det J_*>0$.
Hence, the characteristic matrix is
$$M_k^*=J_*-J_\Delta^*\lambda_k=
\begin{pmatrix}
-a_1u_*-(d+d_{12}v_*)\lambda_k   &  -b_1u_*-d_{12}u_*\lambda_k\\
-b_2v_*          & -a_2v_*-d  \lambda_k
\end{pmatrix},$$
where $\lambda_k$ denotes an eigenvalue of the Laplacian on the domain.
The trace of the characteristic matrix remains negative, while the determinant can be written as a second order polynomial in $\lambda_k$ (as usual in Turing instability analysis)
\begin{equation}\label{detM*ks}
\det M_k^*=d\left(d+d_{12}v_*\right)\lambda_k^2
-\left(d\tr J_*+d_{12}\alpha\right)\lambda_k
+\det J_*,
\end{equation}
where $\alpha:=(2b_2u_*-r_2)v_*$. In order to have Turing instability, the determinant of the characteristic matrix has to be negative for some $\lambda_k$. This is possible when $\alpha>0$. The sign of the quantity $\alpha$ is not fixed but it depends on the parameter values $r_i,\,a_i,\, b_i,\, (i=1,2)$; we have a positive $\alpha$ when
$$ \dfrac{r_1}{r_2}>\dfrac 1 2 \left(\dfrac{b_1}{a_2}+\dfrac{a_1}{b_2}\right).$$
\begin{remark}
Note that the parameter set used in \cite{breden2018existence, breden2013global,iida2006diffusion, izuhara2008reaction} and reported in Table \ref{tab:param} gives $\alpha>0$.
\end{remark}
If we rewrite \eqref{detM*ks} as a second-order polynomial in $d$
\begin{equation}\label{detM*d}
\det M_k^*=\lambda_k^2 d^2 +
(d_{12}v_*\lambda_k^2-\tr J_*)d
-d_{12}\alpha \lambda_k +\det J_*,
\end{equation}
it is possible to locate the bifurcation points solving $\det M_k^*=0$, we obtain
\begin{equation}
d_B=d_B(d_{12},\lambda_k)=\dfrac{-(d_{12}v_*\lambda_k-\tr J_*)+\sqrt{(d_{12}v_*\lambda_k-\tr J_*)^2-4(\det J_*-d_{12}\alpha \lambda_k)}}{2\lambda_k},
\label{d_B}
\end{equation}
depending on the parameters and on the eigenvalue of the Laplacian. With homogeneous Neumann boundary conditions and choosing a 1D domain with length $L_x$, the eigenvalues are
$$\lambda_k=\left( \dfrac{\pi n}{L_x}\right)^2,\quad n\geq0,$$
while on a rectangular 2D domain 
$$\lambda_k=\lambda_{n,m}=\left( \dfrac{\pi n}{L_x}\right)^2+\left( \dfrac{\pi m}{L_y}\right)^2,\quad n,m\geq0.$$
In Table \ref{table:dB1D} we report the bifurcation values $d_B(\lambda_k)$ corresponding to the eigenvalue $\lambda_k$ of the Laplacian for the 1D domain $(0,1)$ and obtained by formula \eqref{d_B}, for different values of $n$, compared with the numerical values estimated by the software \texttt{pde2path}. The obtained values are in good agreement. The scenario is more complicated in a rectangular 2D domain, reported in Table \ref{table:dB1D} for different values of $n$ and $m$. Also in this case we compare the values obtained by formula \eqref{d_B} with the numerical values. The software is not able to detect all the bifurcation values, but it fails when they are too close to each other or with multiplicity more than one. Note that the values are ordered with respect to the bifurcation values $d_B$, which does not translate into an order of the indices $n,\,m$: for instance the couple $(0,3)$ corresponds to a smaller bifurcation value than $(0,7)$. The reason is evident in Figure \ref{table:dB2D}, due to a non-monotonicity of the function $d_B(\lambda_k)$ w.r.t.~$\lambda_k$.

\begin{table}
\centering
\begin{tabular}{ccc}
\toprule
$k$ & $d_B(\lambda_k)$ & \texttt{DB}\\
\midrule
0 & - & -\\
1 & 0.032788 & 0.03279\\
2 & 0.02049 & 0.02046\\
3 & 0.01138 & 0.01133\\
4 & 0.00699 & 0.00693\\
5 & 0.00467 & 0.00460 \\
\bottomrule
\end{tabular}
\caption{Bifurcation values corresponding  to the eigenvalue $\lambda_k$ of the Laplacian for the 1D domain $(0,1)$: values $d_B(\lambda_k)$ are obtained from formula \eqref{d_B}, while values $\texttt{BP}$ are estimated by \texttt{pde2path}, using an uniform mesh with 26 grid points, maximal and minimal step size \texttt{p.nc.dsmax=1e-4} and \texttt{p.nc.dsmin=1e-7}.}
\label{table:dB1D}
\end{table}

\begin{table}
\centering
\begin{tabular}{cccc}
\toprule
$n$&$m$ & $d_B(\lambda_{n,m})$& \texttt{BP}\\
\midrule
     1   &  1  &     0.0329340 &  0.0329346\\
     0   &  4  &     0.032788 &\\
     1   &  0  &     0.032788 &\\
     1   &  2  &     0.032783 & 0.032789\\
     0   &  5  &     0.031545 &\\
     1   &  3  &     0.031545 & 0.031541\\
     1   &  4  &     0.02921 & 0.02920\\
     0   &  6  &     0.027865 & 0.027856\\           
     1   &  5  &     0.02627 & 0.02625\\        
     0   &  7  &     0.02397 & 0.02396\\             
     0   &  3  &     0.0236 & 0.0239\\
     0   &  8  &     0.02049 &\\
     2   &  0  &     0.02049 & 0.02048\\
     2   &  1  &     0.02029 & 0.02028\\
\bottomrule
\end{tabular}
\caption{Bifurcation values corresponding  to the eigenvalue $\lambda_{n,m}$ of the Laplacian for the rectangular 2D domain $[-0.5,0.5]\times[-2,2]$: values $d_B(\lambda_{n,m})$ are obtained from formula \eqref{d_B}, while values $\texttt{BP}$ are estimated by \texttt{pde2path}, using a uniform mesh with 26 grid points on the $x-$edge, maximal and minimal step size \texttt{p.nc.dsmax=1e-4} and \texttt{p.nc.dsmin=1e-7}. With these settings, the software is not able to localize all the expected bifurcation values.}
\label{table:dB2D}
\end{table}

\begin{figure}[!ht]
\centering
 \begin{overpic}[width=0.7\textwidth,tics=10,trim=90 240 100 250,clip]{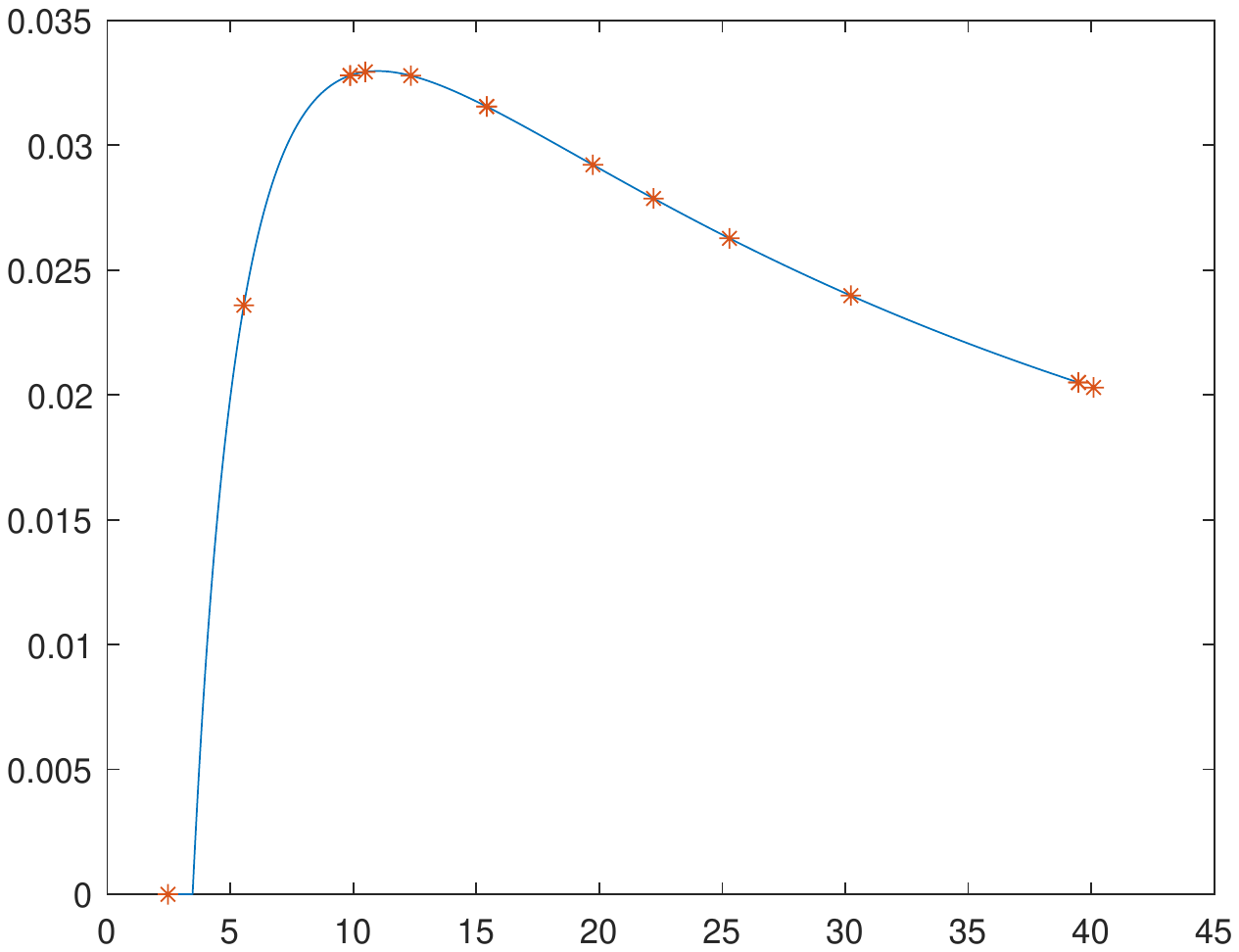}
         \put(0,30){\rotatebox{90}{$d_B(\lambda)$}}
         \put(85,0){$\lambda$}
%
\put(15,10){\small{$\lambda_{0,2}$}}
\put(17,49){\small{$\lambda_{0,3}$}}

\put(30,66){\small{$\lambda_{1,1}$}}
\put(28,52){\small{$\lambda_{0,4},k_{1,0}$}}
\put(32,58){\small{$\lambda_{1,2}$}}
\put(40,64){\small{$\lambda_{0,5}, k_{1,3}$}}
\put(48,60){\small{$\lambda_{1,4}$}}
\put(54,56){\small{$\lambda_{0,6}$}}
\put(59,53){\small{$\lambda_{1,5}$}}
\put(66,50){\small{$\lambda_{0,7}$}}
\put(80,45){\small{$\lambda_{0,8},\lambda_{2,0}$}}
\put(80,36){\small{$\lambda_{2,1}$}}

\put(31,55){\color{gray!30!white}\vector(0,2){7}}
\put(35.5,60){\color{gray!30!white}\vector(0,2){3}}
\put(83,39){\color{gray!30!white}\vector(0,2){3}}
       \end{overpic}
       \caption{Bifurcation value $d_B$ as function of $\lambda$ (solid line) obtained from formula \eqref{d_B}. Markers indicate the bifurcation values corresponding  to the eigenvalue $\lambda_{n,m}$ of the Laplacian for the rectangular 2D domain $[-0.5,0.5]\times[-2,2]$.}
\label{fig:dB}
\end{figure}
\end{appendices}
\end{document}